\documentclass[11pt]{article}
\usepackage{placeins,graphicx,diagbox }

\usepackage{graphicx,bm}
\usepackage{amsfonts,amsmath,fullpage,bbm}
\usepackage{amssymb,multirow,verbatim}
\usepackage{acronym,wrapfig,plain,mathrsfs,enumerate,relsize,color}
\newtheorem{algorithm}{Algorithm}
\usepackage{algorithm}
\newtheorem{assumption}{Assumption}
\usepackage{subfig}
\usepackage{algorithmic}
\usepackage{pifont}

\newcommand{\us}[1]{{\color{black}#1}}

\usepackage{wrapfig} 
\usepackage{times} 
\usepackage{amsmath,mathtools,float}
\usepackage{tikz,pgfplots}
\usepackage{booktabs,caption}
\usepackage[flushleft]{threeparttable}
\usetikzlibrary{arrows,shapes,snakes,shadows,positioning,automata,patterns}
\usetikzlibrary{trees,decorations.pathmorphing,decorations.markings}
\usepackage{amsmath}
\usepackage{multirow}
\usepackage{amssymb,latexsym}  
\usepackage{dsfont} 

\def\tkappa{\tilde{\kappa}}
\usepackage{amsmath,algorithm} 
\usepackage{amssymb}  
\usepackage{color}
\usepackage{cite}

\newcommand{\wh}{\widehat}

\newcommand{\pmat}[1]{\begin{pmatrix} #1 \end{pmatrix}}

\def\Rscr{\mathcal{R}}

\def\tkappa{\tilde{\kappa}}

\newcommand{\blue}[1]{{\color{black}#1}}

\newcommand{\epsilonbar}{{\bar \epsilon}}

\def\Real{\mathbb{R}}

\def\argmin{\mathop{\rm argmin}}

\newtheorem{theorem}{Theorem}
\newtheorem{lemma}{Lemma}

\newtheorem{proposition}{Proposition}
\newtheorem{cor}{Corollary}
\newtheorem{remark}{Remark}

\def\us#1{{{\color{black}#1}}}
\usepackage{verbatim}
\usepackage{ulem}
\newcommand{\vvs}[1]{{\color{black}#1}}

\def\be{\begin{enumerate}}
\def\ee{\end{enumerate}}

\def\argmin{\mathop{\rm argmin}}

 \newcommand{\remove}[1]{}

\def\Real{\mathbb{R}}

\def\argmin{\mathop{\rm argmin}}

\usepackage{tikz}
\usetikzlibrary{decorations.pathreplacing,calc}

\date{}

\begin{document}
\title{Distributed Variable Sample-Size Gradient-response and Best-response Schemes    for   Stochastic Nash  Equilibrium Problems   over Graphs}

\author{  Jinlong  Lei  and  Uday  V. Shanbhag \thanks{ Email:  leijinlong@tongji.edu.cn (Jinlong Lei), udaybag@psu.edu
(Uday  V. Shanbhag). The work  has been partly supported  by NSF grant 1538605
and 1246887 (CAREER).} \thanks{  Lei is  the Department of Control Science and Engineering, Tongji University; She is also with the Shanghai Institute of Intelligent Science and Technology, Tongji University, Shanghai 200092, China. Shanbhag is with
the Department of Industrial and Manufacturing Engineering, Pennsylvania State
University, University Park, PA 16802, USA}}

\maketitle

\begin{abstract} This paper considers a stochastic  NEP  in which each player $i$
minimizes a composite objective $f_i(x) + r_i(x_i)$, where  $f_i$ is an
expectation-valued smooth function and  $r_i$ is a nonsmooth convex
function with an efficient prox-evaluation. In this context, we make the
following contributions. (I)  Under suitable monotonicity assumptions on the
concatenated  gradient map of $f_i$, we derive  {\bf optimal}  rate statements
and oracle complexity bounds for the proposed variable sample-size  proximal
stochastic gradient-response ({\bf VS-PGR}) scheme when the sample-size
increases at a geometric rate. If the sample-size increases at a polynomial
rate  with  degree $v > 0$, the mean-squared error of
the iterates decays at a corresponding polynomial rate while the iteration and
oracle complexities  to  obtain an  $\epsilon$-Nash equilibrium (NE) are
$\mathcal{O}(1/\epsilon^{1/v})$ and $\mathcal{O}(1/\epsilon^{1+1/v})$,
respectively. (II) We then overlay (VS-PGR) with a consensus phase with a view
towards developing distributed protocols for aggregative stochastic NEPs.
In the resulting ({\bf d-VS-PGR}) scheme, when the sample-size  at each iteration grows at a geometric   rate
while the communication rounds per iteration grow at the rate of $ k+1 $, computing an  $\epsilon$-NE requires similar
iteration and oracle complexities to ({\bf VS-PGR}) with a communication
complexity of $\mathcal{O}(\ln^2(1/\epsilon))$. (III) Under a  suitable
contractive property  associated with the proximal best-response (BR)  map,
we design a variable sample-size proximal BR ({\bf VS-PBR}) scheme, where
each  player solves   a sample-average  BR problem. When the
sample-size  increases at a suitable geometric  rate, the
resulting iterates converge at a  geometric  rate while the iteration and
oracle complexity are 	respectively $\mathcal{O}(\ln(1/\epsilon))$  and
$\mathcal{O}(1/\epsilon)$; If the sample-size increases  at a polynomial
rate with degree $v$, the mean-squared error  decays at a
corresponding polynomial rate while the iteration and oracle complexities are
$\mathcal{O}(1/\epsilon^{1/v})$ and $\mathcal{O}(1/\epsilon^{1+1/v})$,
respectively. (IV) Akin to (II), the distributed variant ({\bf d-VS-PBR})
achieves similar iteration and oracle complexities to the centralized   ({\bf
VS-PBR}) with a communication complexity of $\mathcal{O}(\ln^2(1/\epsilon))$
 when the communication rounds per iteration increase at the rate of
$ k+1 $.  
\end{abstract}

\section{Introduction}
Noncooperative games~\cite{fudenberg91game,basar99dynamic}
 consider  the resolution of conflicts  amongst  selfish players, each of
which tries to optimize its payoff,   given its rival strategies.  The Nash equilibrium represents an important solution concept for noncooperative games  \cite{nash50equilibrium}  and has  seen wide applicability in a  breadth
of engineered systems, such as power grids, communication networks,  and
transportation networks (see e.g. \cite{yin09nash2,basar2007control,pan2007games,pang2010design}).
Recently, there has been an interest in ``designing'' games for
distributed control~\cite{li2013designing,arslan2007autonomous}; consequently,
in networked regimes, the role of distributed protocols for computing
equilibria over graphs is of increasing relevance.

 We consider a Nash equilibrium problem (NEP), denoted by $\mathscr{P}$, with a set of  $n $ players indexed by
$i$,  where $i \in\mathcal{N} \triangleq \{1, \cdots,n\}.$ For each $i \in \cal{N}$, the $i$th player  is characterized by  a strategy  $x_i \in \mathbb{R}^{d_i}$  and  an objective  \blue{$ p_i(x_i,x_{-i})$}  dependent  on its own strategy  $x_i$ and the rival strategies  $x_{-i}\triangleq \{x_j\}_{j \neq i}   $.  Let $x\triangleq (x_1,\cdots,x_n) \in
\mathbb{R}^d$ denote  the strategy profile  with $d \triangleq  \sum_{i=1}^n d_i$. Suppose the  $i$th player solves
\begin{align} \tag{\mbox{$\mathscr{P}_i(x_{-i})$}} \label{Ngame}
 \min_{x_i \in \mathbb{R}^{d_i}}  \ p_i(x_i,x_{-i})\triangleq f_i(x_i,x_{-i})+r_i(x_i),\end{align}
where $f_i(x )\triangleq \mathbb{E}\left[\psi_i(x ;\blue{\xi_i(\omega)}) \right] $,
 the random variable $\blue{\xi_i}: { \Omega} \to \blue{\mathbb{R}^{m_i}}$ is   defined on the probability space $({
\Omega}, {\cal F}, \mathbb{P}_i)$,    $\psi_i : \mathbb{R}^d \times \blue{\mathbb{R}^{m_i}}
\to \mathbb{R}$ is a scalar-valued  function, and
$\mathbb{E}[\cdot]$  denotes the expectation with respect to the probability
measure $ \mathbb{P}_i $,  and  $r_i(x_i)$ is a proper, closed, and convex function with an efficient proximal evaluation.   Since $r_i(x_i)$ can capture various forms of nonsmoothness (including $r_i(x_i) \triangleq {\bf 1}_{X_i}(x_i)$ allowing for convex constraint sets $X_i$), $\mathscr{P}$ is a nonsmooth Nash equilibrium problem.   A  Nash equilibrium (NE) of  $\mathscr{P}$  is a  tuple $x^*    \triangleq  \{x_i^*\}_{i=1}^n \in \mathbb{R}^{d} $  such that  for each $ i  \in \mathcal{N}$:
  \begin{align} p_i(x_i^*,x_{-i}^*) \leq p_i(x_i ,x_{-i}^*),  \quad
\forall x_i\in \mathbb{R}^{d_i} .  \nonumber \end{align}
In other  words, $x^*$
is an NE if no player  can  profit from unilaterally deviations.

 Our focus is two-fold: {\em (i) Development of variable sample-size stochastic
proximal gradient-response {\bf (VS-PGR)} and proximal best-response {\bf
(VS-PBR)} schemes with {\bf optimal (deterministic)} convergence  rates  as
well as iteration and oracle complexities; (ii) Extension to distributed
(consensus-based) regimes, referred to as {\bf (d-VS-PGR)} and {\bf
(d-VS-PBR)}, allowing for  resolving aggregative  games where
each player's  payoff depends on its  strategy and an aggregate of all players'
strategies  over   a  static  communication graph and {\bf linear} rates of
convergence are achieved by combining increasing number of consensus steps with
a growing batch-size of sampled gradients.}

  {\bf Prior research.} We discuss some relevant prior  research on gradient and best-response schemes
\blue{along with} their  distributed variants  for  continuous-strategy NEPs as well as  variance-reduced schemes for stochastic optimization.

(i) {\em GR schemes.} Early work considered convex NEPs
where player problems are convex, implying that  an NE~\cite[Chapter 1]{facchinei02finite}  is equivalent to a solution of  the associated variational inequality.  GR schemes have proven  useful in flow
control/routing games  (cf.~\cite{yin09nash2,alpcan2002game,pang2010design}),  often imposing   a
suitable monotonicity property on the concatenated  gradient map.  Merely monotone problems have been addressed via iterative regularization in deterministic~\cite{yin09nash2,kannan12distributed} and
stochastic~\cite{koshal13regularized} regimes,  while
extensions have addressed misspecification~\cite{jiang18distributed} and the
lack of Lipschitzian properties \cite{yousefian15selftuned}.

(ii) {\em BR schemes.} By observing  that an NE is     a fixed point
of the BR map  for convex NEPs (cf.~\cite{basar99dynamic}), we may apply
fixed-point or BR approaches, where  each player selects the strategy that
 maximizes   its payoff,  given rival strategies
(cf.~\cite{fudenberg98theory,basar99dynamic}). There have  been efforts to
extend such schemes to engineered settings (cf.~\cite{scutari2009mimo}), where
the  BR correspondence  can be expressed in a closed  form.      However, BR
schemes do not always lead to convergence to Nash equilibria  even in potential
convex NEPs (see~\cite{facchinei2011decomposition} for a counterexample).
The proximal BR schemes appear to have been  first discussed in \cite{FPang09}
where the  set of fixed points of  the proximal BR map is shown to be
equivalent to the set of Nash equilibria. In~\cite{facchinei2011decomposition},
several  asymptotically convergent   regularized  Gauss-Seidel BR schemes are suggested  for  generalized
potential games.  More recently, sampled BR schemes have been developed~\cite{pang2017two}   to
solve risk-averse  two-stage noncooperative games while the rate statements and
complexity bounds have been provided for a  distinctly different class
(specifically single-loop) of inexact  synchronous, asynchronous,
and delay-tolerant stochastic proximal BR schemes in~\cite{lei2017synchronous}. We
emphasize that in this scheme, each strongly convex expectation-valued
subproblem is solved by a stochastic gradient scheme and the overall iteration
complexity in terms of projected gradient steps is $\mathcal{O}(1/\epsilon)$,
in contrast with $\mathcal{O}(\ln(1/\epsilon))$ complexity obtained in the
present work, albeit in terms of sample-average subproblems.
 Finally, \blue{almost sure}   and mean  convergence of sequences produced by proximal BR
schemes is shown in \cite{lei2017asynchronous} for stochastic
and misspecified potential games.

(iii) {\em Consensus-based distributed schemes for NEPs.}
 In \blue{aggregative} NEPs~\cite{jensen10aggregative}, player payoffs are coupled through an aggregate of player strategies; however, players might not have access to the    aggregate and hence cannot compute payoffs or
gradients, precluding the direct  use of gradient/BR schemes. Inspired by
the consensus-based protocols for distributed optimization \vvs{in convex}
\cite{jakovetic2014fast,nedic2009distributed,lei2016primal,nedic2015distributed} \vvs{and nonconvex~\cite{scutari19distributed,scutari19parallel,lorenzo16next} regimes},
Koshal et al.~\cite{koshal2016distributed} developed distributed synchronous
and  asynchronous algorithms for such games  where  players  utilize  an
estimate of the aggregate  and update it by communicating with their neighbors (also see~\cite{pavel16distributed}).
Deterministic aggregative games  subject to coupling constraints are considered
in \cite{paccagnan2016distributed,parise20distributed, yi19distributed,grammatico17dynamic,belgioioso2017semi}, while
in~\cite{paccagnan2016distributed}, an asymmetric projection algorithm is adopted
for seeking  a  variational generalized Nash  equilibrium (GNE). In
\cite{belgioioso2017semi}, non-differentiable payoff functions  are considered
and  a semi-decentralized algorithm is presented for    finding  a zero of the
associated generalized equation.   However, in both \cite{belgioioso2017semi}
and \cite{paccagnan2016distributed}, an additional central node is required for
updating the Lagrange multiplier  associated with the coupling
constraints.  Distributed primal-dual algorithms are   proposed in
deterministic regimes \cite{zhu2016distributed} while the only known distributed gradient-based scheme for  generalized  stochastic NEPs was developed in  \cite{yu2017distributed}, which  considers asymptotic behavior under constant steps.  We focus on stochastic aggregative NEPs but consider both gradient and BR schemes while providing rate and complexity guarantees.

(iv) {\em Variance-reduced schemes.} There has been an effort to utilize
variance reduction schemes for  solving stochastic programs  within stochastic
gradient-based schemes, where the true gradient is replaced by the average of
an increasing  batch of sampled gradients, leading to a progressive reduction
of the variance of the sample-average  {gradient}.  Thus, such schemes  can
improve  the rates of convergence  or even  allow for recovering  deterministic
convergence rates   (in an expected value sense) if the batch size grows
sufficiently fast, as seen in
convex~\cite{shanbhag15budget,ghadimi2016accelerated,jofre2017variance}
and nonconvex optimization regimes~\cite{lei2020asynchronous,ghadimi2016accelerated,reddi2016proximal}.
However, there has been no known effort to apply such avenues for resolving
stochastic NEPs, particularly via BR schemes.

 {\bf Gaps and novelty.}  Prior
research on stochastic NEPs has largely resided in  standard
GR  approaches (without utilizing variance reduction) with either
little or no available  rate and complexity analysis for either  BR
schemes or distributed variants for  GR and BR schemes. In
addition,  most  prior rate statements show distinct gaps with
deterministic analogs. In this paper, we address the following gaps: (i) {\em
 {BR} schemes.} We provide a novel best-response scheme that can
address stochastic NEPs characterized by a suitable contractive
property; (ii) {\em Variance-reduction schemes.} By overlaying a variable
sample-size framework, both  GR and BR schemes
achieve deterministic rates of convergence with optimal or near-optimal oracle
complexities; (iii) {\em Distributed variants.} Finally, we
extend each scheme to a distributed regime capable of contending with
aggregative NEPs and prove that under suitable communication
requirements, the aforementioned geometric rates of convergence can be
retained.

\begin{table}  [htb]
  \centering
 \begin{tabular}{|c|c|c|c|c|}
        \hline
     Algorithm &   \multicolumn{2}{c|}{ {\bf VS-PGR}} &  \multicolumn{2}{c|}{ {\bf  VS-PBR}}     \\ \hline
     Assumption   &    \multicolumn{2}{c|}{  Strongly monotone} &   \multicolumn{2}{c|}{ Contractive proximal BR Map} \\ \hline
      $S_k$  &   $\lceil \rho^{-(k+1)}\rceil$  &  $\lceil (k+1)^{v}\rceil$
      &  {$\lceil \rho^{-(k+1)}\rceil$}  &  $\lceil (k+1)^{v}\rceil$ \\ \hline
  Rate    $ \mathbb{E}[ \|x_{k }-x^*\|^2  ] $  &  {$\mathcal{O}(\rho^k)$} &
$\mathcal{O}(k^{-v})+\mathcal{O}(q^k), q<1$  &  {$\mathcal{O}(\rho^k)$} & $\mathcal{O}(k^{-v})+\mathcal{O}(a^k), a<1$     \\ \hline
  Iteration  Complexity   &$\mathcal{O}(\ln(1/\epsilon))$ &
$\mathcal{O}((1/\epsilon)^{1/v})$ &
$\mathcal{O}(\ln(1/\epsilon))$ &
$\mathcal{O}((1/\epsilon)^{1/v})$  \\ \hline
   Oracle Complexity & $\mathcal{O}(1/\epsilon)$ & $\mathcal{O}((1/\epsilon)^{1+1/v})$  & $\mathcal{O}(1/\epsilon)$& $\mathcal{O}((1/\epsilon)^{1+1/v})$  \\ \hline
\end{tabular}
 \caption{(VS-PGR) and (VS-PBR) schemes  with $\rho \in (0,1) $ and $v > 0$ } \label{tab1}
\vspace{-0.2in}
\end{table}

\begin{table}  [htb]
  \centering
 \begin{tabular}{|c|c|c|c|c|}
        \hline
     Algorithm &   \multicolumn{2}{c|}{ {\bf d-VS-PGR}} &  \multicolumn{2}{c|}{ {\bf d-VS-PBR}}     \\ \hline
     Assumption   &    \multicolumn{2}{c|}{  Strongly monotone} &   \multicolumn{2}{c|}{ Contractive proximal BR Map} \\ \hline
      $S_k$  &   $\lceil \rho^{-(k+1)}\rceil$  &  $\lceil (k+1)^{v}\rceil$
      &  {$\lceil \rho^{-(k+1)}\rceil$}  &  $\lceil (k+1)^{v}\rceil$ \\ \hline
     Communication    $\tau_k$    & $k+1$ & $\lceil (k+1)^u\rceil $ & $k+1$ & $\lceil (k+1)^u\rceil $ \\ \hline
  Rate    $ \mathbb{E}[ \|x_{k }-x^*\|^2  ] $  &  {$\mathcal{O}(\rho^k)$} &
$\mathcal{O}(k^{-v})$  &  {$\mathcal{O}(\rho^k)$} & $\mathcal{O}(k^{-v})$     \\ \hline
  Iteration  Complexity   &$\mathcal{O}(\ln(1/\epsilon))$ &
$\mathcal{O}((1/\epsilon)^{1/v})$ &
$\mathcal{O}(\ln(1/\epsilon))$ &
$\mathcal{O}((1/\epsilon)^{1/v})$  \\ \hline
   Oracle Complexity & $\mathcal{O}(1/\epsilon)$ & $\mathcal{O}((1/\epsilon)^{1+1/v})$  & $\mathcal{O}(1/\epsilon)$& $\mathcal{O}((1/\epsilon)^{1+1/v})$  \\ \hline
	 Comm.  Complexity & $\mathcal{O}(\ln^2(1/\epsilon))$  &
$\mathcal{O}((1/\epsilon)^{(1+u)/v})$& $\mathcal{O}(\ln^2(1/\epsilon))$   &
$\mathcal{O}((1/\epsilon)^{(1+u)/v})$  \\ \hline
\end{tabular}
 \caption{(d-VS-PGR) and (d-VS-PBR) schemes  with $\rho\in (0,1),$ $v > 0,$ and $ u \in (0,1)$  }\label{tab2}
\vspace{-0.2in}
\end{table}

 {\bf Contributions.}  We summarize  our findings in
Tables~\ref{tab1}--~\ref{tab2} and discuss these next.

{(i). \em VS-PGR.}  In Section   \ref{sec:algorithm}, we  propose a
variable sample-size  proximal GR (VS-PGR) scheme, where
an increasing batch of sampled gradients is utilized at each iteration.
Under a strong monotonicity assumption, the mean-squared error (MSE)   admits
a linear rate of convergence
(\blue{Theorem} \ref{thm1}) when batch-sizes increase geometrically. We further show  in \blue{Theorem} \ref{thm2}  that  the iteration complexity
(no. of proximal evaluations) and oracle complexity (no. of  sampled gradients)
to achieve an $\epsilon-$NE denoted by $x$  satisfying  $ \mathbb{E}
[\| x-x^*\|^2 ] \leq \epsilon$  are respectively   $\mathcal{O}
(\ln(1/\epsilon))$  and  $\mathcal{O} \left((1/\epsilon)^{1+\delta}\right)$  with $\delta\geq 0$.
 In \blue{Corollary}  \ref{cor1},  under suitably chosen
algorithm parameters,   the  iteration and oracle complexity
to obtain an $\epsilon-$NE are shown to be optimal and are bounded by   $\mathcal{O}
(\kappa^2\ln(1/\epsilon))$ and by  $\mathcal{O} \left(
\kappa^2/\epsilon  \right)$,  where $\kappa$ denotes the condition
number.  Finally,  under a polynomially increasing sample-size $ \lceil
(k+1)^v\rceil , v>0$,  we show  in Lemma \ref{poly-lem} that  $ \mathbb{E}[
\|x_{k }-x^*\|^2  ] =\mathcal{O}\left(k^{-v}\right)+ \mathcal{O}\left(q^k\right)$  (where $q< 1$), and
establish that the  iteration and oracle complexity to obtain an
$\epsilon-$NE are $\mathcal{O} ( (1/\epsilon)^{1/v})$ and
$\mathcal{O} \left((1/\epsilon)^{1+1/v}  \right)$, respectively.

{(ii). \em Distributed VS-PGR.} In Section \ref{sec:agg}, we design  a distributed  VS-PGR  scheme (see Algorithm \ref{algo-aggregative}) to compute an  NE  of an aggregative  stochastic
NEP over a static  communication graph where players combine  \blue{variable sample-size}
  proximal  GR  with a consensus update for learning the aggregate.
By suitably increasing the number of consensus steps and sample-size at
each iteration, \blue{(d-VS-PGR)} achieves a linear rate of
convergence (\blue{Theorem} \ref{agg-thm1})  while  in \blue{Theorem} \ref{agg-thm2} and
\blue{Corollary}  \ref{cor2}, the iteration, oracle, and communication complexity to compute
an $\epsilon$-NE are proven to be  $\mathcal{O} (\ln(1/\epsilon))$,
$\mathcal{O} \left( 1/\epsilon \right)$, and $\mathcal{O}
\left(\ln^2(1/\epsilon)\right)$, respectively.  With  polynomially increasing
communication rounds $  \lceil (k+1)^{u} \rceil$ and  sample-size $ \left\lceil
(k+1)^v \right \rceil$ for some $u\in (0,1)$ and $ v>0$, \blue{(d-VS-PGR)}  achieves a polynomial
rate     $ \mathbb{E}[ \|x_{k}-x^*\|^2  ] =\mathcal{O}\left(
k^{-v}\right)$ associated with the iteration, communication,  and oracle
complexity given by  $\mathcal{O} (
(1/\epsilon)^{1/v})$, $\mathcal{O} ( (1/\epsilon)^{(u+1)/v}),$ and
$\mathcal{O} \left((1/\epsilon)^{1+1/v}  \right)$,    respectively.

 {(iii). \em VS-PBR.}  In Section  {\ref{VS-PBR}},  we develop a variable
sample-size proximal  BR   (VS-PBR) scheme (see \blue{Algorithm} \ref
{inexact-sbr-cont})   when the proximal BR map is contractive and
requires that  each  player solves a   sample-average  BR  problem per
step. In \blue{Theorem} \ref{thm4}, the generated iterates converge
to the NE  at a linear rate  in the mean-squared sense when the   sample-size     for
computing the sample-average  payoff  increases geometrically,  leading to an iteration (no. of deter.
opt.  problems solved) and oracle complexity (no. of  samples)  to achieve an
$\epsilon-$NE  of $\mathcal{O}( \ln(1/\epsilon))$
and   $\mathcal{O} ( 1/\epsilon)  $ respectively. Akin to Section  \ref{sec:algorithm}, we  show in
 \blue{Corollary} \ref{poly-lem2} that  when  the sample-size increases at a    polynomial rate of    $ \left\lceil  (k+1)^v \right \rceil$, $ \mathbb{E}[ \|x_{k}-x^*\|^2  ] =\mathcal{O}\left( k^{-v}\right) + \mathcal{O}\left(a^k\right)  $ (where $a < 1$) with the iteration and oracle complexities   are  $\mathcal{O}(1/\epsilon^{1/v})$ and $\mathcal{O}(1/\epsilon^{1+1/v})$,
	respectively.

{(iv). \em Distributed VS-PBR.} In Section~\ref{sec:agg2}, we  design
a distributed  VS-PBR scheme (see \blue{Algorithm}~\ref{algo-aggregative2}) to
compute an NE  of an aggregative  NEP with
contractive proximal BR maps, akin to ({\bf d-VS-PGR}) where the
aggregate is estimated by taking multiple consensus steps   while
the proximal BR is approximated by  solving a  sample-average
BR problem.  When the  number of consensus
steps and  sample-sizes  are raised suitably fast,  the mean-squared error diminishes at a geometric rate (\blue{Proposition} \ref{prp3}). We further show in \blue{Theorem}~\ref{thm-iter-comp} that the iteration, oracle, and communication complexity to compute
an $\epsilon$-NE are    $\mathcal{O} (\ln(1/\epsilon))$, $\mathcal{O}
\left( 1/\epsilon \right)$, and $\mathcal{O}
\left(\ln^2(1/\epsilon)\right)$, respectively.

  {\bf Notation:}  A vector $x$ is assumed to be a column vector while $x^T$ denotes its transpose. $\|x\|$  denotes the Euclidean vector norm, i.e., $\|x\|=\sqrt{x^Tx}$.   We abbreviate  ``almost surely'' by \textit{a.s.} and  $\lceil x \rceil $ denotes   the smallest integer greater than $x$ for $x \in \mathbb{R}$.  For a closed convex function    $r(\cdot)$, the prox.~ operator   is defined  by \eqref{proximal}   for $\alpha>0$:
\begin{equation}\label{proximal}
\textrm{prox}_{\alpha r}(x)  \triangleq  \argmin_{y} \left( r(y)+{1\over 2\alpha} \|y-x\|^2\right) .
\end{equation}

\section{VS-PGR Scheme  and its  distributed variant} \label{sec:monotone}
This section considers the development of a variable sample-size  proximal
stochastic gradient-response scheme for   a class of   strongly  monotone
NEPs.  We derive rate and complexity statements when players can observe rival
strategies in Section~\ref{sec:algorithm} and provide analogous statements for
a distributed variant in Section~\ref{sec:agg} for an aggregative game where
players overlay an additional consensus phase for learning the aggregate


\subsection{Variable sample-size proximal stochastic gradient-response scheme}\label{sec:algorithm}

 \blue{Before the algorithm statements, we firstly} impose the following conditions   on $\mathscr{P}$.

\begin{assumption}~\label{assump-play-prob}
\blue{For each $i\in \mathcal{N},$} (i)  $r_i$ is lower semicontinuous and convex with the effective domain
denoted by $\Rscr_i \triangleq \mbox{dom}(r_i)$;
(ii) for every fixed $x_{-i} \in \Rscr_{-i}  \triangleq \prod_{j \neq i} \Rscr_j$, $f_i(x_i,x_{-i})$ is
  continuously differential and convex in $x_i$ \vvs{on an open set containing $\mathcal{R}_i$};
(iii) for all  $x_{-i} \in  \Rscr_{-i} $ and  any $\blue{\xi_i\in\mathbb{R}^{m_i}}$,
\blue{$  \psi_i(x_i, x_{-i};\xi_i)$ is differentiable  in $x_i$ \us{on an open set containing $\Rscr_i$}.}
\end{assumption}

Define  $ G(x ) = \left(    \nabla_{x_i} f_i(x)    \right) _{i=1}^n  $. The following lemma establishes that a tuple $x^*$ is   an NE      if and only if it is a fixed point of a  suitable map.  
\begin{lemma}[{\bf Equivalence between NE and fixed point  of proximal response  map}] \label{lem1} Given   the stochastic NEP  \blue{$\mathscr{P}$}, suppose Assumption~\ref{assump-play-prob}(i) and \ref{assump-play-prob}(ii)
hold  for each player $i\in \mathcal{N}$. Let $r(x) \triangleq ( r_i(x_i))_{i=1}^n.$
 Then $x^* \in X$ is an NE   if and only if  $x^*$   is a fixed point of  $\textrm{prox}_{\alpha r}(x -\alpha G(x ))$, i.e.,
\begin{align}\label{FP}
x^*=\textrm{prox}_{\alpha r}(x^*-\alpha G(x^*)), \quad  \forall \alpha>0.
\end{align}
\end{lemma}
{\bf Proof.} Note that for \blue{each} $i\in \mathcal{N},$ $f_i(x_i,x^*_{-i})$ and $r_i(x_i)$ is convex in $x_i.$ Then  $x_i^*$ is an optimal solution of  $f_i(x_i,x^*_{-i})+r_i(x_i)$ if and only if for any $\alpha>0$,
$x_i^*=\textrm{prox}_{\alpha r_i}(x_i^*-\alpha  \nabla_{x_i} f_i (x^*)),$ which
by  concatenation  for $i = 1, \cdots, n$ leads to \eqref{FP}. \hfill $\Box$

 Suppose the iteration  index  is given by $k$. Player $i$  at  iteration $k$  holds an estimate $x_{i,k} \in\mathbb{R}^{d_i}$   of the  equilibrium strategy $x_i^*.$ We consider a variable sample-size
generalization of the standard  proximal stochastic gradient method, in which
$S_k$ number of sampled gradients are utilized at  iteration $k.$    For any $i\in \mathcal{N}$, given    $S_k$ realizations  $\nabla_{x_i} \psi_i(x_k;\xi_{i,k}^1), \dots$, $\nabla_{x_i} \psi_i(x_k;\xi_{i,k}^{S_k})$,   $x_{i,0}\in   \Rscr_i ,$
player $i$ updates $x_{i,k+1}$ as follows:
\begin{align}\label{VSSGD}
x_{i,k+1}=\textrm{prox}_{\alpha  r_i} \Bigg[ x_{ i,k}-{\alpha \over S_k}  { \sum_{p=1}^{S_k}  \nabla_{x_i}\psi_i(x_k;\xi_{i,k}^p)  }\Bigg],\end{align}
where  $\alpha>0 $ is the step size. If $w_{k}^p \triangleq \blue{  \big( \nabla_{x_i}\psi_i(x_k;\xi_{i,k}^p)  \big)_{i=1}^n    }-G(x_k)  $  and $\bar{w}_{k,S_k}\triangleq{1\over S_k} \sum_{p=1}^{S_k} w_{k}^p,$
then by  concatenating \eqref{VSSGD} for $i = 1, \cdots, n$, we obtain the compact form:  \begin{align} \tag{VS-PGR} \label{VSSA}
x_{k+1}=\textrm{prox}_{\alpha  r} \left[ x_k-\alpha \left( G(x_k )+\bar{w}_{k,S_k} \right) \right].
\end{align}
We impose  the following conditions on the  gradient mapping   $G(x)$ and  noise $\bar{w}_{k,S_k}$ and rely on $\mathcal{F}_k$, defined as $\mathcal{F}_k\triangleq \sigma\{x_0,x_1,\cdots,x_{k}\}$.
\begin{assumption}~\label{assp-noise}
(i)  $G(x)$ is  $L$-Lipschitz continuous, i.e.,
 $\| G(x)-G(y)\| \leq L \|x-y\|$ for all $x,y\in \Rscr  {\triangleq  \Pi_{j=1}^n \Rscr_j}.$
(ii) $G(x)$ is   $\eta$-strongly monotone, i.e.,
$(G(x)-G(y))^T(x-y) \geq \eta  \|x-y\|^2$ for all $x,y\in \Rscr.$
  (iii) There \blue{exist}  constants $\nu_1, \nu_2 \geq 0$ such that for any $k\geq 0$,
$\mathbb{E}[ \|\bar{w}_{k,S_k} \|^2\mid \mathcal{F}_k] \leq  \tfrac{\nu_1^2\|x_k\|^2 + \nu_2^2}{S_k}$, a.s. .
\end{assumption}

\blue{ Define $T  (x)\triangleq \textrm{prox}_{\alpha r}(x -\alpha G(x )) $.
Then  by the non-expansivity of the proximal operator,
 Assumptions \ref{assp-noise}(i) and \ref{assp-noise}(ii), for any $\alpha<{2\eta  \over   L^2 }$ and $x, y \in \mathcal R$,
\begin{align*}
\|  T(x) -  T(y)\|& \leq  \left \| x -y- \alpha \big(  G(x )   - G(y ) \big)\right\|
   \leq   \sqrt{1-2\alpha\eta +\alpha^2   L ^2  } \|x-y\| <\|x-y\|. 
\end{align*}
This implies that $T(x)$ is a contractive map.  Since    $\mathcal{R}$ is a closed convex  set,
 from \cite[Theorem 2.1.21]{facchinei02finite} it follows that $T(x)$   has a unique  fixed point in $\mathcal{R}$.
 Therefore, the stochastic NEP $\mathscr{P}$  has a unique Nash equilibrium $x^*.$}
We now  establish a simple  recursion  for  the conditional MSE in  terms of $S_k,$
step size $\alpha $, and problem parameters.

\begin{lemma}\label{lem2}
Consider   \eqref{VSSA}. Suppose Assumptions \ref{assump-play-prob}--\ref{assp-noise} hold,
$\{S_k\}$ is a non-decreasing sequence \blue{with} $S_0 \geq 1/\alpha^2$,
\blue{$\tilde{L}  \triangleq \sqrt{ 1  + 2 (1+2\alpha^2)\nu_1^2+2L^2}$, $ \alpha\leq L$,} and
 $\nu^2 \blue{=} {2(1+2\alpha^2)\nu_1^2\|x^*\|^2 + (1+2\alpha^2)\nu_2^2}.$
 \blue{Then the following holds for any $ k\geq 0:$}
\begin{align}\label{prp-bd}
\mathbb{E}[ \|x_{k+1}-x^*\|^2 |\mathcal{F}_k] &
\leq  \left(1- 2\alpha\eta +\alpha ^2 \tilde L^2\right) \|   x_k -x^* \|^2  + \tfrac{\nu^2}{S_k},   ~ a.s. \ .
\end{align}
\end{lemma}
{\bf Proof.} The proof can be found in Appendix \ref{ap:lem2}.\hfill $\Box$

\subsubsection{Geometrically Increasing Sample  Sizes}
 {We  begin with   a  preliminary lemma that will be used in the rate analysis.
\begin{lemma} \label{lem-recur} Let the sequence $\{v_k\}_{k\geq 0}$ with  initial value $v_0 \leq c_0$ satisfy the following  recursion \eqref{lem-recur1} for some $q,\rho\in (0,1)$:
\begin{align}\label{lem-recur1}  v_{k+1} \leq q v_k+c_1 \rho^{k+1} ,\quad \forall k\geq 0. \end{align}
Then for any $k\geq 0, $  we have
\\(i) $v_k \leq \big(c_0+{c_1 \over \max\{ q/\rho, \rho /q \}-1} \big) \max \{\rho,q\}^k $ when $\rho \neq q$;
\\ (ii)  $v_k      \leq   \left(c_0+ {c_1 \over \ln((\tilde{q}/q)^e)}\right)\tilde{q}^k $
for any $\tilde{q}\in (q,1)$ when $\rho=q$.
\end{lemma}
{\bf Proof.}  The proof is given in Appendix \ref{ap:lem-recur}.\hfill $\Box$

Based on Lemma \ref{lem2} and \ref{lem-recur},  we  prove  linear  convergence  of   {the iterates generated by   \eqref{VSSA}  in a  mean-squared sense} with geometrically   increasing batch-size.
\begin{theorem}[{\bf Linear convergence rate of VS-PGR}] \label{thm1}
Let \eqref{VSSA} be applied to $\mathscr{P}$, where
\blue{$S_k   \triangleq     \left \lceil  \alpha^{-2}\rho^{-(k+1)}  \right\rceil$} for some $\rho\in (0,1)$,
 and $\mathbb{E}[ \|x_0-x^*\|^2] \leq C$ for some  $C>0.$
 Suppose Assumptions \ref{assump-play-prob} and \ref{assp-noise} hold.
Let $\alpha \in (0,2\eta/\tilde{L}^2 )$, \blue{where $\tilde{L}$ is defined in Lemma \ref{lem2}.} Then
 $q\triangleq 1-2\alpha\eta +\alpha^2   \tilde{L}^2<1$ and  for any $k\geq 0$.

\noindent (i) If $\rho\neq q,$ then $\mathbb{E}[ \|x_{k }-x^*\|^2]\leq
  C(\rho,q) \max \{\rho,q\}^k  $ where  $  C(\rho,q)  \triangleq  C+{\tfrac{\alpha^2 \nu^2}{\max \{\rho/q, q/\rho\}-1}} .$

\noindent (ii) If $\rho= q,$ then for  any $\tilde{q}\in(q,1),$ $\mathbb{E}[ \|x_{k }-x^*\|^2]\leq \widetilde{D} \tilde{q}^k, $ where $\widetilde{D} \triangleq  C+\tfrac{{\alpha^2\nu^2}}{\ln((\tilde{q}/q)^e)}. $
\end{theorem}
{\bf Proof.}By  definition, $q \in (0,1)$ when $\alpha \in (0,2\eta/  \tilde{L}^2 )$. Then by
taking unconditional expectations on both sides of  Eqn.~\eqref{prp-bd} and using
\blue{$S_k \geq   \alpha^{-2}\rho^{-(k+1)}   $}, we obtain that
\[\mathbb{E}[ \|x_{k+1}-x^*\|^2  ]    \leq  q\mathbb{E}[ \|x_{k }-x^*\|^2  ]
 +   \alpha^2 \nu^2  \rho^{ k+1 } ,\quad \forall k\geq 0 \]
 Then by using Lemma \ref{lem-recur}, we obtain the results. \hfill $\Box$

  Next, we examine the iteration  (no. of proximal evaluations) and oracle complexity (no. of   sampled gradients) of \eqref{VSSA} to compute an  $\epsilon$-Nash equilibrium.  We  refer to  a random strategy
profile $x: { \Omega} \to \Real^{n}$  as an  $\epsilon-$NE if $\mathbb{E}[\|x-x^*\|^2]\leq \epsilon$.
\begin{theorem}[{\bf Iteration and Oracle Complexity}] \label{thm2}
 Suppose Assumptions \ref{assump-play-prob} and \ref{assp-noise} hold.  Let  \eqref{VSSA} be applied to $\mathscr{P}$, where  $\alpha\in (0,2\eta/\blue{\tilde{L}}^2)$   \blue{with $\tilde{L}$  defined in Lemma \ref{lem2}},  \blue{$S_k   \triangleq     \left \lceil  \alpha^{-2}\rho^{-(k+1)}  \right\rceil$ with $\rho\in (0,1)$,} and $\mathbb{E}[ \|x_0-x^*\|^2] \leq C$ for some  $C>0.$
 Set $\tilde{q}\in(q,1)$.
  Then    the number of  proximal evaluations  needed to obtain an $\epsilon-$NE is bounded by $K(\epsilon)$, defined as
\begin{align}\label{rate-prox}
K(\epsilon) \triangleq \begin{cases}  &  {1\over \ln \left({1/q} \right)}  \ln \left( \big(C+{{\alpha^2 \nu^2} \rho \over  q-\rho}\big) \epsilon^{-1}\right) \qquad \quad  {\rm if ~} \rho<q<1,\\
&   {1 \over \ln \left(1/\tilde{q} \right)}
\ln\Big( \big( C+{{\alpha^2 \nu^2} \over \ln((\tilde{q}/q)^e)} \big)  \epsilon^{-1}   \Big) \qquad   ~{\rm if ~} q=\rho,\\
& {1 \over \ln \left({1/\rho} \right)} \ln \left( \big(C+{{\alpha^2 \nu^2} q \over  \rho -q}\big) \epsilon^{-1}\right) \qquad \quad {\rm if ~} q<\rho<1,
\end{cases}
\end{align}
and    the number of sampled gradients required is bounded by $M(\epsilon)$, defined as
\begin{align}\label{rate-SFO}
M(\epsilon) \triangleq \begin{cases}  &
 { 1  \over  \rho \ln (1/\rho)} \left( \big(C+{{\alpha^2 \nu^2} \rho \over  q-\rho}\big) \epsilon^{-1}\right)
^{\tfrac{\ln(1/\rho)}{\ln(1/q)}} +K (\epsilon) \qquad {\rm if ~} \rho<q<1,\\&
{  1 \over  q \ln (1/q)} \left(  \big(C+{{\alpha^2 \nu^2} \over \ln((\tilde{q}/q)^e)}  \big)  \epsilon^{-1} \right)^{\tfrac{\ln(1/q)}{\ln(1/ \tilde{q} )}} + K (\epsilon)  \quad~{\rm ~ if ~} \rho=q, \\
& { 1  \over  \rho \ln (1/\rho)}   \big(C+{{\alpha^2 \nu^2}q \over  \rho-q}\big) \epsilon^{-1} +K(\epsilon)  \qquad \qquad
\qquad  {\rm ~ if ~} q<\rho<1.
\end{cases}
\end{align}
\end{theorem}
{\bf Proof.}We first consider the case $\rho\neq q.$    {From Theorem \ref{thm1}(i), it follows that  for any $k\geq K_1(\epsilon) \triangleq {\ln\left({ C(\rho,q)/\epsilon}\right) \over \ln \left({1/ \max \{\rho,q\}} \right)},$ $ \mathbb{E}[ \|x_{k }-x^*\|^2]\leq \epsilon.$ }
Then by $  C(\rho,q)  =  C+{{\alpha^2 \nu^2} \over \max \{\rho/q, q/\rho\}-1} ,$ we obtain the     bound  on   {the iteration complexity}  defined  in Eqn. \eqref{rate-prox} for cases $ \rho<q<1$ and $q<\rho<1$.
  For any $\lambda > 1$ and  positive integer $K$,  we have that
\begin{align}\label{exponent-sum}
\sum_{k=0}^K \lambda^k &  \leq \int_0^{K+1} \lambda^x dx \leq
\tfrac{\lambda^{K+1}} {\ln(\lambda)}.
\end{align}
 Therefore,  we achieve the following bound on the   number of  samples utilized:
\begin{align*}    & \sum_{k=0}^{ K_1(\epsilon) -1}  S_k    \leq   \sum_{k=0}^{ K_1(\epsilon)-1}  \rho^{-(k+1)}   + K_1(\epsilon)
  \leq  \tfrac{ 1}{\rho \ln (1/\rho)} \rho^{- K_1(\epsilon)} + K_1(\epsilon) .
\end{align*}
Note that for any $0<\epsilon, p<1, c_1>0$,  the following holds:
\begin{equation}\label{rel-exp}
\begin{split} & \rho^{-  \tfrac{\ln( c_1/\epsilon )}{\ln(1/p )} }
	  = \left(e^{\ln(\rho^{-1} )}\right)^{  \tfrac{\ln( c_1/\epsilon )}{\ln(1/p )} }
			  =  e^{\ln(c_1/\epsilon ))^{\tfrac{\ln(1/\rho)}{\ln(1/p )}}}
			  	=   {(c_1/\epsilon )^{\tfrac{\ln(1/\rho)}{\ln(1/p )}}}  .
\end{split}
\end{equation}
Thus,  the number of sampled gradients  required to obtain  an  $\epsilon-$NE  is   bounded by
$  { 1  \over  \rho \ln (1/\rho)}\left( C(\rho,q)\over \epsilon\right)^{\tfrac{\ln(1/\rho)}{\ln(1/\max\{\rho,q\} )}} +K_1(\epsilon) .$
Thus,  we achieve the    bound given in equation  \eqref{rate-SFO} for cases $ \rho<q<1$ and $q<\rho<1$.

We now prove the results for the case $\rho=q.$   {From   Theorem \ref{thm1}(ii) it follows that   for  any $\tilde{q}\in(q,1)$ and $k\geq K_2(\epsilon) \triangleq   {\ln\big({ \widetilde{D}/  \epsilon} \big) \over \ln \left(1/\tilde{q} \right)}  ,$ $\mathbb{E}[ \|x_{k }-x^*\|^2]\leq \epsilon $.}
Then we achieve the    bound given in  Eqn.  \eqref{rate-prox} for  {the case}  $ \rho=q$.
 Therefore,  we may \vvs{obtain a bound on sampled gradients \eqref{rate-SFO} when $ \rho=q$ as follows.
    $\sum_{k=0}^{ K_2(\epsilon) -1}  S_k     \leq  \tfrac{   q^{- K_2(\epsilon)} }{q \ln (1/q)} + K_2(\epsilon)  =  \tfrac{  1}{q \ln (1/q)} \left( \widetilde{D} / \epsilon  \right)^{\tfrac{\ln(1/q)}{\ln(1/ \tilde{q} )}} + K_2(\epsilon) .$}
\hfill $\Box$

  The above theorem establishes that the iteration  and oracle complexity to achieve an $\epsilon-$NE    are    $\mathcal{O} (\ln(1/\epsilon))$  and  $\mathcal{O} ( (1/\epsilon)^{1+\delta})$, where  $\delta=0$ when $\rho\in (q,1),$  $\delta= \tfrac{\ln(q/\rho)}{\ln(1/q)}~{\rm when ~} \rho<q<1$,
and $\delta=\tfrac{\ln(\tilde{q} /q)}{\ln(1/ \tilde{q} )}$ when $ \rho=q.$
In the following, we further examine  the influence of the condition number on the iteration and oracle  complexity.

\begin{cor} \label{cor1}
Let  the scheme \eqref{VSSA} be applied to $\mathscr{P}$, where $\mathbb{E}[ \|x_0-x^*\|^2] \leq C$.
 Suppose Assumptions \ref{assump-play-prob} and \ref{assp-noise} hold. Define the condition  number $\tilde{\kappa}\triangleq \tfrac{\tilde{L}}{\eta}$. Set $\alpha =\tfrac{\eta}{ \tilde{L}^2}  $ and \blue{$S_k   \triangleq     \left \lceil  \alpha^{-2}\rho^{-(k+1)}  \right\rceil$}  with $\rho=1-{1\over 2\tkappa^2}$.
Then   the iteration and oracle  complexity to obtain an $\epsilon-$NE are bounded by
  $\mathcal{O} (\blue{\tkappa}^2\ln(1/\epsilon))$ and by  $\mathcal{O} \left( \blue{\tkappa}^2/\epsilon  \right)$, respectively.
\end{cor}
{\bf Proof.} By $\alpha={\eta\over \blue{\tilde{L}}^2} $  and $\tkappa={\tilde{L}\over \eta}$, we obtain  that
$q=1-2\eta \alpha+\alpha^2\tilde{L}^2=1-{\eta^2\over \tilde{L}^2}=1-{1\over \tkappa^2}.$
 Note that  $\rho>q $ by  $\rho=1-{1\over 2\tkappa^2}$. Thus,
 ${\alpha^2  q \over \rho - q  }    \leq   {2\left(\eta\over \tilde{L}^2\right) ^2    \tkappa^2  }=   2    /\tilde{L}^2  .$
 \blue{By  noticing  that $\ln(1+x) \geq x/(x+1)$ for any $x\geq 0$, and using  $ {1\over \rho} = 1+{1\over 2\tkappa^2-1}$,} we have $\ln \left({1/\rho} \right) \geq  \tfrac{\tfrac{1}{2\tkappa^2-1}}{(1+\tfrac{1}{2\tkappa^2-1})}=\tfrac{1}{2 \tkappa^2}.$
Thus, by  Eqns. \eqref{rate-prox}  and \eqref{rate-SFO} for   the case $ q<\rho<1$,  the results hold by the following:
\begin{align*}
  K(\epsilon)& =\tfrac{\ln \big(C+{\alpha^2 \nu^2 q /(  \rho -q)}\big) +\ln(1/\epsilon )}{\ln \left({1/\rho} \right)}
 \leq 2\left(  \ln \big(C+2\nu^2 /\tilde{L}^2\big)  +\ln(1/\epsilon )  \right)  \tkappa^2
 =\mathcal{O} (\blue{\tkappa}^2)\ln(1/\epsilon), \\
M(\epsilon)&=\tfrac{  C+{\alpha^2 \nu^2 q /(  \rho -q)}}{\rho \ln (1/\rho)} \epsilon^{-1}   +K(\epsilon)
 \leq \tfrac{ C+2\nu^2 /\tilde{L}^2}{\epsilon}2 \tkappa^2 \left(1+{1\over 2\tkappa^2-1} \right)   +K(\epsilon)  =\mathcal{O} (\tfrac{\blue{\tkappa}^2}{\epsilon}) . ~
\end{align*}
\hfill $\Box$

\subsubsection{Polynomially Increasing Sample-Size}\label{VSPGR-ploy}

We now investigate the convergence properties of the scheme \eqref{VSSA} with polynomially increasing sample size.
 {We first  prove a preliminary result.}

{
 \begin{lemma}\label{bd-cqv}
Consider the function  {$d(x) = q^{x^u} x^{v}$ where $q \in (0,1)$, $x >0$, $u\in (0,1]$, and $v>0$.} Then $d(x)$ is unimodal on $\Real_+$ with a unique maximizer given by $x^* = \tfrac{v}{\ln(1/q)}$. Furthermore, $q^x \leq c_{q,v}  {x^{-v}}$ for all $x \in \Real_+$ where $c_{q,v} \triangleq   e^{- v/u }{\left(\tfrac{v}{u\ln(1/q)}\right)^{ v/u }}.$
\end{lemma}
{\bf Proof.} We begin by noting that $d'(x) =   \ln(q)  u  q^{x^u}x^{u-1}x^v+v  q^{x^u}  x^{v-1}= q^{x^u}  x^{v-1}(v-u\ln(1/q) x^u  ) $ and $d'(x^*) = 0$ if $x^* = \left(\tfrac{v}{u\ln(1/q)}\right)^{1/u}.$ Unimodality follows by noting that  $d(0) = 0$, $d'(x) > 0$ if $x \in (0,x^*)$, and $d'(x) < 0$ when $x > x^*$. It follows that
$$   {c_{q,v} \triangleq \max_{x  \geq   0}d(x)=q^{(x^*)^u} (x^*)^v} = q^{\tfrac{v/u }{\ln(1/q) }}{\left(\tfrac{v}{u\ln(1/q)}\right)^{ v/u }} = e^{- v/u }{\left(\tfrac{v}{u\ln(1/q)}\right)^{ v/u }}. \qquad \Box$$}

\begin{proposition}\label{poly-lem}
Let    \eqref{VSSA} be applied to $\mathscr{P}$, where \blue{$S_k   \triangleq     \left \lceil  \alpha^{-2}(k+1)^v \right\rceil$} with  $v>0$, and $\mathbb{E}[ \|x_0-x^*\|^2] \leq C$.  Suppose Assumptions~\ref{assump-play-prob} and \ref{assp-noise} hold. Let $\alpha \in (0,2\eta/\tilde{L}^2 )$  and $q\triangleq 1-2\alpha\eta +\alpha^2\tilde{L}^2 $
\blue{with $\tilde{L}$   defined in Lemma \ref{lem2}.} Then
\begin{align}\label{poly-rate}
&\mathbb{E}[ \|x_k-x^*\|^2  ]      \leq q^{k} \left(C+ \alpha^2\nu^2 \tfrac{e^{2v}q^{-1}-1}{1-q }\right)+ \tfrac{ 2\alpha ^2\nu^2 q^{-1}}{\ln(1/q) } k^{-v} ,\quad \forall k\geq 0 .
\end{align}
In addition, the  iteration and oracle complexity to obtain an $\epsilon-$NE are  {$\mathcal{O} (v (1/\epsilon)^{1/v})$ and   $\mathcal{O} \left(e^{v} v^v(1/\epsilon)^{1+1/v}  \right)$}, respectively.
\end{proposition}
{\bf Proof.}The proof can be found in Appendix \ref{ap:poly-lem}.  \hfill $\Box$

 \begin{remark}\label{rem-poly}  \vvs{By} Proposition  \ref{poly-lem}, \vvs{the constant in both complexity bounds grows at an exponential rate with $v$} while the rate improves at a polynomial rate. Choosing $v$ requires \vvs{trading off} available computational resources with the \vvs{the cost of}  generating   sample-averaged gradients with large sample-sizes,
  a focus of ongoing research.
\end{remark}

\subsection{Distributed VS-PGR for Aggregative Games}\label{sec:agg}
Next,  we consider a structured nonsmooth stochastic  aggregative game $\mathscr{P}^{\rm agg}$, where player $i\in \mathcal{N}$ solves the following parametrized problem:
\begin{align} \tag{\mbox{$\mathscr{P}^{\rm agg}_i(x_{-i})$}} \label{Ngame_agg}
\min_{x_i \in \mathbb{R}^d}  \ p^{\rm agg}_i(x_i,x_{-i})\triangleq f_i(x_i, \bar{x} )+r_i(x_i), \end{align}
where   $\bar{x} \triangleq  \sum_{i=1}^n x_i$ denotes the aggregate of all players'  strategies and  $f_i(x_i, \bar{x} ) \triangleq \mathbb{E}\left[\psi_i(x_i,x_i+\bar{x}_{-i} ;\blue{\xi_i} ) \right] $    is expectation-valued with  $\bar{x}_{-i}  \triangleq   \sum_{j=1,j\neq i}^n x_j$ and the random variable \blue{$\xi_i: { \Omega} \to\mathbb{R}^{m_i}$.}
  We impose the following assumptions   on $\mathscr{P}^{\rm agg}$.
\begin{assumption}~\label{assump-agg1}   \blue{(i)The function $r_i$ is lower semicontinuous and convex with effective domain
denoted by $\mathcal{R}_i$, which is  required to be  compact};
(ii)   For any \blue{$x_{-i} \in  \Rscr_{-i} $, $f_i(x_i,x_i+\bar{x}_{-i})$} is continuously differentiable  and convex in $x_i \in \mathcal{R}_i$;
(iii) For any \blue{$x_{-i} \in  \Rscr_{-i} $} and  any $\blue{\xi_i\in\mathbb{R}^{m_i}}$,
 \blue{$  \psi_i(x_i, x_i+\bar{x}_{-i};\xi_i)$ is differentiable  in $x_i \in \mathcal{R}_i$.}
\end{assumption}


\subsubsection{Algorithm Design}
In this part, we design a distributed algorithm to  compute an NE of  $\mathscr{P}^{\rm agg}$, where each player may exchange information with its
local neighbors, and subsequently update  its estimate of  the aggregate
and the equilibrium strategy.  The interaction  among players  is defined  by an undirected  graph $\mathcal{G}=(\mathcal{N},\mathcal{E} )$, where  $\mathcal{N} \triangleq\{1,\dots,n\}$  is the set of  players and  $\mathcal{E}$ is  the set of undirected edges between players. The set of  neighbors of  player  $i$  is defined as  $\mathcal{N}_i=\{ j\in \mathcal{N}:
(i,j)\in \mathcal{E}\}$, and player $i$ is assumed to be a neighbor of itself.
  Define  the adjacency matrix $A=[a_{ij}]_{i,j=1}^n$,
where $a_{ij}>0$  if $j\in \mathcal{N}_i$ and $a_{ij}=0$, otherwise.    A
path in $\mathcal{G}$ with length $p$ from $i_1$ to $i_{p+1}$ is a  sequence of
distinct nodes, $i_1i_2\dots i_{p+1}$, such that  $(i_m, i_{m+1})\in
\mathcal{E}$, for all $m=1,\dots,p$. The graph $\mathcal{G}$ is termed {\it
connected} if there is a  path between  any two distinct players $i,j\in\mathcal{N}$.  Though  each player does not have access to all players' strategies, it may
estimate   the aggregate  $\bar{x} $  by communicating with its  neighbors.

 \blue{Suppose the iteration index is given $k.$}
Player $i$  at time $k$ holds an estimate $x_{i,k}$ for its equilibrium
strategy  and an estimate $v_{i,k}$ for the average of the aggregate.  To
overcome the fact that the communication network is sparse, we assume that to
compute $v_{i,k+1},  $    players communicate  $\tau_k$ rounds rather than once  at
major iteration $k+1 $. The   strategy of each player is updated by a variable
sample-size proximal stochastic gradient scheme characterized by  \eqref{alg-agg-strategy1}   dependent on the constant step size  $\alpha>0$    \blue{and    $S_k $  number of   sampled gradients    $\nabla_{x_i}\psi_i \big(x_{i,k}, n\hat{v}_{i,k};\xi_{i,k}^p\big) $, where     $  \big \{\xi_{i,k}^p\big\}_{ p=1}^{ S_k} $ denote the  independent and identically distributed (i.i.d.)
 random  realizations of $\xi_i$.}  We specify the scheme in Algorithm \ref{algo-aggregative}.
\begin{algorithm}[htbp]
\caption{Distributed VS-PGR for Aggregative   SNEPs} \label{algo-aggregative}
 {\it Initialize:} Set $k =0$, and  $ v_{ i,0} =x_{ i,0}  \in \mathcal{R}_i$  for  any $i \in\mathcal{N}$.
Let $\alpha>0$, $\{\tau_k\}$ and $\{S_k\}$ be deterministic sequences.

{\it Iterate until  $k \geq K$. }
 \\  {\bf Consensus.} $\hat{v}_{i,k} := v_{i,k} $ for any $ i\in \mathcal{N}$ and
repeat the following  update by $\tau_k$ times:
$$\hat{v}_{i,k} := \sum_{j\in \mathcal{N}_i} a_{ij}\hat{v}_{j,k},~ \forall i\in \mathcal{N}.$$

{\bf Strategy Update.} For every $ i\in \mathcal{N}$:
\begin{align}
x_{i,k+1} &:=\textrm{prox}_{\alpha  r_i}  \left[ x_{i,k}-  {\alpha \over S_k} \sum_{p=1}^{S_k}
 \nabla_{x_i}\psi_i \left(x_{i,k}, n\hat{v}_{i,k};\xi_{i,k}^p\right)  \right], \label{alg-agg-strategy1}
\\ v_{i,k+1} &:= \hat{v}_{i,k}+x_{i,k+1}-x_{i,k},\label{alg-agg-average}
\end{align}
where $S_k$  is the number of sampled gradients used   at   time  $k$  and $  \xi_{i,k}^p, p=1, \cdots,S_k,$
denote the i.i.d. random  realizations of $\xi$.
\end{algorithm}

\blue{For any $x\in \mathcal{R}$, define
\begin{align}
 &   \phi_i(x)= \nabla_{x_i} f_i\left(x_i,     \sum_{i=1}^n x_i\right), \quad  F_i (x_i,\bar{x}) \triangleq \phi_i(x),
 \quad  \phi(x)\triangleq  \left(\phi_i(x) \right)_{i=1}^n  . \label{equi-phiF}
\end{align}
 We assign a map $F_i(x_i,z):\mathcal{R}_i\times \mathbb{R}^{d} \to \mathbb{R}^{d}$ to each player $i$,
where   $F_i: \mathcal{R}_i \times \sum_{j=1}^n \mathcal{R}_j \to \mathbb{R}^{d}$ is defined as in \eqref{equi-phiF}.}
Define $ e_{i,k}    \triangleq {1\over S_k} \sum_{p=1}^{S_k}  \nabla_{x_i}\psi_i(x_{i,k}, n\hat{v}_{i,k};\xi_{i,k}^p)    -\blue{F_i}(x_{i,k}, n\hat{v}_{i,k})  .$ Then  \eqref{alg-agg-strategy1} can be rewritten as:
\begin{align}
x_{i,k+1} & =\textrm{prox}_{\alpha  r_i}  \left[ x_{i,k}-\alpha \left(\blue{F_i}(x_{i,k}, n\hat{v}_{i,k}) +e_{i,k}\right)\right]. \label{alg-agg-strategy}
\end{align}
We impose the following conditions on $\mathcal{G}$,   gradient mapping, and  observation noises.
\begin{assumption}~\label{assump-agg2}
(i) The undirected graph $\mathcal{G}$ is   connected  and the adjacency matrix $A$ is  symmetric with row sums  equal to  one.  (ii) The mapping $\phi(x) $ is $\eta_{\phi}-$strongly monotone, i.e., $(\phi(x)-\phi(y))^T(x-y) \geq \eta_{\phi}  \|x-y\|^2, \quad  \forall x,y\in \mathcal{R}.$ (iii) The mapping $\phi(x)$ is $L_{\phi}$-Lipschitz continuous over $ \mathcal{R} $, i.e., $\| \phi(x)-\phi(y)\| \leq L_{\phi} \|x-y\|,\quad \forall x,y\in \mathcal{R}.$  (iv)  \blue{For each $i\in \mathcal{N} $ and  any fixed $x_i\in \mathcal{R}_i$,
 $ F_i(x_i, z)  $  is  Lipschitz continuous in $z $ over any compact set, i.e.,
for any positive constant $c_z$, there exists a constant $L_{i}$  possibly depending on $c_z$  such that
for all $  z_1,z_2\in \blue{\mathbb{R}^d}$ with $\|z_1\|\leq c_z$ and $\|z_2\|\leq c_z$: $$\|  F_i(x_i, z_1)- F_i(x_i, z_2)\| \leq L_{i} \|z_1-z_2\|.$$
\noindent (v)
For each $i\in \mathcal{N},$ there \blue{exist  positive constants $\nu_{i,1}$ and $\nu_{i,2}$} such that for any $k\geq 0,$
$\mathbb{E}[ \|e_{i,k}\|^2| \mathcal{F}_k] \leq \tfrac{\nu_{i,1}^2\|x_{i,k}\|^2+ \nu_{i,2}^2}{S_k},\quad ~ a.s.$ }
\end{assumption}

\blue{\begin{remark} We  give a simple example to show how to define the map $F_i(x_i,z)$ for each player $i$
and  any $z\in\mathbb{R}^d.$ Consider the Nash-Cournot equilibrium problem to be implemented   in Section \ref{sec:sim}:
$\min_{x_i\in X_i} f_i(x)=c_i(x_i)- d^T x_i+x_i^T   B \sum_{j=1}^n x_j.$
We have that   $\phi_i(x)=  \nabla c_i(x_i)     - d+    B  x_i+  B \sum_{j=1}^n x_j   $.
By setting  $F_i(x_i,z)=   \nabla c_i(x_i) -d+    B x_i+   B z  $,
it is obvious  that $ F_i(x_i,\bar{x})=\phi_i(x)$ and     Assumptions \ref{assump-agg2}(iv) holds.
  \end{remark}}

\subsubsection{Preliminary  Results}
Define $A(k)\triangleq A^{\tau_k}$. Then by  Assumption  \ref{assump-agg2}(i),
$A(k) $  is    symmetric  with row sums equaling  one.  We  now recall some prior results.

\begin{lemma} \label{lem-pre}  (i) \cite[Prop.~1]{nedic2009distributed}   Suppose Assumption
\ref{assump-agg2}(i) holds. Then  there exists a constant $\theta>0$ and
$\beta\in (0,1)$ such that for any $ i,j\in \mathcal{N}, $
$\Big | \big[A^k\big]_{ij}-{1\over n} \Big| \leq \theta \beta^k, \forall k\geq 1.$

\noindent (ii) \cite[Lemma 2]{koshal2016distributed}  If $y_k\triangleq  \sum_{i=1}^n v_{i,k}/n$, then
$y_k= \sum_{i=1}^n x_{i,k}/n$.
\end{lemma}

 We now introduce the transition matrices $\Phi(k,s)$ from time instance  $s$ to   $k \geq s$, defined  as $\Phi(k,k)  =A(k), ~ \Phi(k,s)  =A(k)A(k-1)\cdots A(s)$ for any $ 0\leq s<k.$  We  may  then establish an  upper bound  on the consensus error,
 \blue{and   provide a supporting lemma. The proofs of Lemmas \ref{lem6} and \ref{lem7} can be found in Appendix  \ref{ap:lem7}.}

\begin{lemma} \label{lem6}Suppose  Assumptions  \ref{assump-agg1}(i) and   \ref{assump-agg2}(i) holds.
Let Algorithm \ref{algo-aggregative} be applied to $\mathscr{P}^{\rm agg}$.  Then
the following holds for any $k\geq 0:$
\begin{align}\label{agg-bd0}
\|y_k-\hat{v}_{i,k}\|&\leq \theta D_{\mathcal{R}}   \beta^{\sum_{p=0}^k \tau_p}  + 2    \theta D_{\mathcal{R}} \sum_{s=1}^k  \beta^{\sum_{p=s}^k \tau_p}  ,
\end{align}
where $D_{\mathcal{R}}\triangleq \sum_{j=1}^n \max\limits_{x_j\in \mathcal{R}_j} \| x_{j}\|$,
and  the constants $\theta$ and $\beta$ are defined in Lemma \ref{lem-pre}(i).
\end{lemma}

\begin{lemma}\label{lem7}
Define $\tau_k  \triangleq \lceil (k+1)^{u} \rceil$  for some $u\in (0,1].$ Let $\beta\in(0,1)$.
Then  the following holds for any $k\geq 1:$
\begin{align*}
\sum_{s=1}^k  \beta^{\sum_{p=s}^k \tau_p}  &  \leq e \big(  \ln(\beta^{-1/(u+1)}) \big)^{ -1 \over u+1  } \beta^{{(k+1)^{u+1}  \over u+1}}+\beta^{{(k+1)^{u+1}  -k^{u+1}  \over u+1}}  \Big(1+  {  u+1\over k^u\ln(1/\beta)}\Big).
\end{align*}
\end{lemma}

\subsubsection{Convergence Analysis}
\begin{proposition}\label{prp2}
Suppose  Assumptions \ref{assump-agg1} and  \ref{assump-agg2} hold. Let Algorithm \ref{algo-aggregative} be applied to $\mathscr{P}^{\rm agg}$, where $\tau_k  \triangleq k+1$ and \blue{$S_k  \triangleq \left\lceil \alpha^{-2}\rho^{-(k+1)} \right \rceil$} for some $\rho\in (0,1)$. Set     $\gamma\triangleq \max\{\rho,  \beta\}$,
\blue{$\tilde{L}_{\phi} \triangleq \sqrt{ 1/2+    (1+2\alpha^2) {\bar{\nu}_1^2}  +  2L_{\phi}^2}$},
$\varrho_{\phi}\triangleq  1-2\alpha\eta_{\phi} +2\alpha^2  \tilde{L}_{\phi}^2 ,$
  $\bar{\nu}^2 =(1+2\alpha^2)(2\bar{\nu}^2_1\|x^*\|^2+\bar{\nu}^2_2)$ with  $\bar{\nu}_1 \triangleq \max_{1\leq i \leq n} \nu_{i,1}$ and $\bar{\nu}^2_2 \triangleq \sum_{i =1}^{n} \nu^2_{i,2}$,
$C_1\triangleq \theta D_{\mathcal{R}} $ and   $C_2\triangleq   2 \theta D_{\mathcal{R}}  \left(  e   \sqrt{1/\ln(\beta^{-1/2}) }+ { 2+\ln(1/\beta)\over \beta^{1/2}\ln(1/\beta)}  \right) $ with  $ D_{\mathcal{R}} \triangleq \sum_{j=1}^n \max\limits_{x_j\in \mathcal{R}_j} \| x_{j}\|$,  where the constants $\theta$ and $\beta$ are given in Lemma \ref{lem-pre}(i).
  Then     for any $ k\geq 0$:
\begin{align}\label{agg-prp1}
&\mathbb{E}[\|x_{k+1}-x^*\|^2  ]  \leq \varrho_{\phi} \mathbb{E}[\|x_{k }-x^*\|^2  ]  +C_3\gamma^{k+1}, \quad \\
\label{prp-defc}
\mbox{ where }
C_3 & \triangleq  \alpha^2  \bar{\nu}^2+ 4\alpha n  D_{\mathcal{R}}  \left(C_1  + C_2 \right)  \sum_{i=1}^nL_i   +  4\alpha^2n^2\beta   \left(C_1^2 + C_2 ^2 \right)     \sum_{i=1}^nL_i^2.
\end{align}
\end{proposition}
{\bf Proof.}The proof can be found in Appendix  \ref{ap:prp2}.\hfill $\Box$

Based on    the recursion \eqref{agg-prp1} in  Proposition ~\ref{prp2}  and by      using Lemma \ref{lem-recur},  we obtain the  linear    convergence of Algorithm \ref{algo-aggregative} with geometrically increasing sample-sizes and communication rounds increasing at a linear rate given by $\tau_k = k+1$.
\begin{theorem}[{\bf Linear rate of convergence}] \label{agg-thm1}
Suppose  Assumptions \ref{assump-agg1} and  \ref{assump-agg2} hold. Let Algorithm \ref{algo-aggregative} be applied to $\mathscr{P}^{\rm agg}$, where $\tau_k \triangleq  k+1$,  \blue{$S_k  \triangleq \left\lceil \alpha^{-2}\rho^{-(k+1)} \right \rceil$} for some $\rho \in (0,1)$, and   $\mathbb{E}[ \|x_0-x^*\|^2] \leq C$.
\blue{Let  $\gamma,$ $\varrho_{\phi} $, $\tilde{L}_{\phi}$, and $C_3$  be defined in  Proposition ~\ref{prp2}. Suppose  $\alpha \in \left (0,\eta_{\phi}/ \tilde{L}_{\phi}^2 \right)$.}  Then  $\varrho_{\phi} \in (0,1) $ and    for any $ k\geq 0$:

\noindent (i) If $\gamma\neq \varrho_{\phi},$ then $\mathbb{E}[ \|x_{k }-x^*\|^2]\leq
   {\left(C+{C_3 \over  \max\{\varrho_{\phi}/\gamma, \gamma/\varrho_{\phi}\}-1} \right)} \max \{\varrho_{\phi},\gamma\}^k $.

\noindent (ii) If $\gamma=\varrho_{\phi} ,$ then for  any $\tilde{\varrho}_{\phi}\in(\varrho_{\phi},1),$ $\mathbb{E}[ \|x_{k }-x^*\|^2]\leq \left(C+{C_3 \over \ln((\tilde{\varrho}_{\phi}/\varrho_{\phi})^e)} \right)\tilde{\varrho}_{\phi}^k  $.
\end{theorem}

Similarly to Theorem \ref{thm2},  we may derive bounds on the iteration,
oracle, communication complexity  (no. of communication rounds) to compute an $\epsilon$-NE.
\begin{theorem} \label{agg-thm2}
\blue{Suppose the conditions in Theorem \ref{agg-thm1} hold.} Then the  iteration, communication, and oracle
complexity   to obtain an $\epsilon-$NE  are  {respectively} bounded by
$K(\epsilon)$,  ${K(\epsilon) (K(\epsilon)+1)\over 2},$ and $M(\epsilon)$, where $K(\epsilon)$ and $M(\epsilon)$ are defined as follows.
\begin{align}\label{agg-rate-prox}
K(\epsilon) \triangleq \begin{cases}    {1\over \ln \left({1/\varrho_{\phi}} \right)}  \ln \left(
\big(C+{C_3 \over  \varrho_{\phi}/\gamma -1} \big) \epsilon^{-1}\right)
& {~\rm if ~} \gamma<\varrho_{\phi}<1,\\
   {1 \over \ln \left(1/\tilde{\varrho}_{\phi} \right)}
\ln\Big( \big(C+{C_3 \over \ln((\tilde{\varrho}_{\phi}/\varrho_{\phi})^e)} \big) \epsilon^{-1}   \Big) & ~{\rm if ~} \gamma=\varrho_{\phi},\\
 {1 \over \ln \left({1/\gamma} \right)} \ln\left(
\big(C+{C_3 \over  \gamma/\varrho_{\phi}  -1} \big) \epsilon^{-1}\right) &    {\rm if ~} \varrho_{\phi}<\gamma<1,
\end{cases}
\end{align} \vspace{-0.15in}
\begin{align}\label{agg-rate-SFO}
M(\epsilon) \triangleq \begin{cases}
  { 1 \over  \rho \ln( 1/\rho)}\left(
\big(C+{C_3 \over  \varrho_{\phi}/\gamma -1} \big) \epsilon^{-1}\right)^{\tfrac{\ln(1/\rho )}{\ln(1/\varrho_{\phi})}}  +K (\epsilon)
 &  {\rm if ~} \gamma<\varrho_{\phi} <1,\\
  { 1 \over  \rho \ln( 1/\rho)}\Big( \big(C+{C_3 \over \ln((\tilde{\varrho}_{\phi}/\varrho_{\phi})^e)} \big) \epsilon^{-1}   \Big) ^{\tfrac{\ln(1/\rho)}{\ln(1/ \tilde{\varrho}_{\phi} )}}  + K (\epsilon) & {\rm if ~} \gamma=\varrho_{\phi}, \\
   { 1 \over  \rho \ln( 1/\rho)}  \big(C+{C_3 \over  \gamma/\varrho_{\phi}  -1} \big) \epsilon^{-1} +K (\epsilon)   & {\rm if ~} \varrho_{\phi}<\gamma<1.
\end{cases}
\end{align}
\end{theorem}
{\bf Proof.} {Based on the geometric   rate established in Th.~\ref{agg-thm1}, we can establish the iteration complexity ($K(\epsilon)$ defined in \eqref{agg-rate-prox}) and oracle complexity ($M(\epsilon)$ defined in \eqref{agg-rate-SFO})  in the same way as that of  Theorem \ref{thm2}.}
 Since  $\tau_k=k+1$, the   communication complexity required to obtain an $\epsilon$-NE is bounded by
$ \sum_{k=0}^{K(\epsilon)-1} \tau_k=   \sum_{k=1}^{K(\epsilon) } k={K(\epsilon)(K(\epsilon)+1)\over 2}$.\hfill $\Box$

 We now prove  that  the {\em optimal} oracle complexity $\mathcal{O}\left( {1/ \epsilon}\right)$ is
obtainable  under suitable algorithm  parameters.

\begin{cor} \label{cor2} Let Algorithm \ref{algo-aggregative} be applied to $\mathscr{P}^{\rm agg}$,
 where  $\mathbb{E}[ \|x_0-x^*\|^2] \leq C$.
 Suppose  Assumptions \ref{assump-agg1} and  \ref{assump-agg2} hold.  Set   $\alpha={\eta_{\phi} \over  2\blue{\tilde{L}_{\phi}}^2}$,  $\tau_k=k+1$, and  \blue{$S_k  \triangleq \left\lceil \alpha^{-2}\rho^{-(k+1)} \right \rceil$}  with  $\rho \triangleq  \max \left \{1-{\eta_{\phi}^2\over a \blue{\tilde{L}_{\phi}}^2},  \beta\right\} $ for $a>2$.  Then    the iteration, communication,  and oracle complexity to obtain an $\epsilon-$NE  are    $\mathcal{O} ( \ln(1/\epsilon))$,
   $\mathcal{O} ( \ln^2(1/\epsilon))$, and   $\mathcal{O} \left(1/\epsilon  \right)$, respectively.
\end{cor}
{\bf Proof.}By  $\alpha ={\eta_{\phi}\over 2 \blue{\tilde{L}_{\phi}}^2}  $, we obtain  that  $\varrho_{\phi}= 1-2\alpha\eta_{\phi} +2\alpha^2 \blue{\tilde{L}_{\phi}}^2 =1-{\eta_{\phi}^2 \over  2\blue{\tilde{L}_{\phi}}^2}.$
 Note that  $\gamma=\max\{\rho,\beta \}=\rho>\varrho_{\phi}  $ by the fact that
  $\rho \geq 1- {\eta_{\phi}^2 \over  a\blue{\tilde{L}_{\phi}}^2}> 1-{\eta_{\phi}^2 \over  2\blue{\tilde{L}_{\phi}}^2}$.
Thus, by using  \eqref{agg-rate-prox}  and \eqref{agg-rate-SFO} for  the case $\gamma >\varrho_{\phi}  $,   $\gamma=\rho$,
we obtain that  $K(\epsilon) = {1 \over \ln \left({1/\rho} \right)} \ln\left(   { \epsilon^{-1}  \over  \rho \ln( 1/\rho)}  \big(  C+{C_3 \over   \rho/\varrho_{\phi}-1}  \big)  \right)$ and $M(\epsilon)=  {   C+{C_3 \over   \rho/\varrho_{\phi}-1}   \over  \rho \ln( 1/\rho)}    \epsilon^{-1}   +K (\epsilon)   . $
\hfill $\Box$

\begin{remark} Our work is not the first to utilize increasing
communication rounds. In \cite{jakovetic2014fast}, a distributed accelerated gradient algorithm  is
employed
 in convex settings, where the rate is $\mathcal{O}(1/k^2)$ (optimal)  while the
total number of communications rounds  is $\mathcal{O}(k\ln(k))$ up to time
$k$.  {In~\cite{hamedani2017multi}, for a distributed primal-dual algorithm for a constrained strongly convex problem, a  non-asymptotic convergence rate
$\mathcal{O}(1/k^2)$ is derived, requiring $\mathcal{O}(k\ln(k))$  local communications
after $k$ steps.}  Our scheme (Algorithm \ref{algo-aggregative}) requires
$k+1$ communication rounds at iteration $k$ and a total  $\mathcal{O}(k^2)$ up
to time instance $k$   to recover the optimal  {geometric  convergence
rate} but does so in a stochastic game-theoretic regime. To the best of  our
knowledge, the optimal  communication complexity for aggregative game in
deterministic regimes  is still an open question and this paper is amongst the
first to establish the communication complexity in stochastic NEPs via stochastic gradient-based techniques.
\end{remark}

 We now   explore the performance of Algorithm \ref{algo-aggregative} with  polynomially increasing  sample sizes and communication rounds. \blue{The proof can be found in Appendix \ref{app:cor5}.}
\begin{cor}  \label{cor5}
Suppose  Assumptions \ref{assump-agg1}--\ref{assump-agg2} hold. Let
Algorithm \ref{algo-aggregative} be applied to $\mathscr{P}^{\rm agg}$,
where $\mathbb{E}[ \|x_0-x^*\|^2] \leq C$, $\tau_k \triangleq   \lceil (k+1)^{u}
\rceil$  with $u\in (0,1)$, and   \blue{$S_k  \triangleq \left\lceil \alpha^{-2}(k+1)^v  \right \rceil$} with $ v>0$.
Let  $\alpha \in \left (0,\eta_{\phi}/ \blue{\tilde{L}_{\phi}}^2
\right)$  \blue{with $\tilde{L}_{\phi}$   defined in  Proposition ~\ref{prp2}.}
  Then for all $k\geq 1$, $ \mathbb{E}[ \|x_{k}-x^*\|^2  ] =\mathcal{O}\left(
k^{-v}\right)$.  In addition,   the iteration, communication,  and oracle
complexity to obtain an $\epsilon-$NE  are   $\mathcal{O} ( (1/\epsilon)^{1/v})$, $\mathcal{O} ( (1/\epsilon)^{(u+1)/v}),$ and
$\mathcal{O} \left((1/\epsilon)^{1+1/v}  \right)$,    respectively.
\end{cor}

\section{VS-PBR Scheme and the  distributed variant}\label{sec:BR}
In this section, after providing some background in Section~\ref{vs-pbr_back}, we consider a class of stochastic NEPs where the  proximal  BR map  is contractive~\cite{FPang09}.  In Section~\ref{VS-PBR}, we canalze the rate and complexity for a  variable sample size   proximal BR scheme, where at each iteration, each player solves a
sample-average BR  problem. In Section~\ref{sec:agg2}, analogous rate and complexity statements are provided for a distributed variant of ({\bf VS-PBR}).

\subsection{Background  on proximal best-response maps}\label{vs-pbr_back}
For  any  tuple $y\in \mathbb{R}^d,$ let the proximal BR map $\widehat{x}(y)$
of  $\mathscr{P}$ be defined as follows:
\begin{align}\label{eq-prox-resp}
\widehat{x}(y) \triangleq\argmin_{x \in \mathbb{R}^d}  \left[  \sum_{i=1}^n  p_i(x_i,y_{-i}) + {\mu \over 2} \|x-y\|^2 \right]\quad
{\rm ~for ~some~} \mu>0.
\end{align}
It is  clear that the objective function is  separable in $x_i$ and \eqref{eq-prox-resp}
reduces to a set of player-specific proximal   BR problems, where    {player $i$ solves}  the following  problem:
\begin{align}\label{prox-BRi}\widehat{x}_{i}(y)    \triangleq \argmin_{x_i \in \mathbb{R}^{d_i}}
\left [  \mathbb{E}\left[ \psi_i(x_i,y_{-i};\xi ) \right]+ r_i(x_i) +{\mu \over 2} \|x_i-y_i\|^2 \right] .
\end{align}
We impose the following assumption on problem \eqref{Ngame}.
\begin{assumption}~\label{assp-compact}
For each player $i\in \mathcal{N}$,
(i) Assumption \ref{assump-agg1}(i) holds;
\\ (ii)  for every fixed $x_{-i} \in \Rscr_{-i}$, $f_i(x_i,x_{-i})$ is twice continuously differentiable
  and convex in $x_i \in \mathcal{R}_i$,
and  $\nabla_{x_i} f_i(x_i,x_{-i})$ is $L_{fi}$-Lipschitz continuous  in $x_i$, i.e.,
$\|\nabla_{x_i}  f_i(x_i, x_{-i})-\nabla_{x_i}  f_i(x_i', x_{-i})\| \leq  L_{fi} \|x_i-x_i'\|$
for all $x_i,x_i'\in\mathcal{R}_i$ ;
 (iii) for   any $x_{-i} \in  \Rscr_{-i} $ and   $ \xi_i\in\mathbb{R}^{m_i}$, $  \psi_i(x_i, x_{-i};\xi_i)$
  is differentiable  in $x_i \in  \Rscr_i$ and there exists some $\nu_i>0$:
$\mathbb{E}[\| \nabla_{x_i} f_i(x_i,x_{-i})-\nabla_{x_i} \psi_i(x_i,x_{-i};\xi_i)\|^2] \leq \nu_i^2,
 \quad \forall x \in \mathcal{R} .$
\end{assumption}

Then by \cite[Proposition 12.5]{FPang09},  $x^*$ is an NE of the game $\mathscr{P}$ if and only if $x^*$ is a fixed point of the proximal best-response map $ \widehat{x}(\bullet),$ that is, if and only if $x  ^*= \widehat{x}(x^*)$. By Assumption \ref{assp-compact},  the  second derivatives  of the   functions $f_i,~\forall i \in \mathcal{N}$ on $\mathcal{R}$ are bounded. Analogous to the avenue adopted in ~\cite{FPang09},   we may define  \begin{align}
&  \Gamma \triangleq \pmat{\tfrac{\mu}{\mu+\zeta_{1,\min}} &
	\tfrac{\zeta_{12,\max} }{\mu+\zeta_{1,\min}} & \hdots &
		\tfrac{\zeta_{1n,\max}}{\mu+\zeta_{1,\min}}\\
\tfrac{\zeta_{21,\max}}{\mu+\zeta_{2,\min}} &
	\tfrac{\mu}{\mu+\zeta_{2,\min}} & \hdots &
		\tfrac{\zeta_{2n,\max}}{\mu+\zeta_{2,\min}}\\
	\vdots & & \ddots & \\
\tfrac{\zeta_{n1,\max}}{\mu+\zeta_{n,\min}} &
	\tfrac{\zeta_{n2,\max}}{\mu+\zeta_{n,\min}} & \hdots &
		\tfrac{\mu}{\mu+\zeta_{n,\min}}}, \label{matrix-hessian}
\end{align}
where $\zeta_{i,\min} \triangleq \inf_{x \in \mathcal{R}} \lambda_{\min}  \left (\nabla^2_{x_i} f_i(x) \right)  \mbox{ ~and ~} \zeta_{ij,\max} \triangleq \sup_{x \in \mathcal{R}}  \| \nabla^2_{x_ix_j} f_i(x) \|  ~  \forall j \neq i.$ Then by \cite[Theorem 4]{pang2017two}, we may obtain the following relation:
\begin{align}
  \pmat{\|\wh{x}_1(y') -\wh{x}_1(y)\| \\
				\vdots \\
\|\wh{x}_n(y') - \wh{x}_n(y)\|}  \leq   \Gamma   \pmat{\|y_1'- y_1\| \\
				\vdots \\
\|y'_n - y_n\|}  .\label{cont-prox-best-resp}
\end{align}
If the spectral radius    $ \rho(\Gamma)< 1$,  then  the  proximal best-response map is contractive w.r.t. some monotonic norm.
Sufficient conditions for the contractive property of the  proximal BR   map $ \wh{x}(\bullet)$ can be found in \cite{FPang09,pang2017two}.

\subsection{Variable sample-size proximal  BR schemes}\label{VS-PBR}
Suppose at iteration $k,$ we have   $S_k$  \blue{i.i.d.} realizations  $\xi_{i,k}^1, \cdots ,\xi_{i,k}^{S_k}$ of the random vector $\xi_i.$ For any $x_i\in \mathcal{R}_i,$ we approximate $f_i(x_i,y_{-i,k})$ by its sample-average  ${1\over S_k}\sum_{p=1}^{S_k} \psi_i(x_i,y_{-i,k}; \xi_{i,k}^p)$
and  solve the  sample-average  BR problem \eqref{SAA}, leading to
Algorithm~\ref{inexact-prox-br}.
\begin{algorithm} [htbp]
\caption{{Variable-size proximal  best-response  scheme}} \label{inexact-sbr-cont}
 Set $k:=0$.  Given $K > 0$, let  $ y_{i,0}=x_{i,0} \in X_i$ for $i = 1, \hdots, n$.
\begin{enumerate}
\item[(1)] For $i = 1, \hdots, n$,  player  $i$ updates estimate $x_{i,k+1}$ as
 \begin{equation}\label{SAA}
\begin{split}
x_{i,k+1}& =  \argmin _{x_i \in  \Real^{d_i}}  \Bigg[ {1\over S_k}\sum_{p=1}^{S_k} \psi_i(x_i,y_{-i,k}; \xi_{i,k}^p)
 +r_i(x_i)  + {\mu \over 2} \|x_i-y_{i,k}\|^2 \Bigg].
\end{split}
\end{equation}
\item[(2)] For $i = 1, \hdots, n$,  $y_{i,k+1} :=x_{i,k+1}$;
\item[(3)]  $k:=k+1$ and return to (1)  if $k < K$.
\end{enumerate}
\label{inexact-prox-br}
\end{algorithm}

Denote by   $\varepsilon_{i,k+1} \triangleq  x_{i,k+1}- \wh{x}_{i}(y_k)  $   the inexactness associated with the approximate  proximal BR solution. We   now give the bound of   $ \mathbb{E}\left[\|\varepsilon_{i,k+1}\|^2 \right]   $ regarding the   inexactness sequence  in the following lemma. \blue{The proof is given in Appendix \ref{ap:lem4}.}
\begin{lemma}\label{lem4} Suppose Assumption   \ref{assp-compact}   holds.
Let   Algorithm \ref{inexact-sbr-cont} be applied to $\mathscr{P}$.
 Define $C_{i,b} \triangleq  {  \mu \over \mu^2+L_{fi}^2}\Big(1-L_{fi}/ \sqrt{   \mu^2+L_{fi}^2} \Big)^{-1}$.
 Then   for each $i=1,\cdots,n:$
 $\mathbb{E}[ \|x_{i,k+1}-\wh{x}_{i}(y_k) \|^2]\leq \tfrac{\nu_i^2  C_{i,b}^2}{S_k}$ for all $k\geq 0.$
\end{lemma}

Based on this lemma,  we obtain a  linear   rate of convergence with a suitably selected sample size $S_k$.
\blue{The proof of Proposition \ref{prp-linear}  can be found in Appendix \ref{ap:prp-linear}.}

\begin{proposition}[{\bf  Linear rate of convergence}] \label{prp-linear}~
Suppose  Assumption \ref{assp-compact}  holds  and     $a  \triangleq \| \Gamma\|<1$, where $\Gamma$ is defined in \eqref{matrix-hessian}.  Define $C_{ns}=\max_i\nu_i^2  C_{i,b}^2  $ with  $C_{i,b} \triangleq  {  \mu \over \mu^2+L_{fi}^2}\Big(1-L_{fi}/ \sqrt{   \mu^2+L_{fi}^2} \Big)^{-1}.$  Let   Algorithm \ref{inexact-sbr-cont} be applied to $\mathscr{P}$, where $\mathbb{E}[\|x_{0}-x ^*\|^2] \leq C $  and  $S_k=\left\lceil {  C_{ns}\over\eta^{2(k+1)}}\right\rceil$ for some $\eta\in (0,1)  $. Then  the following  hold.

\noindent (i) If $a \neq \eta,$ then $\mathbb{E}[ \|x_{k }-x^*\|^2]\leq
\left( \sqrt{C}+ {\sqrt{n}\over  \max\{a/\eta,\eta/a\}-1} \right)^2  \max \{a, \eta \}^{2k}.$

\noindent (ii) If $\eta=a,$ then for   any  $ \tilde{a} \in (a,1)$, $
\mathbb{E}[\|x_{k}-x ^*\|^2]  \leq  \left(   \sqrt{C }+ \sqrt{n} / \ln((\tilde{a}/a)^e) \right)^2  \tilde{a} ^{2k} .$
\end{proposition}

Note that  $x_{i,k+1}$ defined by \eqref{SAA} requires solving a  deterministic optimization.
In the following, we establish the iteration complexity  (no. of deterministic optimization problems  solved)  and oracle complexity  to obtain an   $\epsilon-$NE.

\begin{theorem} \label{thm4}
Suppose  Assumption \ref{assp-compact}  holds,  $a  \triangleq \|
	\Gamma\|<1$, $C_{i,b} \triangleq  {  \mu \over
		\mu^2+L_{fi}^2}\big(1-L_{fi}/ \sqrt{   \mu^2+L_{fi}^2} \big)^{-1}$, and  $C_{ns} \triangleq \max_i\nu_i^2  C_{i,b}^2 $. Let   Algorithm \ref{inexact-sbr-cont} be applied to $\mathscr{P}$,  where $\mathbb{E}[\|x_{ 0}-x ^*\|^2] \leq C $  and  $S_k=\left\lceil {  C_{ns}\over  \eta^{2(k+1)}}\right\rceil$   for some $\eta\in (0,1) $.   Let $ \tilde{a} \in (a,1)$, and $D= 1/ \ln(( \tilde{a} /a)^e) $. Then
the  iteration and oracle complexity to obtain  an  $\epsilon-$NE
 are bounded by $K_b(\epsilon)$ and $M_b(\epsilon)$
respectively,  each of which is defined as follows.
\begin{align}\label{br-rate-prox}
K_b(\epsilon) & \triangleq \begin{cases}  &  {1\over \ln \left( 1/a \right)} \ln \left (  {\sqrt{C}+  \eta \sqrt{n}/(a- \eta )   \over  \sqrt{\epsilon } } \right) \quad \quad {\rm if ~} \eta<a,\\
&  	\tfrac{1}{\ln(1/  \tilde{a} )}  \ln\left ( { \sqrt{C }+ \sqrt{n} D  \over   \sqrt{\epsilon }} \right) \qquad  \qquad ~{\rm if ~} \eta=a,\\
&  {1\over \ln \left( 1/\eta \right)} \ln\left (  {   \sqrt{C}+  a\sqrt{n}/(  \eta -a)   \over  \sqrt{ \epsilon}} \right) \quad \quad {\rm if ~} \eta>a ,
\end{cases}\\
\label{br-rate-SFO}
M_b(\epsilon) & \triangleq \begin{cases}  &
  {C_{ns}\over\eta^2 \ln(1/\eta^{2 })} \left ({ \left( \sqrt{C}+  \eta \sqrt{n}/(a- \eta )  \right)^2  \over  \epsilon } \right)^{\ln (1/\eta ) \over\ln \left(1/a \right)}   +K (\epsilon) \quad {\rm if ~} \eta<a,\\&
{C_{ns}\over a^2 \ln(1/a^{2 })} \left ( { (\sqrt{C }+ \sqrt{n} D)^2 \over   \epsilon } \right)^{\ln (1/a ) \over\ln(1/  \tilde{a} ) }+K (\epsilon)\qquad  \quad{\rm ~ if ~}  \eta=a, \\
&    {C_{ns}\over  \eta^2 \ln(1/\eta^{2 })} \left (  {  \left( \sqrt{C}+  a\sqrt{n}/(  \eta -a)   \right)^2 \over   \epsilon} \right)   +K (\epsilon)\quad  \qquad ~~{\rm if ~} \eta>a.
\end{cases}
\end{align}
\end{theorem}
{\bf Proof.}The proof can be found in Appendix \ref{ap:thm4}.
\hfill $\Box$

  The above theorem establishes that when the number of scenarios  increases at a geometric rate,  the iteration  and oracle complexity to achieve an $\epsilon-$NE    are   respectively  $\mathcal{O} (\ln(1/\epsilon))$  and  $\mathcal{O} ( (1/\epsilon)^{1+\delta})$, where  $\delta=0$ when $\eta\in (a,1),$  $\delta= \tfrac{\ln(a/\eta)}{\ln(1/a)}~{\rm when ~} \eta<a<1$, and $\delta=\tfrac{\ln( \tilde{a}  /a)}{\ln(1/ \tilde{a} )}$ when $\eta=a.$
Similarly to the  discussions in Section \ref{VSPGR-ploy}, we now
establish  the rate and complexity properties of  Algorithm \ref{inexact-sbr-cont}   with polynomially increasing
 sample-sizes. \blue{The proof of the following Corollary \ref{poly-lem2}  can be found in Appendix \ref{ap:poly-lem2}.  }
\begin{cor}\label{poly-lem2}
 Let   Algorithm \ref{inexact-sbr-cont} be applied to $\mathscr{P}$,  where $\mathbb{E}[\|x_{ 0}-x ^*\|^2] \leq C $  and   $S_k \triangleq  \lceil (k+1)^{ v}\rceil $ for some $v>0$.
Suppose  Assumption \ref{assp-compact}  holds  and     $a  \triangleq \| \Gamma\|<1$, where $\Gamma$ is defined in \eqref{matrix-hessian}. Then   we obtain the polynomial rate of convergence  $ \mathbb{E}[ \|x_{k}-x^*\|^2  ] =\mathcal{O}\left( k^{- v}\right)$  and establish that   the  iteration and oracle complexity bounds to obtain an $\epsilon-$NE are
{{$\mathcal{O} (v (1/\epsilon)^{1/v})$ and   $\mathcal{O} \left(e^{v} v^v(1/\epsilon)^{1+1/v}  \right)$}}, respectively.
\end{cor}

\subsection{Distributed VS-PBR for Aggregative Games}\label{sec:agg2}
We   propose a distributed VS-PBR scheme   to solve   the aggregative game $\mathscr{P}^{\rm agg}$ formulated in Section \ref{sec:agg}. Suppose at iteration $k,$  each player updates its belief of  the aggregate by multiple
consensus steps,  utilizes   $S_k$ realizations to  approximate the cost and solve the  sample-average proximal BR problem \eqref{alg-agg-strategy2}. We then obtain   Algorithm  \ref{algo-aggregative2}.

\begin{algorithm}[htbp]
\caption{Distributed  VS-PBR for Aggregative  Stochastic NEPs} \label{algo-aggregative2}
 {\it Initialize:} Set $k =0$   and  $ v_{ i,0} =x_{ i,0 } \in \mathcal{R}_i$  for  each $i \in\mathcal{N}$.
Let $\alpha>0$ and $\{\tau_k \}$ be a deterministic sequence.

{\it Iterate until  $k > K$}

 {\bf Consensus.} $\hat{v}_{i,k} := v_{i,k} $ for each $ i\in \mathcal{N}$ and
repeat $\tau_k$ times:
$$\hat{v}_{i,k} := \sum_{j\in \mathcal{N}_i} a_{ij}\hat{v}_{j,k}, \quad \forall i\in \mathcal{N}.$$

 {\bf Strategy Update.}  For  every $ i\in \mathcal{N}$,
  \begin{align}
x_{i,k+1}& =  \argmin _{x_i \in  \Real^{d}}  \Big[ {1\over S_k}\sum_{p=1}^{S_k} \psi_i(x_i,\blue{x_i-x_{i,k}+}n\hat{v}_{i,k}; \xi_{i,k}^p)
 +r_i(x_i)  + {\mu \over 2} \|x_i-x_{i,k}\|^2 \Big] , \label{alg-agg-strategy2}
\\ v_{i,k+1} &:= \hat{v}_{i,k}+x_{i,k+1}-x_{i,k},\label{alg-agg-average2}
\end{align}
where $  \xi_{i,k}^p, p=1, \cdots,S_k,$ denote the i.i.d. random  realizations of $\xi_i$.
\end{algorithm}

\subsubsection{Rate Analysis}

We impose the following assumptions on the $\mathscr{P}^{\rm agg}.$
{
\begin{assumption}\label{ass-agg-noise}
For each $i\in \mathcal{N} $,   (i) Assumption \ref{assump-agg1}(i) holds;
(ii)  for any $y\in \mathbb{R}^d$, $f_i(x_i,x_i+y)$  is twice continuously differentiable  and convex in $x_i \in \mathcal{R}_i$;
(iii) for any  $y\in \mathbb{R}^d$, there exists a constant $L_a>0$ such that for any $x_i, x_i' \in \mathcal{R}_i:$
$\| g_i(x_i,  y )- g_i(x_i',  y )\| \leq L_a \| x_i-x_i' \| ,$
where $g_i(x_i,y)\triangleq \nabla_{x_i}f_i(x_i,x_i+ y )$;
 (iv) for any  $x_i\in \mathcal{R}_i,$  $  g_i(x_i,  y) $  is $L_i$-Lipschitz continuous in $y\in \mathbb{R}^d$, i.e.,
 \[\|  g_i(x_i,  y_1)- g_i(x_i,  y_2)\| \leq L_{gi} \|y_1-y_2\|,  \quad \forall y_1,y_2\in \mathbb{R}^d;\]
(v) for  any $y\in \mathbb{R}^{d} $, $  \psi_i(x_i, x_i+y;\xi)$ is differentiable  in $x_i\in \Rscr_i$  such that  for some $\nu_i>0$,
   $\mathbb{E}[\| \nabla_{x_i} f_i(x_i,x_i+y)-\nabla_{x_i} \psi_i(x_i,x_i+y;\xi)\|^2] \leq \nu_i^2,
    \quad \forall x_i \in \mathcal{R}_i, y\in \mathbb{R}^{d}.$
\end{assumption}

\begin{remark} We  give a simple example to illustrate  Assumptions   \ref{ass-agg-noise}(i)-(iv) hold.
 Consider the Nash-Cournot equilibrium problem  where player $i$ aims to solve
$\min_{x_i\in X_i} f_i(x_i,\bar{x})=c_i(x_i)- d^T x_i+x_i^T   B \sum_{j=1}^n x_j .$ Then
$f_i(x_i,x_i+y)=c_i(x_i)- d^T x_i+x_i^T   B  (x_i+y),$ and    $g_i(x_i,y)=  \nabla c_i(x_i)     - d+    2B  x_i+  B y  $.
It is obvious   that  Assumption  \ref{ass-agg-noise}(iv)  holds, while  Assumption   \ref{ass-agg-noise}(iii) holds when
 each $\nabla c_i(x_i)$ is Lipschitz continuous  in $x_i \in \mathcal{R}_i$.
  \end{remark}

For any $x_i \in \mathcal{R}_i $ and any $z  \in \mathbb{R}^d,$ we define a proximal BR map
\blue{to $\mathscr{P}^{\rm agg}$.}
\begin{align}\label{def-map-T}
T_i(y_i,z)    \triangleq  \argmin_{x_i \in \mathbb{R}^d}
\left [  f_i(x_i, x_i +z -y_i )+ r_i(x_i) +{\mu \over 2} \|x_i-y_i\|^2 \right],\quad \mu>0 .
\end{align}
Then $T_i(x_i,z) $ is uniquely defined  by   Assumptions  \ref{ass-agg-noise}(i) and  \ref{ass-agg-noise}(ii).
The Lipschitz continuity   of  $T_i(y_i,z)$ is proved in the next lemma ({Proof in Appendix \ref{ap:Lip-Tmap}}).

\begin{lemma}\label{Lip-Tmap}~
Suppose   Assumptions   \ref{ass-agg-noise}(i)-(iv) hold.
 Define $L_t \triangleq {\mu \max_{i\in \mathcal{N}}L_{gi} \over  \mu^2+L_a^2 } \big(1-L_a /\sqrt{ \mu^2+L_a^2}\big)^{-1}.$
There  for each $i\in \mathcal{N}$ and any $y_i\in \mathcal{R}_i$, the following holds:
	\begin{align}\label{T-Lipschitz}
\|T_i(y_i,z_1)  -T_i(y_i,z_2)  \|  \leq L_t\|z_1-z_2\| ,
\quad \forall  z_1,z_2  \in \mathbb{R}^d.
\end{align}
\end{lemma}

Similarly to Lemma \ref{lem4}, we can obtain the following result.

\begin{lemma}\label{agg-variance}  Let   Algorithm \ref{algo-aggregative2} be applied to $\mathscr{P}^{\rm agg}$.
 Suppose   Assumption      \ref{ass-agg-noise}(i), (ii), (iii), and (v)  hold.
  Define $\varepsilon_{i,k+1}\triangleq x_{i,k+1}-T_i(x_{i,k}, n\hat{v}_{i,k}) $.  Then    for each $i=1,\cdots,n,$
$ \mathbb{E}[ \|\varepsilon_{i,k+1} \|^2]\leq {\nu_i^2  C_r^2 \over S_k} $ with  $C_r \triangleq  {  \mu \over \mu^2+L_a^2}\left(1-L_a/ \sqrt{   \mu^2+L_a^2} \right)^{-1}$.
\end{lemma}
This allows for obtaining a  linear   rate with a suitably selected sample size $S_k$ and  the number of communication rounds $\tau_k$. The proof of can be found in Appendix \ref{app:prp3}.
\begin{proposition}\label{prp3}
Suppose   $a=\|\Gamma\|<1$ with $\Gamma$   defined by \eqref{matrix-hessian}, Assumptions \ref{assump-agg2}(i) and   \ref{ass-agg-noise}  hold.     Let Algorithm \ref{algo-aggregative2} be applied to $\mathscr{P}^{\rm agg}$, where  $\mathbb{E}\left [ \|x_0-x^*\|^2\right] \leq C,$ $\tau_k=k+1$,
 and $S_k=\left\lceil { C_r^2  \max_i\nu_i^2\over  \eta^{2(k+1)}}\right\rceil$  for some $\eta\in (0,1)$.
  Define $C_4\triangleq \sqrt{n}  + n^{3\over 2} L_t \left(  C_1  + C_2 \right) $, and  $\gamma\triangleq \max\{\eta,  \beta\}, $
where $C_1$ and $C_2$ are defined in Proposition \ref{prp2}, and $\beta$ is  given in Lemma \ref{lem-pre}(i).
  Then  the following hold for any  $k\geq 0$:

\noindent (i) If $a \neq \gamma,$ then $\mathbb{E}[ \|x_{k }-x^*\|^2]\leq
Q^2 \max \{a, \gamma \}^{2k},$ where $Q\triangleq   \sqrt{C}+ {C_4\over \max\{a/\gamma,\gamma/a\}-1}  .$

\noindent (ii) If $\gamma=a,$ then for $ \tilde{a} \in (a,1)$, $ \mathbb{E}[\|x_{k}-x ^*\|^2]  \leq  \left(\widetilde{Q}\right)^2 \tilde{a}^{2k} ,$ where $\widetilde{Q}   \triangleq    \sqrt{C }+\tfrac{C_4}{\ln((\tilde{a}/a)^e}.$

\end{proposition}

\subsubsection{Iteration and Oracle Complexity}

We now  establish the iteration, oracle, and communication  complexity of
Algorithm~\ref{algo-aggregative2} to achieve  an $\epsilon-$NE.  {The
derivation of Theorem~\ref{thm-iter-comp} is similar to Theorem
\ref{thm4} and its proof is omitted.}

\begin{theorem} \label{thm-iter-comp}   Let Algorithm \ref{algo-aggregative2} be applied to $\mathscr{P}^{\rm agg}$, where
requirements of Proposition \ref{prp3}  hold. Suppose  $a \neq \gamma $.
  Then the number of  optimization problems solved  by player $i$ to obtain an  $\epsilon-$NE  is  bounded by  $  K (\epsilon) \triangleq \left\lceil \tfrac{\ln(Q/\sqrt{\epsilon})}{\ln(1/ \max \{a, \gamma \} )}\right\rceil,$ while  the communication and oracle  complexity  to obtain an  $\epsilon-$NE are ${K(\epsilon)(K(\epsilon)+1)\over 2}  $ and  $\mathcal{O}\left(
	  (1/ \epsilon)^{{\ln (1/\eta) \over\ln \left({1/ \max \{a,\gamma\}} \right)}} \right)$.
\end{theorem}
%
%
%
%
%

\begin{remark}
(i)   From Theorem \ref{thm-iter-comp} we conclude that the  iteration complexity in terms of the deterministic optimization solver, oracle complexity in terms of sampled gradient, and communication complexity to compute an
$\epsilon$-NE are    $\mathcal{O} (\ln(1/ \epsilon ))$, $\mathcal{O} \left((1/\epsilon)^{1+\delta}\right)$ with $\delta\geq 0$, and $\mathcal{O} \left(\ln^2(1/ \epsilon )\right)$, respectively.
\\
(ii)  Let  $\eta\in \left(\max\{a, \beta \},1\right)$. Then
we obtain  the {\em optimal} oracle complexity $\mathcal{O} \left(1/\epsilon\right)$. \\
(iii) Similarly to Corollary \ref{cor5},  when $\tau_k= \lceil (k+1)^{u} \rceil$ and  $S_k=\left\lceil (k+1)^v \right \rceil$ for  $u\in (0,1)$ and $ v>0$,  we  can obtain the  convergence rate $ \mathbb{E}[ \|x_k-x^*\|^2  ] =\mathcal{O}\left( k^{-v}\right)$  and  establish that   the iteration, communication,  and oracle complexity to obtain an $\epsilon-$NE  are    $\mathcal{O} ( (1/\epsilon)^{1/v})$, $\mathcal{O} ( (1/\epsilon)^{(u+1)/v}),$ and   $\mathcal{O} \left((1/\epsilon)^{1+1/v}  \right)$,    respectively. In fact, as done earlier, we may clarify the dependence of the constants on $u$ and $v$.
\end{remark}

\section{Numerical Simulations}\label{sec:sim}

In this section, we  empirically  validate the performance of   the proposed algorithms   on   the    networked  Nash-Cournot equilibrium problem~\cite{kannan12distributed,yu2017distributed,koshal2016distributed},
where  firms compete in quantity produced. It  is a classical example of an  aggregative game  where the inverse-demand function depends on the sum of production by   all   firms. We assume that  there are $n$  firms, regarded as the    set  of players  $\mathcal{N}=\{1,\dots,n\}$
  competing over   $L$  spatially distributed markets (nodes) denoted by $\mathcal{L}=\{1,\cdots,L\}$.
 For any $i$, the $i$th firm needs to determine a continuous-valued
  nonnegative quantity of products to be produced
and delivered to the markets, which is defined as $x_i=(x_{i}^1, \cdots, x_{i}^L) \in \mathbb{R}^L$, where  $x_{i}^l $
denotes  the sales of firm $i$ at  the market $l$.
 Furthermore, the $i$th firm is characterized by   a random linear  production cost
function $c_i(x_i;\xi_i)=(c_i+\xi_i)  \sum_{l=1}^L x_{i}^l $  for some   parameter  $c_i>0 $ where  $\xi_i$ is a mean-zero random variable. We further assume that  the price  $p_l$  of products sold in     market $l\in \mathcal{L}$  is determined by the linear inverse demand (or price) function  corrupted by noise $p_l(\bar{x}_l;\zeta_l)=d_l+\zeta_l-  b_l \bar{x}_l,$
where $\bar{x}_l=\sum_{i=1}^n x_{i}^l$ is  the total sales  of products at the market $l,$
    the positive parameter $d_l$ indicates the price when
the production of the good is zero,  the positive parameter $b_l$ represents the slope of the inverse demand
function,  and the random disturbance $\zeta_l$ is zero-mean.
 Consequently,   firm $i $ has an expectation-valued  cost function defined  as
 $ f_i(x)=\mathbb{E} \big[ c_i(x_i;\xi_i)-\sum_{l =1}^L p_l(\bar{x}_l;\zeta_l)x_i^l  \big].$
 Suppose    firm  $i \in \mathcal{N}$  has finite  production capacity $ X _i=\{x_i\in \mathbb{R}^{L}:  x_i \geq 0,   x_i^l \leq \textrm{cap}_{il}\}.$  Then the objective of firm $i$ is to find a feasible strategy that optimizes its cost,
  i.e.,  $\min_{x_i \in X_i} f_i(x_i,x_{-i}).$

{\bf Numerical settings.} In the numerical study, we consider a network with $n $ firms
and $L = 10$ markets with   the   parameters in the payoffs   set  as  $d_l\sim U(40,50), b_l\sim U(1,2), c_i\sim U(3,5) $ for all $i\in \mathcal{N}$ and $ l\in \mathcal{L},  $
where $U(\underline{u},\bar{u})$ denotes the   uniform distribution over an interval $[\underline{u},\bar{u}]$ with $\underline{u}<\bar{u}$. In  the stochastic
settings, the random variables  are assumed to be $\xi_i \sim U(-c_i/5, c_i/5 ),    \zeta_l \sim   U(-a_l/5, a_l/5 ),$ respectively. We further set  $\textrm{cap}_{il}=2 $ for each $ i\in \mathcal{N} $ and $ l\in \mathcal{L}$.

{\em We first validate the performance of the{ \bf distributed VS-PGR}  scheme.}
We  consider   four kinds  of undirected connected  graphs: (i) Cycle graph, which consists of a single cycle and every   node has exactly two edges incident with it; (ii) Star graph, where there is a center node connecting to all every other node; (iii) Erd\H{o}s--R\'{e}nyi graph, which   is constructed by connecting nodes randomly and each edge is included in the graph with   probability $2/n$  independent from every other edge;   (iv) Complete graph,
 each node  has an edge connecting it to every other node. We define a doubly stochastic matrix $A=[a_{ij}]_{i,j=1}^n$
with  $a_{ii} =1-{d(i)-1\over d_{\max}},  a_{ij}=  {1\over d_{\max}}  {\rm~if~} (j,i) \in \mathcal{E} {\rm ~and ~} j\neq i,a_{ij}=0,$ otherwise,
where $d(i)=| \mathcal{N}_i|$ denotes the number of neighbors of player $i $ and $ d_{\max}=\max_{i\in \mathcal{N}} d(i).$
We implement Algorithm \ref{algo-aggregative} with  $\tau_k=\lceil \log(k) \rceil  $ and  $S_k=\big\lceil \rho^{-(k+1)} \big \rceil$, and terminate it when the total number of samples utilized  reached $10^6$ and  report the  empirical    error  of $ { \mathbb{E}[ \| x_{k}-x^*\|] \over   \|  x^*\|  }$ by averaging across   50 sample paths.
The simulation results are demonstrated in Table \ref{TAB1}. As expected,  Algorithm \ref{algo-aggregative}  with  complete graph has fastest convergence rate  and  the empirical    error  at the termination increases with the size of the network.
It  also indicates that the constant step-size  should not be taken  too large, otherwise it  might  lead to non-convergence,
see e.g.   $\alpha=0.02$ in   the case  $n=50$.  In Figure \ref{fig_tra}, we  further display  trajectories  of   the iterates generated by  the centralized  VS-PGR   over the complete graph  and its distributed variant over the ER graph.  Though  the centralized   has faster  convergence rate
than its distributed variant,   it requires much more  rival information (or communications).

  \begin{table}[htb]
\centering
 \begin{tabular}{|c|c|c|c|c|c|c| c|}
 \hline   $\alpha$ & $\rho$ & $n$ &   Cycle  &
  Star  &   E-R  &Complete   \\ \hline
 {\multirow{4}{*}{ $ 0.01$}}      &  {\multirow{2}{*}{$ 0.98 $ }}  &$ 20$  &  3.16e-04& 1.15e-01 &7.5e-02 &2.96e-04       \\  \cline{3-7}
& &  $ 50$       &  1.55e-01&4.73e-01 &3.68e-01&1.1e-03 \\  \cline{2-7}
 &   {\multirow{2}{*}{$ 0.985$}}   & $ 20$  &  1e-03& 1.15e-01 &7.47e-02  &2.36e-04       \\ \cline{3-7}
&& $ 50$   &  1.49e-01&4.73e-01 &3.67e-01&4.78e-04   \\
 \hline
 {\multirow{4}{*}{ $ 0.02$}}      &  {\multirow{2}{*}{$ 0.98 $ }}  &$ 0$        &  9.07e-04& 1.15e-01 &7.47e-02 &2.96e-04 \\ \cline{3-7}
&& $ 50$   &  2.67e-01&5.27e-01 &4.37e-01&2.07e-01
 \\  \cline{2-7}
 &  {\multirow{2}{*}{$ 0.985 $ }}  &$ 20$   &  1.2e-03& 1.15e-01 &7.47e-02 &3.65e-04   \\  \cline{3-7}
&& $ 50$ & 2.67e-01&5.27e-01 &4.37e-01&2.07e-01     \\
 \hline
 \end{tabular}
      \caption{Empirical    Error} \label{TAB1} \vspace{-0.2in}
\end{table}

We now investigate how does the network structure influence the convergence properties.
Set   $n=20$ and run  Algorithm \ref{algo-aggregative} over   the cycle,  star, and   Erd\H{o}s--R\'{e}nyi   graphs with  $\tau_k=k+1$, $\alpha=0.01,$ and  $S_k=\left\lceil \beta^{-(k+1) } \right \rceil$, where  the network connectivity parameter $\beta$ are   respectively $ 0.967, 0.95, 0.986$ for cycle, star and ER graphs.
The simulation results are demonstrated in Figure \ref{fig} with the left and right figure respectively  displaying  the rate of convergence and oracle complexity.  It is shown  that  the star graph with $\beta=0.95$ has the fastest convergence rate while  the ER graph with $ \beta=0.986$  has the  slowest convergence rate,
this is consistent with Theorem \ref{agg-thm1} that smaller $\beta$ may lead to faster rate of convergence (since $\gamma=\beta$). It is also worth noting that  for obtaining an  $\epsilon$-NE,  the ER graph  requires the smallest number of samples, while    the star graph has the worst oracle complexity.  These findings support the
theoretic  results   in  Theorem \ref{agg-thm2}   that  larger $\beta$ may lead to    better oracle complexity.

\begin{figure}[htbp]
\centering
\begin{minipage}[t]{0.33\textwidth}
\centering
   \includegraphics[width=2.2in]{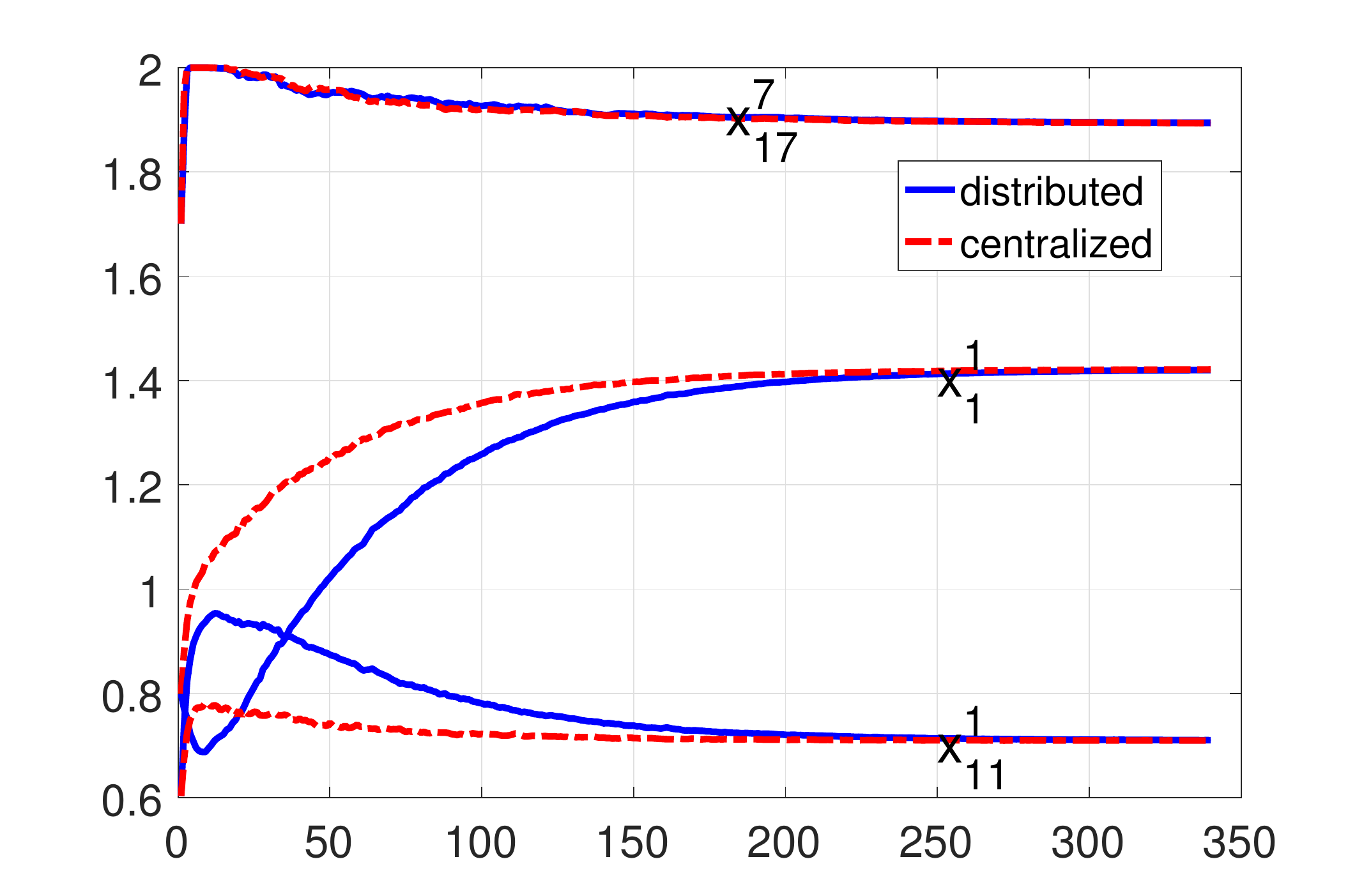} \\
 \caption{Trajectories}\label{fig_tra}
\end{minipage}
\begin{minipage}[t]{0.66\textwidth}
\centering
 \includegraphics[width=4.4in]{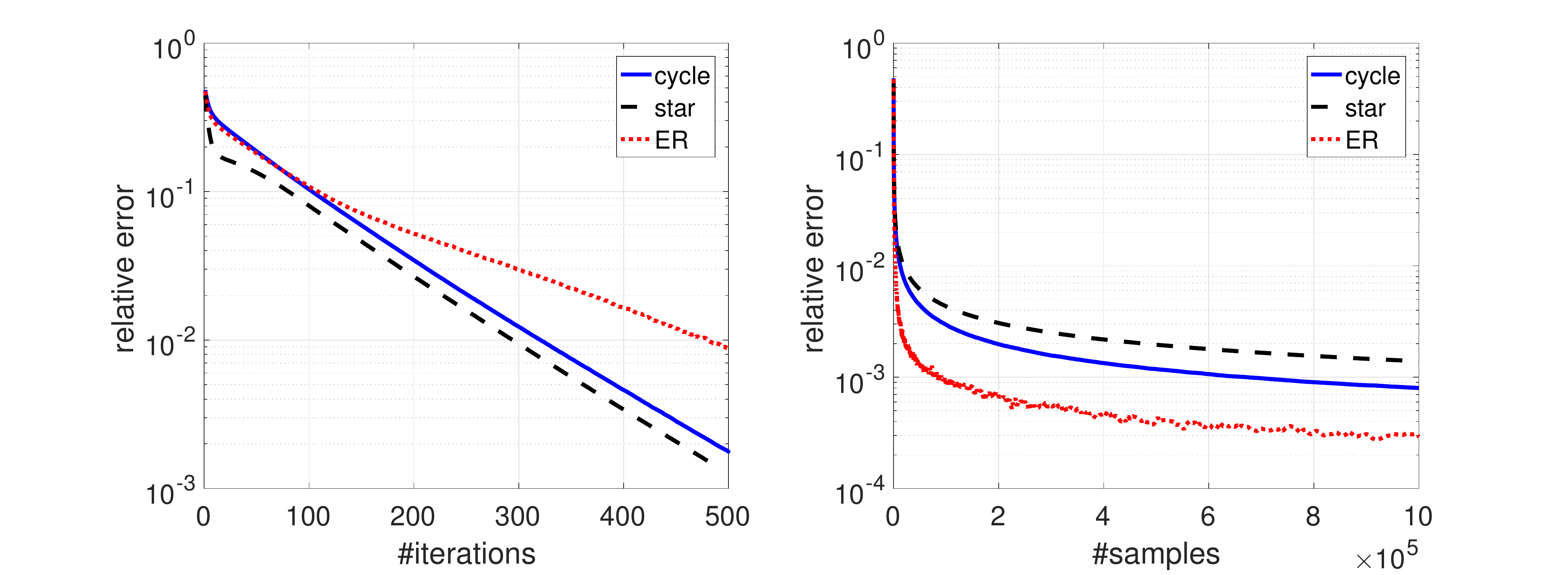}
 \caption{Rate and oracle complexity of  Algorithm \ref{algo-aggregative}}  \label{fig}
\end{minipage}
\end{figure}

We  then  run the VS-PGR algorithm (namely,  Algorithm \ref{algo-aggregative} over a complete graph)  with geometrically  and polynomially increasing sample-sizes and demonstrate the results in Figure \ref{fig_ply}.
The results in Fig. \ref{fig_ply}(a) show  the rate of convergence,  implying that with a small number of proximal evaluations  there is no big difference on convergence rate, while the  algorithm with  geometrically increasing sample-size  will  outperform   the algorithm with polynomially increasing sample-size if  more proximal evaluations are available.
 Fig. \ref{fig_ply}(b)  and  Fig. \ref{fig_ply}(c)  demonstrate the  total number of samples required    to obtain an $\epsilon$-NE, where  it is shown in  Fig. \ref{fig_ply}(b)  that with low accuracy $\epsilon$ the polynomial sample-size with smaller  degree $v$ appears to have better oracle complexity, while  for a high accuracy $\epsilon$,  the geometrically     and     polynomially   increasing (with larger $v$)   sample-size  may have   better oracle complexity.
The numerical results are consistent with the discussions in Remark \ref{rem-poly}.

\begin{figure}[!htb]
  \centering
    \includegraphics[width=6.6in]{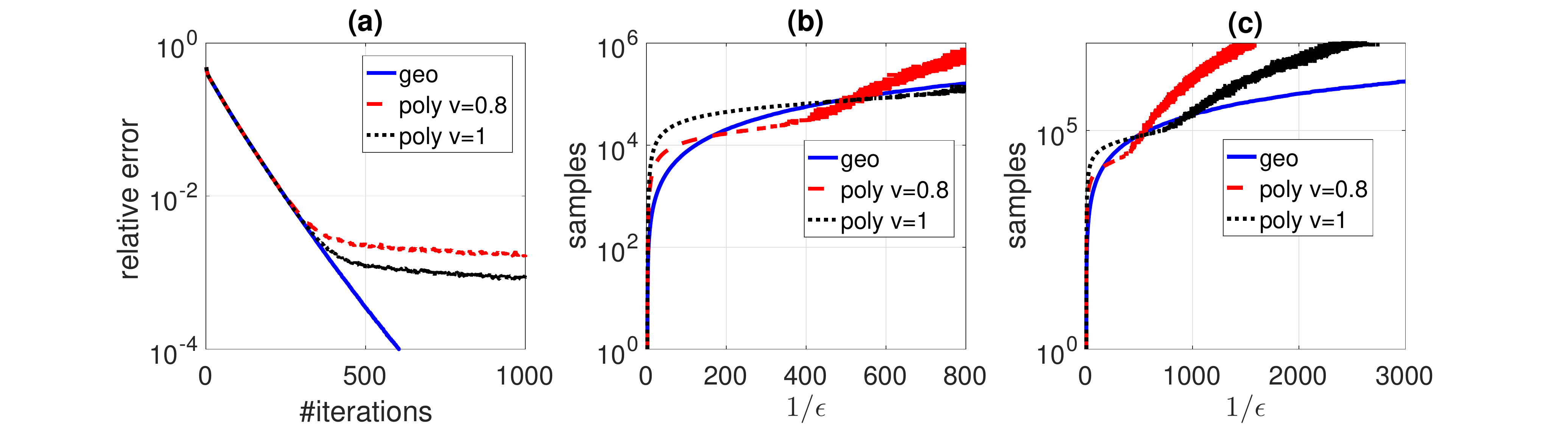}
 \caption{VS-PGR with geometrically  and polynomially increasing sample-size}  \label{fig_ply}
\end{figure}

{\em  Comparison  with  stochastic gradient descent (SGD):}
We set   $n=20$ and  compare    Algorithm \ref{algo-aggregative} and SGD by running  both schemes over  the    Erd\H{o}s--R\'{e}nyi   graph  up to  $10^6$ samples.    We show the results in  Table~\ref{tab3} and Figure \ref{fig2},  where   SGD-$t$ denotes the minibatch SGD algorithm that  utilizes $t$ samples  at each iteration while in
 Algorithm \ref{algo-aggregative} we set  $S_k=\left\lceil  \beta ^{-(k+1)/2} \right \rceil$ and $\tau_k=k+1$.  Though it is seen from  Table~\ref{tab3}   that SGD can obtain slightly better  empirical error,   Algorithm \ref{algo-aggregative}  can significantly reduce the computation time and   the rounds of communication.
We can  also observe from the  iteration complexity    demonstrated   in Figure~\ref{fig2}   that
Algorithm  \ref{algo-aggregative} requires fewer  proximal evaluations than SGD  for   approximating an NE with the same accuracy.

  \begin{table}[htb]
 \newcommand{\tabincell}[2]{\begin{tabular}{@{}#1@{}}#2\end{tabular}}
\centering
 \begin{tabular}{|c|c|c|c|}
 \hline    &     \tabincell{c} {   Algorithm  \ref{algo-aggregative}}   &   \tabincell{c} { SGD-16 , $\tau_k=\ln(k)$} &     \tabincell{c}{SGD-8 , $\tau_k=\ln(k)$}
 \\ \hline  \hline
 emp.err& 5.74e-04 & 2.58e-04 &2.48e-4\\ \hline
  prox.eval &   469 &  6.25e+4  &1.25e+5 \\ \hline
comm.&   1.11e+5  &6.55e+5  &1.41e+6  \\ \hline
 CPU(s)   & 1.8 & 14.67&28.76\\ \hline
 \end{tabular}
 \captionof{table}{Comparison of Algorithm\ref{algo-aggregative} and   SGD} \label{tab3}
\end{table}

\begin{figure}[htb]
\centering
\begin{minipage}[t]{0.49\textwidth}
\centering
   \includegraphics[width=2.6in]{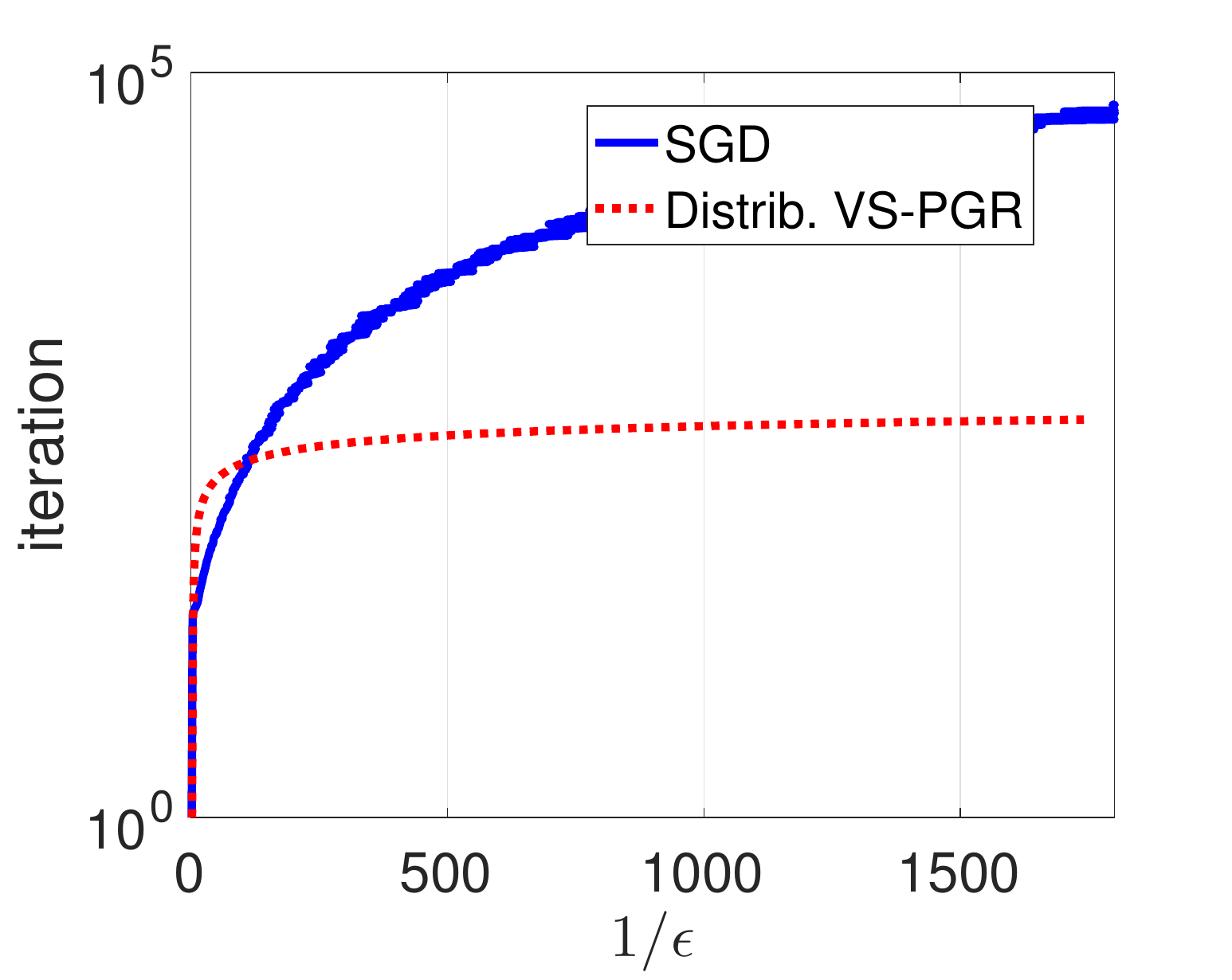} \\
  \caption{Iteration Complexity} \label{fig2}
\end{minipage}
\begin{minipage}[t]{0.02\textwidth}
\end{minipage}
\begin{minipage}[t]{0.49\textwidth}
\centering
 \includegraphics[width=2.6in]{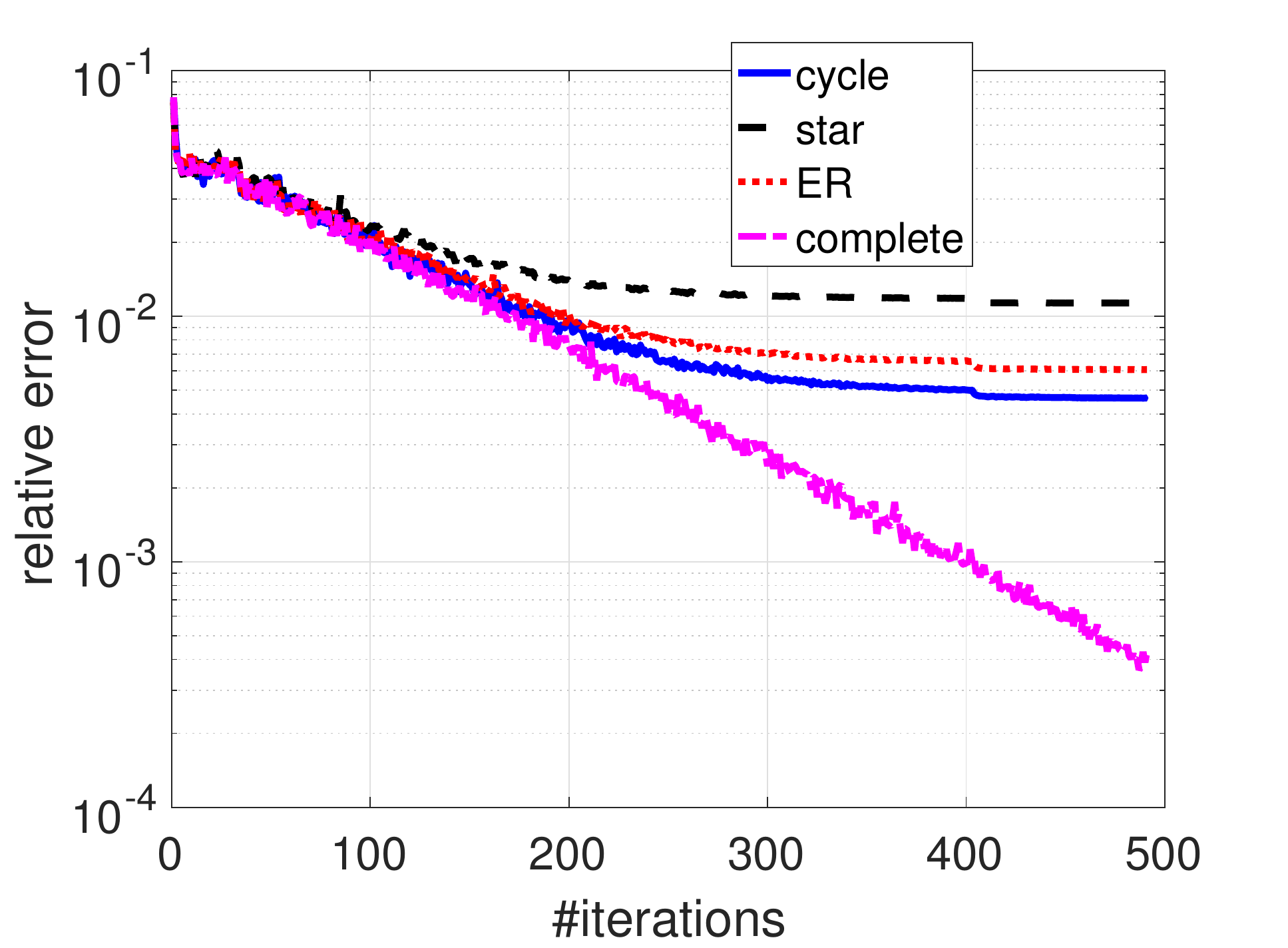}
\caption{   Rate of  Algorithm  \ref{algo-aggregative2} }   \label{pic_rate}
\end{minipage}
\end{figure}

{\em We now validate the performance of the {\bf distributed  VS-PBR} scheme.}
Suppose that for each firm $i$,   there exists a random quadratic  production cost
function $c_i(x_i;\xi_i)=(c_i+\xi_i)  \sum_{l=1}^Lx_i^l +{\rho_i \over 2} x_i^Tx_i$  for some   $c_i>0 $ and  random disturbance  $\xi_i$ with mean zero. We choose some parameters such that Assume that $\lambda_{\min} ( \rho_i \mathbf{I}_L+2diag(b)) >(n-1) \lambda_{\max}(diag(b))
 $ for each $i \in \mathcal{N}.$ Then by the definition of $\Gamma$ in  \eqref{matrix-hessian},   $\| \Gamma\|_{\infty}<1$  and  thus  the proximal BR map is contractive.
Set $n=13 , L=6$, and $\textrm{cap}_{il}=2 $ for any $ i\in \mathcal{N} $ and $ l\in \mathcal{L}$.  We  then run  Algorithm \ref{algo-aggregative2}  with  $\tau_k=\lceil \log(k) \rceil  $, $\mu=20,$ and  $S_k=\left\lceil 0.98^{-(k+1) } \right \rceil$, and demonstrate the    convergence rate   in Fig.  \ref{pic_rate},   showing  that a  better network connectivity may lead to a  faster  rate.

{\em Comparison of distributed VS-PBR and VS-PGR.}
Let  the network  be randomly generated by the  Erd\H{o}s--R\'{e}nyi   graph.
We run Algorithm   \ref{algo-aggregative} with $\alpha=0.04$  and  Algorithm   \ref{algo-aggregative2} with $\mu=30,$
where  $\tau_k=k+1$   and  $S_k=\left\lceil 0.98^{-(k+1) } \right \rceil$. The numerical results for both  schemes are shown in Figure \ref{fig_BPR},  from which it is seen that distributed VS-PBR  has  faster convergence   rate  since  it  requires solving a deterministic optimization problem   per iteration while   distributed VS-PGR merely takes a  proximal gradient step.  Furthermore, the  demonstrated  oracle and communication  complexity  show   that distributed VS-PBR  necessitates less samples and communication  rounds than distributed   VS-PGR to obtain an $\epsilon$-NE.

 \begin{figure}[!htb]
       \centering
 \includegraphics[width=7in]{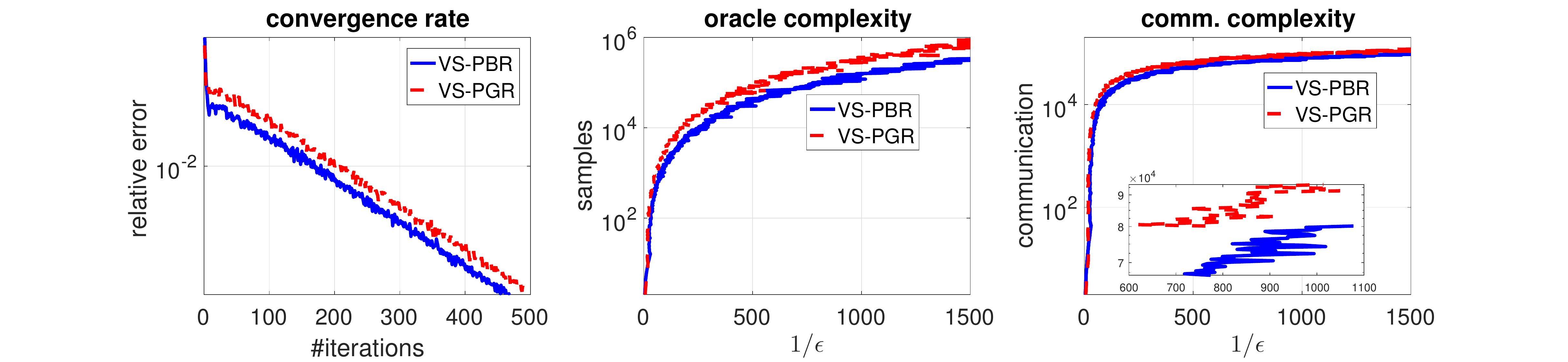}
   \caption{Comparison of distributed VS-PBR and VS-PGR. }   \label{fig_BPR}
 \end{figure}

\section{Concluding Remarks}\label{sec:conclusion}
 Stochastic NEPs and their networked variants represent an important
generalization of deterministic NEPs. While there have been recent
advances in the computation of equilibria in both distributed and stochastic
regimes, at least three gaps currently exist in the available rate statements:
(i) Gap between between rate statements for deterministic and stochastic Nash
games; (ii) Little available by way of implementable best-response schemes in
stochastic regimes; (iii) Lack of computational, oracle, and communication
complexity statements for distributed gradient and best-response schemes.
Motivated by these gaps, we consider four distinct schemes for the resolution
of  a class of stochastic convex NEPs  where each player-specific
objective is a sum of an expectation-valued smooth function and  a convex
nonsmooth function: (i)  VS-PGR scheme for strongly
monotone stochastic NEPs; (ii)  VS-PBR  for stochastic NEPs with
contractive proximal  BR maps,   (iii) Distributed VS-PGR for strongly
monotone aggregative NEPs, and  (iv) Distributed VS-PBR    for stochastic
aggregative games with contractive proximal BR maps.  Under suitable
geometrically increasing sample-size and  linearly increasing  consensus steps
for the distributed variants, we show that all schemes generate sequences that
converge at the ({\bf optimal}) geometric  rate and derive bounds on the
computational, oracle, and communication complexity.  We further quantify the rate and complexity bounds
of the schemes when the sample-sizes and the rounds of
communications increase at  prescribed polynomial rates.

 \appendix

  \section{Proof of results stated in Section \ref{sec:monotone}}

  \subsection{Proof of  Lemma \ref{lem2}} \label{ap:lem2}
  Consider the reformulation  ({VS-PGR}).
 By  using  the nonexpansive property of the proximal operator and   Eqn. \eqref{FP},
  $\|x_{k+1}-x^*\|^2$  can be bounded as follows:
\begin{align}\label{bd-xs}
\|x_{k+1}-x^*\|^2& =\big \| \textrm{prox}_{\alpha  r}  \left[ x_k-\alpha \left( G(x_k )+\bar{w}_{k,S_k} \right) \right]
 - \textrm{prox}_{\alpha  r}  \left[ x^*-\alpha   G(x^* ) \right] \big\|^2\notag
\\& \leq \big \|   x_k -x^*-\alpha \left( G(x_k )-G(x^*)\right)-\alpha\bar{w}_{k,S_k} \big\|^2\notag
\\& \leq \|   x_k -x^* \|^2- 2\alpha \left( G(x_k )-G(x^*)  \right)^T (  x_k -x^*   ) +\alpha ^2\left \|  G(x_k )-G(x^*) \right \|^2  \notag \\&\quad ~- 2\alpha (x_k -x^*)^T \bar{w}_{k,S_k} +2\alpha^2 \left( G(x_k )-G(x^*)\right)^T \bar{w}_{k,S_k}
+ \alpha ^2\|\bar{w}_{k,S_k} \|^2
\\&  \leq  \|   x_k -x^* \|^2- 2\alpha\eta \|x_k-x^*\|^2  +\alpha^2 \| G(x_k )-G(x^*)\|^2 \blue{+\alpha^2 \| x_k -x^*\|^2} \notag
\\&\quad \blue{~+\|\bar{w}_{k,S_k}\|^2+\alpha^2 \| G(x_k )-G(x^*)\|^2+\alpha^2 \| \bar{w}_{k,S_k}\|^2 + \alpha ^2\|\bar{w}_{k,S_k} \|^2}\notag
\\& \blue{\leq  \left(1- 2\alpha\eta+\alpha^2+2\alpha ^2L^2\right) \|x_k-x^*\|^2   + (1+2\alpha ^2)\|\bar{w}_{k,S_k}) \|^2}, \notag
\end{align}
\blue{where the  third inequality  holds by using Assumption \ref{assp-noise}(ii)
and $2ab\leq a^2+b^2$, while  the last inequality follows from Assumption  \ref{assp-noise}(i).}
{Note by Assumption \ref{assp-noise}(iii) and $S_k\geq 1/\alpha^2$  that for any $k\geq 0:$
\begin{align*}
 \mathbb{E}[ \|\bar{w}_{k,S_k} \|^2\mid \mathcal{F}_k] & \leq  \tfrac{2\nu_1^2\|x_k-x^*\|^2 +2\nu_1^2\|x^*\|^2+ \nu_2^2}{S_k}
\leq 2\alpha^2\nu_1^2\|x_k-x^*\|^2   + \tfrac{ 2\nu_1^2\|x^*\|^2+ \nu_2^2}{S_k}, ~ a.s.
\end{align*}}
Then by taking expectations  conditioned on $\mathcal{F}_k$  on both sides of  \blue{\eqref{bd-xs}}, recalling that $x_k$ is adapted to  $\mathcal{F}_k$,  \blue{and using the above inequality, we obtain that}
 \begin{align*}
\mathbb{E}[ \|x_{k+1}-x^*\|^2 |\mathcal{F}_k] &  \leq  \left(1- 2\alpha\eta +\blue{(2L^2+1)}\alpha^2+
2\alpha^2(1+2\alpha^2)\nu_1^2  \right) \|   x_k -x^* \|^2
\\\quad +   \tfrac{2(1+2\alpha^2)\nu_1^2\|x^*\|^2 + (1+2\alpha^2)\nu_2^2}{S_k}
&   \leq  \left(1- 2\alpha\eta +\alpha ^2{\tilde L}^2\right) \|   x_k -x^* \|^2  + \tfrac{\nu^2}{S_k},
\end{align*}
\blue{where the last inequality holds by  the definitions of  $\nu$ and $ \tilde{L}$.} \hfill $\Box$

  \subsection{Proof of  Lemma \ref{lem-recur}} \label{ap:lem-recur}
  Based on \eqref{lem-recur1}, we obtain that
\begin{align} v_{k+1} \leq q^{k+1} v_0 + c_1 \sum_{m=1}^{k+1} q^{k+1-m}\rho^{ m}. \label{bd-linear}
\end{align}
(i)  Consider $q\neq \rho$.  For $\rho<q,$  we obtain the following bound:
$$ \sum_{m=1}^{k } q^{k+1 -m}\rho^{ m}= q^{k+1}\sum_{m=1}^{k } (\rho/q)^{ m} \leq \tfrac{ \rho/ q}{1 -\rho/q}q^k=\frac{1}{ q/\rho-1} q^{k+1}.$$
Similarly,  $ \sum_{m=1}^{k } q^{k+1 -m}\rho^{ m} \leq  \tfrac{1}{\rho/q-1} {\rho}^{k+1}$ when $ q<\rho.$
Therefore,  \[ \sum_{m=1}^{k } q^{k+1 -m}\rho^{ m} =\tfrac{ \max \{\rho,q\}^k}{\max\{ q/\rho, \rho /q \}-1}   {\rm~when~}  \rho \neq q.\]
  Hence (i)  follows  by \eqref{bd-linear} and  $v_0\leq c_0$.

\noindent (ii)      Recall from \cite[Lemma 2]{lei2017synchronous} that
$kq^k\leq \tilde{q}^k/ \ln((\tilde{q}/q)^e) $ for any $\tilde{q}\in(q,1).$
This together with \eqref{bd-linear}  and $ \rho=q $   implies that
$v_k      \leq q^{k }c_0 +c_1  kq^{ k}\leq  (c_0+ \tfrac{c_1}{\ln((\tilde{q}/q)^e)})\tilde{q}^k. $ \hfill $\Box$

\subsection{Proof of Proposition \ref{poly-lem}}\label{ap:poly-lem}
By taking unconditional expectations on both sides of  Eqn. \eqref{prp-bd},
 and using $\blue{S_k\geq \alpha^{-2}(k+1)^v} $, we obtain that for any $k\geq 0:$
\[\mathbb{E}[ \|x_{k+1}-x^*\|^2  ]     \leq  q\mathbb{E}[ \|x_{k }-x^*\|^2 ] +  \alpha ^2\nu^2 (k+1)^{-v} . \]
 Hence
\begin{align}
\mathbb{E}[ \|x_{k+1}-x^*\|^2  ]    \leq q^{k+1}\mathbb{E}[ \|x_{0 }-x^*\|^2  ] +\alpha ^2\nu^2 \sum_{m=1}^{k+1} q^{k+1-m} m^{-v} . \label{bd-poly}
\end{align}
Since $q\in (0,1) $ and $ v>0$, we have that  $q^{-m}m^{-v} \leq   \int_{m }^{m+1} q^{-t} (t-1)^{-v} dt$  for any $ m\geq 2 $. Thus,
\begin{equation}\label{poly-sum}
\begin{split}  \sum_{m=1}^{k+1} q^{k+1-m} m^{-v}  & =q^{k+1} \sum_{m=1}^{{\lceil  2v/ \ln(1/q) \rceil}} q^{ -m} m^{-v}+q^{k+1}\sum_{m=\lceil 2v/ \ln(1/q)  \rceil +1}^{k+1} q^{ -m} m^{-v} \\&  \leq   q^{k+1} \sum_{m=1}^{{\lceil  2v/ \ln(1/q) \rceil}}  q^{-m}+q^{k+1}\int_{\lceil 2v/ \ln(1/q)   \rceil+1}^{k+2}  \tfrac{(q^{-1})^t}{(t-1)^v } dt
\\&  \leq  q^{k}\tfrac{(1/q)^{ {\lceil  2v/ \ln(1/q) \rceil}}-1}{q^{-1}-1}+q^k\int_{\lceil  2v/ \ln(1/q)   \rceil }^{k+1 }  \tfrac{(q^{-1})^t}{t^v } dt.
\end{split}
\end{equation}   Integrating by parts, we obtain that
\begin{align}\label{integ-par}
    \int_a^b  \tfrac{(q^{-1})^t}{t^v } dt =   \int_a^b  \tfrac{1}{t^v}   \left( \tfrac{(q^{-1})^t}{ \ln(q^{-1})} \right)' dt  =
 \tfrac{(q^{-1})^t}{t^v \ln(q^{-1})}\Big|_{a}^b + \int_a^b   \tfrac{ v}{t \ln(q^{-1})} \tfrac{(q^{-1})^t}{t^{v }} . \end{align}
Note that $\tfrac{v}{t \ln(1/q)} \leq \tfrac{1}{2}$ when $t\geq \lceil 2v/ \ln(1/q)  \rceil $. Therefore,
 by setting $a={\lceil 2v/ \ln(1/q)  \rceil }$, $ b={k+1} $ in \eqref{integ-par},   the following holds:
\begin{equation}\label{poly-int}
\begin{split}   & \int_{\lceil 2v/ \ln(1/q)  \rceil }^{k+1}  \tfrac{(q^{-1})^t}{t^v } dt \leq
 \tfrac{(q^{-1})^t}{ t^v \ln(q^{-1})}\Big| _{\lceil 2v/ \ln(1/q)  \rceil }^{k+1}  +{1 \over2}  \int_{\lceil 2v/ \ln(1/q)  \rceil }^{k+1}     { (q^{-1})^t \over t^v}  \\
& \Rightarrow  \int_{\lceil 2v/ \ln(1/q)  \rceil }^{k+1}  \tfrac{(q^{-1})^t}{t^v } dt   \leq
\tfrac{2(q^{-1})^t}{t^v \ln(q^{-1})}\Big| _{\lceil 2v/ \ln(1/q)  \rceil }^{k+1}
\leq \tfrac{2q^{-(k+1)} (k+1)^{-v}}{\ln(1/q) }  . \end{split}
\end{equation}
Note that  $(1/q)^{ 2v/ \ln(1/q)}=\left(e^{\ln(1/q)}\right)^{ 2v/ \ln(1/q)}=e^{2v}$ and hence
$(1/q)^{\lceil  2v/ \ln(1/q) \rceil} \leq e^{2v}/q$.
Then by substituting  \eqref{poly-int} into \eqref{poly-sum}, we have that
\begin{align}\label{lem3-poly-bd}  \sum_{m=1}^{k+1} q^{k+1-m} m^{-v}    \leq  q^{k+1} \tfrac{e^{2v}q^{-1}-1}{1-q}+\tfrac{2q^{-1} (k+1)^{-v}}{\ln(1/q) }.
\end{align}
This incorporated with  \eqref{bd-poly}   and    $\mathbb{E}[ \|x_0-x^*\|^2] \leq C$   {produces \eqref{poly-rate}}.

 Since $q\in (0,1)$ and $v>0$,  by Lemma~\ref{bd-cqv} {with $u=1$}, $q^k \leq c_{q,v}  k^{-v}   $ with  $c_{q,v} \triangleq   {e^{-v}\left(\tfrac{v}{\ln(1/q)}\right)^v} $. Then  {by \eqref{poly-rate}, we conclude that for any $ k\geq 1$,}
$$\mathbb{E}[ \|x_{k }-x^*\|^2  ]      \leq  \left(   Cc_{q,v} + \alpha^2\nu^2c_{q,v}  {e^{2v} q^{-1}-1 \over  1-q} + { 2\alpha ^2\nu^2 q^{-1}\over  \ln(1/q) }  \right) k^{-v}  \triangleq  C_v  k^{-v}  .$$
  {Then for any $  k\geq  K(\epsilon)\triangleq  \left( { C_v  \over \epsilon }\right)^{1/v}$, $  \mathbb{E}[ \|x_{k}-x^*\|^2  ] \leq \epsilon$.}  {By noting  that $C_v=\mathcal{O}(e^v v^v)$,   the iteration complexity is  $\mathcal{O} (  {v} (1/\epsilon)^{1/v})$.}
Therefore, the  {number of sampled gradients required to obtain}   an  $\epsilon-$NE  is   bounded by
\begin{align*}  \sum_{k=0}^{ K(\epsilon) -1 }    (k+1)^v  & = (K(\epsilon))^v+ \sum_{k=1}^{ K(\epsilon) -1 }    k^v  \leq  (K(\epsilon))^v+  \int_{1}^{ K(\epsilon)  } t^v dt \\&= \tfrac{ C_v}{\epsilon }+\tfrac{t^{v+1} }{ v+1}\Big |_{1}^{ K(\epsilon) }
 =\tfrac{ C_v}{\epsilon } + (v+1)^{-1} \left(\tfrac{ C_v}{\epsilon }\right)^{1+{1\over v}} . \end{align*}
Therefore, the oracle complexity is $\mathcal{O} \left(  {e^v v^v}\left({ 1/\epsilon }\right)^{1+{1\over v}}  \right).$
\hfill $\Box$

  \subsection{Proof of  Lemmas  \ref{lem6} and  \ref{lem7}} \label{ap:lem7}
\
\newline

{\bf Proof of Lemma \ref{lem6}.}   By  the consensus step in
Algorithm~\ref{algo-aggregative}, we note  that $\hat{v}_{i,k}=\sum_{j=1}^n [A(k)]_{ij} v_{j,k}$.  Then by \cite[(16) in Lemma 4]{koshal2016distributed}, we obtain the following 
\begin{align}\label{agg-bd1}
\|y_k-\hat{v}_{i,k}\|&\leq \sum_{j=1}^n \left| {1\over n}-[\Phi(k,0)]_{ij}\right| \|v_{j,0}\|  +\sum_{s=1}^k \sum_{j=1}^n
\left| {1\over n}-[\Phi(k,s)]_{ij}\right| \|x_{j,s}-x_{j,s-1}\|.
\end{align}
By definition of $\Phi(k,s)$ and $A(k)$, we have that $\Phi(k,s)=A^{\sum_{p=s}^k \tau_p}.$ Then by  {Lemma \ref{lem-pre}(i) it follows that for any $i,j\in \mathcal{N}$:} $
\left| {1\over n}-[\Phi(k,s)]_{ij}\right| \leq \theta \beta^{\sum_{p=s}^k \tau_p},~ \forall k\geq s  .$
Then by substituting  {this} bound into   \eqref{agg-bd1} we obtain that
\begin{align*}
\|y_k-\hat{v}_{i,k}\|&\leq \theta \beta^{\sum_{p=0}^k \tau_p}\sum_{j=1}^n \|v_{j,0}\|  +   \theta \sum_{s=1}^k \beta^{\sum_{p=s}^k \tau_p}   \sum_{j=1}^n ( \|x_{j,s}\| +\|x_{j,s-1}\|),
\end{align*}
and hence by defining   $D_{\mathcal{R}} \triangleq  \sum_{j=1}^n \max\limits_{x_j\in \mathcal{R}_j} \| x_{j}\|$, we obtain \eqref{agg-bd0}.
\hfill $\Box$

{\bf Proof of Lemma \ref{lem7}.}    Since   $\tau_k= \lceil (k+1)^{u} \rceil$ and $u>0$, we obtain that
\begin{align} \label{sum-tau-bd}
\sum_{p=s}^k \tau_p \geq   \sum_{p=s+1}^{k+1}  p^{u}   \geq \int_{s}^{k+1} t^u dt
 ={t^{u+1} \over u+1} \Big |_{s }^{k+1}
={(k+1)^{u+1}-s^{u+1} \over u+1} .\end{align}
Then by $\beta\in (0,1)$, the following holds with $b   \triangleq \beta^{-1/(u+1)}:$
\begin{align} \label{cor5-poly-bd1}\sum_{s=1}^k  \beta^{\sum_{p=s}^k \tau_p} \leq \beta^{\tfrac{(k+1)^{u+1}}{u+1}} \sum_{s=1}^k b^{s^{u+1}  }=\beta^{\tfrac{(k+1)^{u+1}}{u+1}} \left( b^{k^{u+1}  }+\sum_{s=1}^{k-1} b^{s^{u+1}  }\right).
\end{align}
By defining $t$ as $t=s^{u+1} $, implying that $s=t^{1/(u+1)} $  and $ds=\tfrac{1}{u+1} t^{-u/(u+1)} dt$. Then from  $b>1 $ it follows that
\begin{align} \label{cor5-poly-bd2}
\sum_{s=p}^{k-1} b^{s^{u+1}  }  \leq\int_{ p }^{k} b^{s^{u+1}} ds=\tfrac{1}{u+1} \int_{p^{u+1}}^{k^{u+1}} \tfrac{b^t}{t^{u/(u+1)}} dt.
\end{align}
 {Using \eqref{integ-par}  with $v=u/(u+1)$ and $q=b^{-1}$, we obtain that
$$   \int  \tfrac{b^t}{t^{u/(u+1)}} dt =\tfrac{
b^t}{ t^{u/(u+1)}\ln(b)} + \tfrac{ u}{(u+1) } \int \tfrac{1}{t\ln(b)} \tfrac{b^t dt}{t^{ u/(u+1)}}. $$}
Define   $k_0 \triangleq  \left \lfloor \ln(b)^{-1/(u+1)}   \right\rfloor$.  Then $(k_0+1)^{u+1} \geq  \tfrac{{1}}{\ln(b)} $ and
 $ { 1\over t  \ln(b)}\leq  1~{\rm if }~ t \geq  (k_0+1)^{u+1} $. Thus,
\begin{align*}  &  \int_{(k_0+1)^{u+1} }^{k^{u+1}}   {b^t \over t^{u/(u+1)}} dt  \leq {
b^t\over t^{u/(u+1)}\ln(b)} \Big |_{(k_0+1)^{u+1} }^{k^{u+1}}+  {u\over u+1}\int_{(k_0+1)^{u+1} }^{k^{u+1}}   {b^t \over t^{u/(u+1)}} dt
\\&  \Rightarrow  {1\over u+1} \int_{(k_0+1)^{u+1} }^{k^{u+1}}   {b^t \over t^{u/(u+1)}} dt    \leq  {
b^t\over t^{u/(u+1)}\ln(b)} \Big |_{(k_0+1)^{u+1} }^{k^{u+1}}
<  {  b^{k^{u+1}}\over k^u\ln(b)}.
\end{align*}
This  together  with \eqref{cor5-poly-bd1}, \eqref{cor5-poly-bd2}, and $b= \beta^{-1/(u+1)}>1$ implies that
\begin{align*}  \sum_{s=1}^k  \beta^{\sum_{p=s}^k \tau_p} &\leq \beta^{{(k+1)^{u+1}  \over u+1}}\sum_{s=1}^{k_0}  b^{s^{u+1}  }+\beta^{{(k+1)^{u+1}  \over u+1}}
\left( \sum_{k_0+1}^{k-1}  b^{s^{u+1}  }  +  b^{k^{u+1}  }\right)
\\&  \leq \beta^{{k^{u+1}  \over u+1}} k_0  b^{k_0^{u+1}  }+\beta^{{k^{u+1}  \over u+1}}b^{k^{u+1}} \left(1+  {  1\over k^u\ln(b)}\right).
\end{align*}
  Then by the fact that $ b^{k_0^{u+1}} \leq  b^{1/\ln(b)}    =  e $ since  $k_0^{u+1} \leq 1/\ln(b) $,
  and using $b=\beta^{-1/(u+1)}$,  we  proved the lemma.      \hfill $\Box$

 \subsection{Proof of Proposition \ref{prp2}}\label{ap:prp2}
Similarly to Lemma \ref{lem1}, $x^* \in X$ is an NE   if and only if  $x^*$   is a fixed point of  $\textrm{prox}_{\alpha r}(x -\alpha \phi(x )) $.   Then by using \eqref{alg-agg-strategy}, \blue{\eqref{equi-phiF}, $2ab\leq a^2+b^2$,}
  and the non-expansive property of the proximal operator, we  have that
\begin{align}\label{bd-relation}
\notag  \ \|x_{i,k+1}-x_i^*\|^2   &\leq \|x_{i,k}-x_{i}^*-\alpha \big(\blue{F_i}(x_{i,k}, n\hat{v}_{i,k}) \notag  -\blue{\phi_i(x^*)}\big) -\alpha e_{i,k}\|^2
  \\  \notag & \leq   \|x_{i,k}-x_{i}^*\|^2+\alpha^2
{\big \| \blue{F_i}(x_{i,k}, n\hat{v}_{i,k}) -\blue{\phi_i(x^*)} \big\|^2}+\alpha^2 \| e_{i,k}\|^2
\\ &  \notag \quad \blue{-2  \alpha ( x_{i,k}-x_{i}^* )^Te_{i,k}+2\alpha^2 \big( F_i (x_{i,k}, n\hat{v}_{i,k}) - \phi_i(x^*)  \big)^Te_{i,k}}
\\ &\quad-2\alpha{ (x_{i,k}-x_{i}^*)^T\big(\blue{F_i}(x_{i,k}, n\hat{v}_{i,k}) -\blue{\phi_i(x^*)}   \big)}
 \\ \notag & \leq    (1+\alpha^2) \|x_{i,k}-x_{i}^*\|^2+2\alpha^2
\underbrace{\big \| \blue{F_i}(x_{i,k}, n\hat{v}_{i,k}) -\blue{\phi_i(x^*)} \big\|^2}_{\small \rm Term \ 1}
 \\  \notag & \quad-2\alpha \underbrace{ (x_{i,k}-x_{i}^*)^T\big(\blue{F_i}(x_{i,k}, n\hat{v}_{i,k}) -\blue{\phi_i(x^*)}  \big)}_{\small \rm Term \ 2}+ (1+2\alpha^2)  \| e_{i,k}\|^2.
\end{align}

We now estimate   Term 1 and Term 2 of   Eqn.  \eqref{bd-relation}.
 \blue{Note by Assumption  \ref{assump-agg1}(i) and Lemma \ref{lem-pre}(ii) that $\{y_k\}$ is a bounded sequence,
 hence from Lemmas \ref{lem6} and \ref{lem7}  it is seen that $\hat{v}_{i,k}$ is bounded as well. Then by} using the triangle   inequality,  $(a+b)^2\leq 2(a^2+b^2) $,  Assumption \ref{assump-agg2}(iv), \blue{and \eqref{equi-phiF},}
 Term 1 may be  bounded as follows:
\begin{align}\notag \label{bd-term1}
 {\rm Term \ 1 } & \leq 2 \big \| \blue{F_i}(x_{i,k}, n\hat{v}_{i,k}) -\blue{F_i}(x_{i,k}, ny_k)\|^2  + 2\|\blue{F_i}(x_{i,k},ny_k)-\phi_i(x^*)\big\|^2
\\ &  \leq 2n^2L_i^2 \|  \hat{v}_{i,k}-y_k\|^2  +  2\|\blue{\phi_i(x_k)-\phi_i(x^*)}\big\|^2 .
\end{align}
By definition    $ D_{\mathcal{R}} \triangleq \sum_{j=1}^n \max\limits_{x_j\in \mathcal{R}_j} \| x_{j}\|$ and Assumption \ref{assump-agg2}(iv),  the following holds:
\begin{align*}
  & (x_{i,k}-x_i^*)^T (\blue{F_i}(x_{i,k},n\hat{v}_{i,k})-\blue{F_i}(x_{i,k},ny_k) )
  \\ & \geq -\| x_{i,k}-x_i^*\| \|\blue{F_i}(x_{i,k},n\hat{v}_{i,k})-\blue{F_i}(x_{i,k},ny_k)   \|
\\& \geq -nL_i \|x_{i,k}-x_i^*\|   \|  \hat{v}_{i,k}-y_k\|  \geq -2nL_iD_{\mathcal{R}}  \|  \hat{v}_{i,k}-y_k\|.
\end{align*}
\blue{Then by using \eqref{equi-phiF} and   Lemma \ref{lem-pre}(ii)},  Term 2 is  lower bounded by the following:
\begin{align}\label{bd-term2}
{\rm Term \ 2}& =  (x_{i,k}-x_i^*)^T (\blue{F_i}(x_{i,k},\hat{v}_{i,k})  -\blue{F_i}(x_{i,k},ny_k) ) \notag
\\& \quad + (x_{i,k}-x_i^*)^T (\blue{F_i}(x_{i,k},ny_k)-\phi_i(x^*))\notag\\
 & \geq -2nL_iD_{\mathcal{R}}  \|  \hat{v}_{i,k}-y_k\|    + (x_{i,k}-x_i^*)^T (\blue{\phi_i(x_k)}-\phi_i(x^*)) .
\end{align}
\blue{Note by   Assumption \ref{assump-agg2}(v) and $ S_k \geq \alpha^{-2}$  that
\begin{align}\label{bd-eik}
 \mathbb{E}[ \|e_{i,k} \|^2\mid \mathcal{F}_k] &
 \leq \tfrac{2\nu_{i,1}^2\|x_{i,k}-x_i^*\|^2  + 2\nu_{i,1}^2\|x_i^*\|^2+ \nu_{i,2}^2}{S_k} \notag
 \\ & \leq  2 \alpha^2\nu_{i,1}^2\|x_{i,k}-x_i^*\|^2   + \tfrac{ 2\nu_{i,1}^2\|x_i^*\|^2+ \nu_{i,2}^2}{S_k}, ~ a.s.
\end{align}}
Since $\{\tau_k\}$ is a deterministic sequence,   $\hat{v}_{i,k}$ is adapted to
$\mathcal{F}_k $  by Algorithm \ref{algo-aggregative}.   Then by taking expectations  conditioned on $\mathcal{F}_k$ of the   inequality \eqref{bd-relation},  by substituting  \blue{\eqref{bd-term1}, \eqref{bd-term2}, and\eqref{bd-eik}}, we obtain that
\begin{align*}
\mathbb{E}[\|x_{i,k+1}-x_i^*\|^2 | \mathcal{F}_k]
& \leq  \left(1+\alpha^2+\blue{ 2\alpha^2(1+2\alpha^2)\nu_{i,1}^2    }\right)  \|x_{i,k}-x_{i}^*\|^2   +  4\alpha^2 \blue{\|\phi_i(x_k)-\phi_i(x^*)\big\|^2}
 \\&  \quad + 4 \alpha^2n^2L_i^2 \|  \hat{v}_{i,k}-y_k\|^2  -2\alpha(x_{i,k}-x_i^*)^T \blue{ (\phi_i(x_k)-\phi_i(x^*)) }
  \\& \quad +4\alpha n L_iD_{\mathcal{R}}  \| \hat{v}_{i,k}-y_k\|
  \blue{+{ (1+2\alpha^2) (2\nu_{i,1}^2\|x_i^*\|^2+ \nu_{i,2}^2) \over S_k}} .
\end{align*}
Summing the above inequality over $i=1,\cdots,n$, \blue{by recalling the definitions $\bar{\nu}_1 =\max_{1\leq i \leq n} \nu_{i,1}$ and $\bar{\nu}^2_2 = \sum_{i =1}^{n} \nu^2_{i,2}$, we obtain the following,
\blue{where the last inequality follows from  Assumptions \ref{assump-agg2}(ii) and  \ref{assump-agg2}(iii) and   the definitions of  $\varrho_{\phi} $  and $\bar{\nu}$.}}
\begin{align}\label{bd-relation3}
 \mathbb{E}[\|x_{k+1}-x^*\|^2 | \mathcal{F}_k] &   \leq  \left(1+\alpha^2+  2\alpha^2(1+2\alpha^2) {\bar{\nu}_1^2} \right) \|x_{k}-x^*\|^2  +   4\alpha^2  \|\phi (x_k)-\phi (x^*)\big\|^2   \notag
   \\& \notag \quad+4\alpha n D_{\mathcal{R}} \sum_{i=1}^nL_i   \|  \hat{v}_{i,k}-y_k\| + 4\alpha^2n^2 \sum_{i=1}^nL_i^2 \|  \hat{v}_{i,k}-y_k\|^2
    \\ & \quad  -2\alpha(x_k-x^*)^T   (\phi (x_k)-\phi (x^*))   +\tfrac{ (1+2\alpha^2)(2\bar{\nu}^2_1\|x^*\|^2+\bar{\nu}^2_2)}{S_k}
\\& \notag \leq  \blue{ \Big(\underbrace{1-2\alpha \eta_{\phi}+2\alpha^2\big(1/2+    (1+2\alpha^2) {\bar{\nu}_1^2}
 + 2 L_{\phi}^2\big)}_{\triangleq \varrho_{\phi} } \Big) \|x_{k}-x^*\|^2}
 \\& \notag \quad\blue{+4\alpha n D_{\mathcal{R}} \sum_{i=1}^nL_i   \|  \hat{v}_{i,k}-y_k\| + 4\alpha^2n^2 \sum_{i=1}^nL_i^2 \|  \hat{v}_{i,k}-y_k\|^2+\tfrac{\bar{\nu}^2}{S_k}} .
\end{align}
%

 We now estimate the  bound for  $ \|  \hat{v}_{i,k}-y_k\| $ based on  Lemmas~\ref{lem6}--\ref{lem7}.
  By  setting  $u=1$,     from $\beta\in (0,1)$ and   $\sum_{p=0}^k \tau_p=(k+1)(k+2)/2$ it follows that
\begin{equation*}
\begin{split}
&  \|y_k-\hat{v}_{i,k}\|
\\& \leq 2 \theta D_{\mathcal{R}}  \left(  e   \sqrt{1/\ln(\beta^{-1/2}) }   \beta^{{(k+1)^2  \over 2}}+\beta^{{2k+1 \over 2}} \Big({ 2\over k\ln(1/\beta)} +1\Big) \right)
  + \theta D_{\mathcal{R}}  \beta^{(k+1)(k+2)/2}    \\&\leq 2 \theta D_{\mathcal{R}}  \left(  e   \sqrt{1/\ln(\beta^{-1/2}) }+ { 2+\ln(1/\beta)\over \beta^{1/2}\ln(1/\beta)}  \right) \beta^{ k+1} + \theta D_{\mathcal{R}}  \beta^{(k+1)(k+2)/2}
   \\&= C_1\beta^{(k+1)(k+2)/2} + C_2  \beta^{ k+1}, \quad \forall k\geq 1.
\end{split}
\end{equation*}
By \eqref{agg-bd0} it is seen that the above inequality also holds for $k=0$, and hence
\begin{equation}\label{bd-consensus-err}
\begin{split}
  \|y_k-\hat{v}_{i,k}\|  \leq C_1\beta^{(k+1)(k+2)/2} + C_2  \beta^{ k+1},  \quad \forall k \geq 0.
\end{split}
\end{equation}
Since $(a+b)^2\leq 2(a^2+b^2) $,  $S_k \geq   \alpha^{-2}  \rho^{-(k+1)}$,   by taking the unconditional expectations on both sides of  \eqref{bd-relation3},  we obtain that
\begin{align*}
\notag   &\mathbb{E}[\|x_{k+1}-x^*\|^2  ] \leq \varrho_{\phi} \mathbb{E}[\|x_{k }-x^*\|^2 ]
 +4\alpha n D_{\mathcal{R}}  \left(C_1\beta^{(k+1)(k+2)/2} + C_2  \beta^{ k+1 } \right)  \sum_{i=1}^nL_i \\
& + 4\alpha^2n^2 \left(C_1^2\beta^{(k+1)(k+2)} + C_2^2  \beta^{ 2(k+1) }\right)     \sum_{i=1}^nL_i^2+\alpha^2 \rho^{k+1} \bar{\nu}^2,  \quad \forall k\geq 0.
\end{align*}
This  implies    \eqref{agg-prp1} by the definition of $C_3$ in \eqref{prp-defc} and $\beta\in(0,1)$.\hfill $\Box$

\subsection{Proof of Corollary \ref{cor5}}\label{app:cor5}
 We first  estimate the  bound for  $ \|  \hat{v}_{i,k}-y_k\| $.
 {Since $t^{u+1} $ with $ u>0$ is  convex in $x>0,$ we have that
 $(k+1)^{u+1}  -k^{u+1}\geq \nabla x^{u+1} \big|_{x=k}=(u+1) k^u .$
   Hence ${(k+1)^{u+1}    \over u+1} \geq \tfrac{k^{u+1}}{u+1}+k^u $ and ${(k+1)^{u+1}  -k^{u+1}  \over u+1} \geq k^u.$} Then by using \eqref{agg-bd0}, Lemma \ref{lem7},  and  \eqref{sum-tau-bd},   from $\beta\in (0,1)$ it follows that for any $k\geq 1:$
\begin{align*}
   \|y_k-\hat{v}_{i,k}\| & \leq \theta D_{\mathcal{R}}  \big( 1+  2e \big(  \ln(\beta^{-1/(u+1)}) \big)^{ -1 \over u+1  } \big) \beta^{{k^{u+1}    \over u+1}+k^u}
\\& + 2 \theta D_{\mathcal{R}}   \beta^{k^u}  \big(1+  {  u+1\over k^u\ln(1/\beta)}\big)\leq \widetilde{C}_{\beta} \beta^{k^u} , \end{align*}
where  $ \widetilde{C}_{\beta}\triangleq  \theta D_{\mathcal{R}} \beta^{1\over  u+1 }  \big( 1+  2e \left(  \ln(\beta^{-1/(u+1)}) \right)^{ -1\over u+1  } \big)  + 2 \theta D_{\mathcal{R}}    { u+1+\ln(1/\beta) \over \ln(1/\beta)} $.
By \eqref{agg-bd0} it is seen that the above inequality also holds for $k=0$. Hence  $ \|y_k-\hat{v}_{i,k}\|  \leq \widetilde{C}_{\beta} \beta^{k^u}$ for any $ k \geq 0.$  Then using \eqref{bd-relation3},
\blue{$S_k \geq \alpha^{-2} (k+1)^v  $,} and   $\varrho_{\phi}= 1-2\alpha\eta_{\phi} +2\alpha^2 L_{\phi}^2 $, we have
\begin{equation} \begin{split}\label{cor5-bd}
  & \mathbb{E}[\|x_{k+1}-x^*\|^2 ] \leq \varrho_{\phi} \mathbb{E}[ \|x_{k}-x^*\|^2]   +\alpha^2 (k+1)^{-v}  \blue{\bar{\nu}^2}   \\&  +4\alpha n D_{\mathcal{R}}\widetilde{C}_{\beta} \beta^{k^u}  \sum_{i=1}^nL_i   + 2\alpha^2n^2\widetilde{C}_{\beta}^2 \beta^{2k^u}  \sum_{i=1}^nL_i^2
 \\   &\leq \varrho_{\phi} \mathbb{E}[ \|x_{k}-x^*\|^2] +  \widetilde{C}_{\beta} ^o \beta^{k^u} +\alpha^2 (k+1)^{-v}   \blue{\bar{\nu}^2}
 \\ & =\varrho_{\phi}^{k+1}\mathbb{E}[ \|x_{0}-x^*\|^2] +  \alpha^2 \blue{\bar{\nu}^2}  \sum_{m=0}^{k+1} \varrho_{\phi} ^{k+1-m} m^{-v}+  \widetilde{C}_{\beta} ^o  \sum_{m=0}^{k} \varrho_{\phi} ^{k-m} \beta^{m^u}   ,
\end{split}
\end{equation}
where $\widetilde{C}_{\beta} ^o \triangleq 4\alpha n D_{\mathcal{R}}\widetilde{C}_{\beta}   \sum_{i=1}^nL_i       + 2\alpha^2n^2\widetilde{C}_{\beta}^2  \sum_{i=1}^nL_i^2 $. We now estimate the uppeer bound of the last term in the above inequality.
Since $ \varrho_{\phi} \in (0,1)$ and $\beta \in ги0,1)$, $
 \sum_{m=0}^{k} \varrho_{\phi} ^{-m}    \beta^{m^u}  \leq \int_{0}^{k+1} \varrho_{\phi} ^{-t} \beta^{t^u} dt.$
 Integrating by parts,
 \begin{align*}  &  \int  \varrho_{\phi} ^{-t} \beta^{t^u}  dt    =  \int \beta^{t^u}  \left( { \varrho_{\phi}^{-t}\over \ln( \varrho_{\phi}^{-1})} \right)' dt  = {  \beta^{t^u} \varrho_{\phi}^{-t} \over \ln( \varrho_{\phi}^{-1})} - \int { \varrho_{\phi}^{-t}\over \ln( \varrho_{\phi}^{-1})}  (\beta^{t^u} )' dt  \\&= { \beta^{t^u} \varrho_{\phi}^{-t}\over \ln( \varrho_{\phi}^{-1})}  + u\ln(1/\beta) \int { \varrho_{\phi}^{-t}\over \ln( \varrho_{\phi}^{-1})}   \beta^{t^u}  t^{u-1}dt
  = { \beta^{t^u} \varrho_{\phi}^{-t}\over \ln( \varrho_{\phi}^{-1})}  +  \int {u \ln(1/\beta) \over \ln( \varrho_{\phi}^{-1}) t^{1-u}}   \varrho_{\phi}^{-t}\beta^{t^u}  dt  \end{align*}
Note that  $ {\ln(1/\beta) u\over \ln( \varrho_{\phi}^{-1}) t^{1-u}}  \leq {1\over 2}$ when $t\geq t_0  \triangleq\left \lceil
( 2u \ln(1/\beta)/ \ln(1/\varrho_{\phi})   )^{1/(1-u)} \right \rceil  $.
Then \begin{align*}    & \int_{t_0}^{k+1} \varrho_{\phi} ^{-t} \beta^{t^u}  dt   \leq { \beta^{t^u} \varrho_{\phi}^{-t}\over \ln( \varrho_{\phi}^{-1})} \Big|_{t_0}^{k+1} +  \int_{t_0}^{k+1} {1\over 2}   \varrho_{\phi}^{-t}\beta^{t^u}  dt .
\\& \Rightarrow  \int_{t_0}^{k+1}  \varrho_{\phi} ^{-t} \beta^{t^u}  dt  \leq  {2 \beta^{t^u} \varrho_{\phi}^{-t}\over \ln( \varrho_{\phi}^{-1})}  \big|_{t=t_0}^{k+1}
\leq  {2 \beta^{(k+1)^u} \varrho_{\phi}^{-(k+1)}\over \ln( \varrho_{\phi}^{-1})} .  \end{align*}
Note that  $\sum_{m=0}^{t_0-1} \varrho_{\phi} ^{ -m} \beta^{m^u} \leq \sum_{m=0}^{t_0-1} \varrho_{\phi} ^{ -m} \leq {\varrho_{\phi} ^{-t_0} -1 \over \varrho_{\phi} ^{-1} -1 } .$
This incorporated  with \eqref{cor5-bd} and \eqref{lem3-poly-bd} implies that
$ \mathbb{E}[\|x_{k+1}-x^*\|^2 ]  =\mathcal{O} \left((k+1)^{-v}+ e^{2v}  \varrho_{\phi}^{k+1}+  \beta^{(k+1)^u} +  \varrho_{\phi}^{- t_0} \varrho_{\phi}^{k+1}  \right)  .$
Thus, by using  Lemma \ref{bd-cqv}, we obtain that
\begin{align*}\mathbb{E}[\|x_{k+1}-x^*\|^2 ]
&=\mathcal{O} \left( \Big( e^v v^v    + e^{- v/u }{\left(\tfrac{v}{u\ln(1/q)}\right)^{ v/u }} + v^v  e^{-v}
    \varrho_{\phi}^{-t_0} \Big)(k+1)^{-v} \right)
\\&=  \mathcal{O} \left(
  e^v {\left(\tfrac{v}{u }\right)^{ v/u }}      (k+1)^{-v} \right)=\mathcal{O} \left( C_{v,u}    (k+1)^{-v} \right).
\end{align*}

 Then for any $  k\geq  K(\epsilon)\triangleq  \left( { C_{v,u}  \over \epsilon }\right)^{1/v}$, $  \mathbb{E}[ \|x_{k}-x^*\|^2  ] \leq \epsilon$.}  {By noting  that $ C_{v,u} =\mathcal{O}(e^v {\left(\tfrac{v}{u }\right)^{ v/u }} )$,   the iteration complexity is  $\mathcal{O} (  {\left(\tfrac{v}{u }\right)^{ 1/u }} (1/\epsilon)^{1/v})$.
Therefore, the   number of sampled gradients required to obtain    an  $\epsilon-$NE  is   bounded by
\begin{align*}  \sum_{k=0}^{ K(\epsilon) -1 }    (k+1)^v   \leq    \int_{1}^{ K(\epsilon)  } t^v dt  =  {t^{v+1} \over v+1}\Big |_{1}^{ K(\epsilon) } \leq  (v+1)^{-1} \left({  C_{v,u}   \over \epsilon }\right)^{1+{1\over v}} . \end{align*}
Hence,  the oracle complexity is $\mathcal{O} \left(  {e^v {\left(\tfrac{v}{u }\right)^{ v/u }} }\left({ 1/\epsilon }\right)^{1+{1\over v}}  \right).$
The number of  communication rounds required to obtain an $\epsilon-$NE  is   bounded by
  $    \sum\limits_{k=0}^{ K(\epsilon) -1 }  \tau_k   \leq    \int_{1}^{ K(\epsilon)  } t^u dt
  \leq  {K(\epsilon)^{u+1} \over u+1} =
\mathcal{O} \Big(  e^{u} {\left(\tfrac{v}{u }\right)^{ (u+1)/u }}  \big({1\over \epsilon}\big)^{(u+1)/v}\Big).  \qquad   \Box    $

\section{Proof of results stated in Section \ref{sec:BR}}

\subsection{Proof of Lemma \ref{lem4}}\label{ap:lem4}
 Define \[ \bar{w}_{i,k}(x_i)\triangleq {1\over S_k}\sum_{p=1}^{S_k}  \nabla_{x_i} \psi_i(x_i,y_{-i,k}; \xi_{i,k}^p)-\nabla_{x_i}f_i(x_i,y_{-i,k}).\]
By applying  the optimality condition on \eqref{SAA} and  \eqref{prox-BRi} ,  $x_{i,k+1}$  and $\wh{x}_{i}(y_k) $
are respectively a fixed point of  ${\rm prox}_{\alpha r_i} \big[ x_i-\alpha\big( \nabla_{x_i}f_i(x_i,y_{-i,k})+\mu(x_i-y_{i,k})+ \bar{w}_{i,k}(x_i)\big)\big]$ and  ${\rm prox}_{\alpha r_i} \left[ x_i-\alpha  (\nabla_{x_i} f_i(x_i,y_{-i,k}) +\mu(x_i-y_{i,k}) )\right] $ for any $\alpha>0$.  Then by the nonexpansive property of the proximal operator, we have the following:
\begin{align}\label{br-bd1}
 &\|x_{i,k+1}-\wh{x}_{i}(y_k) \|  \leq \alpha\|\bar{w}_{i,k}( x_{i,k+1})\|+
 \\&\Big\|  (1-\alpha \mu) (x_{i,k+1}- \wh{x}_{i}(y_k)) -\alpha\left(\nabla_{x_i} f_i( x_{i,k+1},y_{-i,k})- \nabla_{x_i} f_i(\wh{x}_{i}(y_k) ,y_{-i,k}) \right)\Big\|   \notag
\end{align}
Note by  Assumption  \ref{assp-compact}(ii)  that 
\begin{align*}
& \qquad  \Big\|  (1-\alpha \mu) (x_{i,k+1}- \wh{x}_{i}(y_k)) -\alpha\left(\nabla_{x_i} f_i( x_{i,k+1},y_{-i,k})-  \nabla_{x_i} f_i(\wh{x}_{i}(y_k) ,y_{-i,k}) \right)\Big\| ^2
 \\& = (1-\alpha \mu)^2 \Big\|  x_{i,k+1}- \wh{x}_{i}(y_k) \Big\|^2 + \alpha^2 \left\|\nabla_{x_i} f_i( x_{i,k+1},y_{-i,k})-  \nabla_{x_i} f_i(\wh{x}_{i}(y_k) ,y_{-i,k})  \right \|^2
 \\&\quad -2 \alpha (1-\alpha \mu) \left(\nabla_{x_i}f_i( x_{i,k+1},y_{-i,k})- \nabla_{x_i} f_i(\wh{x}_{i}(y_k) ,y_{-i,k}) \right)^T(x_{i,k+1}- \wh{x}_{i}(y_k) )
\\& \leq \left( (1- \alpha \mu)^2+\alpha^2L_{fi}^2 \right)\left\|  x_{i,k+1}- \wh{x}_{i}(y_k) \right\|^2
 {\rm~when~} 0<\alpha \mu<1,
\end{align*}
which incorporated with \eqref{br-bd1} implies that for any $0<\alpha<1/\mu :$
\begin{align*}
&\|x_{i,k+1}-\wh{x}_{i}(y_k) \| \leq  \sqrt{ (1- \alpha \mu)^2+\alpha^2L_{fi}^2}  \|  x_{i,k+1}- \wh{x}_{i}(y_k)\| + \alpha  \|\bar{w}_{i,k}( x_{i,k+1})\|   .
\end{align*}
In the above inequality,  by setting  $\alpha={\mu \over\mu^2+L_{fi}^2},$ we obtain that
 $ \|x_{i,k+1}-\wh{x}_{i}(y_k) \|\leq  C_{i,b}\|\bar{w}_{i,k}( x_{i,k+1})\|  .$
By Assumption \ref{assp-compact}(iii), there holds $ \mathbb{E}[\|\bar{w}_{i,k}(x_{i,k+1})\|^2] \leq {\nu_i^2\over S_k}$ and
 $\mathbb{E}[\|x_{i,k+1}-\wh{x}_{i}(y_k) \|^2]\leq   C_{i,b}^2 \mathbb{E}[\|\bar{w}_{i,k}(x_{i,k})\|^2]\leq {\nu_i^2  C_{i,b}^2  \over S_k}.  $ \hfill $\Box$

  \subsection{Proof of Proposition  \ref{prp-linear}}\label{ap:prp-linear}
By    $x_i^* = \wh{x}_i (x^*)$,  using    the triangle inequality and   $y_k=x_k$, we obtain that
  $\|x_{i,k+1}- x_i^*\|
			 \leq  \|x_{i,k+1} -\wh{x}_i (x_k)\|
			+\| \wh{x}_{i}(x_k) -
			\wh{x}_i (x^*)\|.$
 Then by   the triangle inequality,   \eqref{cont-prox-best-resp}, and  $a  \triangleq \| \Gamma\|<1$, we have  the following bound:
\begin{align*}
& v_{k+1} \triangleq \left\|   \pmat{ \|x_{1,k+1} -{x}_{1}^*\| \\
				\vdots \\
 \|x_{n,k+1}  - {x}_{n}^*\|} \right \|   \leq
a \left\|  \pmat{
	 \|x_{1,k}-x_1^*\|\\
				\vdots \\
	 \|x_{n,k}- x_n^*\|} \right \| + \left\|  \pmat{
	 \|x_{1,k+1} -  \wh{x}_{1}(x_k)\|\\
				\vdots \\
	 \|x_{n,k+1} - 	 \wh{x}_{n}(x_k)\|} \right\|.
\end{align*}
Therefore, the following holds:
\begin{align*}
& v_{k+1}^2  \leq a^2 v_k^2+ 2av_k\left\|  \pmat{
	 \|x_{1,k+1} -
			 \wh{x}_{1}(x_k)\|\\
				\vdots \\
	 \|x_{n,k+1} -
			 \wh{x}_{n}(x_k)\|} \right\|+\left\|  \pmat{
	 \|x_{1,k+1} -
			 \wh{x}_{1}(x_k)\|\\
				\vdots \\
	 \|x_{n,k+1} -
			 \wh{x}_{n}(x_k)\|} \right\|^2.
\end{align*}
  Since  $S_k=\big\lceil { \max_i\nu_i^2  C_{i,b}^2  \over \eta^{2(k+1)}}\big\rceil \geq {\nu_i^2  C_{i,b}^2 \over \eta^{2(k+1)}} $, by  Lemma \ref{lem4},   $  \mathbb{E}[ \|x_{i,k+1}-\wh{x}_{i}(y_k) \|^2]\leq  \eta^{2(k+1)} $. Then
by taking unconditional expectations on both sides  of  the  above inequality, and using the H\"{o}lder's inequality
$\mathbb{E}[\|XY\|] \leq   ( \mathbb{E}[\|X \|^2] )^{1\over 2} ( \mathbb{E}[\|Y \|^2] )^{1\over 2}
$,  we obtain that
\begin{align}
 & \mathbb{E}[\|x_{k+1}-x^*\|^2] \leq a^2  \mathbb{E}[\|x_{k }-x^*\|^2] + 2a  \sqrt{n}\eta^{k+1} \sqrt{ \mathbb{E}[\|x_{k }-x^*\|^2] }+n\eta^{2(k+1)} \notag
  \\& =\left(a \sqrt{ \mathbb{E}[\|x_{k }-x^*\|^2] }+\sqrt{n}\eta^{ k+1}\right)^2 \notag
\\& \Rightarrow  \sqrt{ \mathbb{E}[\|x_{k+1 }-x^*\|^2] } \leq  a \sqrt{ \mathbb{E}[\|x_{k }-x^*\|^2] }+\sqrt{n}\eta^{ k+1}   . \label{bd-mse}
\end{align}
 Based on the recursion \eqref{bd-mse},  by using Lemma \ref{lem-recur}  we obtain the results.  \hfill $\Box$

 \subsection{Proof of Theorem \ref{thm4}}\label{ap:thm4}
We  first validate  the case  when $\eta=a.$   Note that
 \[  \tilde{a} ^{ k } \leq
\tfrac{\sqrt{\epsilon}}{(\sqrt{C }+ \sqrt{n} D) }  ,\quad \forall  k \geq K_{b_0}(\epsilon) =
 	\tfrac{\ln \left ((\sqrt{C }+ \sqrt{n} D) / \sqrt{\epsilon} \right)}{\ln(1/  \tilde{a} )}.\]
Then by Proposition \ref{prp-linear}(ii),  we obtain that $\mathbb{E} [\| x_k-x^*\|^2 ] \leq \epsilon $
for any $ k \geq K_{b_0}(\epsilon) .$
Then the bound given by \eqref{br-rate-prox} for $\eta=a$ holds.
By using  \eqref{exponent-sum} and \eqref{rel-exp}, we may bound the {oracle evaluations required  to obtain an $\epsilon$-NE by} \begin{align*}
\sum_{k=0}^{K_{b_0}(\epsilon) -1 } S_k& \leq   \sum_{k=0}^{K_{b_0}(\epsilon) -1 } \left(\tfrac{C_{ns}}{ a^{2(k+1)}}+1\right)
 \leq \tfrac{C_{ns}}{a^2 \ln(1/a^{2 })} a^{-2K_{b_0}(\epsilon)}+K_{b_0}(\epsilon)
\\&=\tfrac{C_{ns}}{a^2 \ln(1/a^{2 })} \left ((\sqrt{C }+ \sqrt{n} D) / \sqrt{\epsilon} \right)^{\ln (1/a^2) \over\ln(1/  \tilde{a} ) }+K_{b_0}(\epsilon).
\end{align*}
This is   the bound given in \eqref{br-rate-SFO} for the case  $  \eta=a$.
We now prove the results for  $\eta \neq a.$
 From  Proposition \ref{prp-linear}(i) it follows that
$ \mathbb{E}[ \|x_{k }-x^*\|^2]\leq \epsilon$ for $k\geq K_{b_1}(\epsilon)
 \triangleq \tfrac{\ln\left({ \left( \sqrt{C}+ \tfrac{\sqrt{n}}{\max\{a/\eta,\eta/a\}-1} \right) /\sqrt{\epsilon}}\right)}{\ln \left({1/ \max \{a,\eta\}} \right)}.$
Then we achieve the bound given in \eqref{br-rate-prox} for $\eta>a$ and $\eta<a.$
 {Similarly,} by  \eqref{exponent-sum} and \eqref{rel-exp}, we may bound the number of   sampled gradients    by
\begin{align*}
\sum_{k=0}^{K_{b_1}(\epsilon) -1 } S_k & \leq \tfrac{C_{ns} \eta^{-2K_{b_1}(\epsilon)}}{ \eta^2 \ln(1/\eta^{2 })}+K_{b_1}(\epsilon)
\\& =\tfrac{C_{ns}}{\eta^2 \ln(1/\eta^{2 })} \big(  \big( \sqrt{C}+ \tfrac{\sqrt{n}}{\max\{a/\eta,\eta/a\}-1} \big)  / \sqrt{\epsilon} \big)^{\tfrac{\ln (1/\eta^2)}{\ln \left({1/ \max \{a,\eta\}} \right)}}+K_{b_1}(\epsilon),
 \end{align*}
 giving us the required result.  \hfill $\Box$

 \subsection{Proof of Corollary \ref{poly-lem2}}\label{ap:poly-lem2}   Define
 $C_{mc}=\sqrt{\sum_{i=1}^n\nu_i^2   C_{i,b}^2}
 {\rm~with~}C_{i,b} \triangleq  {  \mu \over \mu^2+L_i^2}\Big(1-L_{fi}/ \sqrt{   \mu^2+L_{fi}^2} \Big)^{-1}.$
  Similar to Eqn. \eqref{bd-mse}, by using  $S_k=\lceil (k+1)^v \rceil$  and Lemma \ref{lem4},
\begin{align*}
\sqrt{ \mathbb{E}[\|x_{k+1 }-x^*\|^2] } &\leq  a \sqrt{ \mathbb{E}[\|x_{k }-x^*\|^2] }+  C_{mc}(k+1)^{-v/2}
\\& \leq  a^{k+1} \sqrt{ \mathbb{E}[\|x_{0 }-x^*\|^2] }+C_{mc}\sum_{m=1}^{k+1} a^{k+1-m} m^{-v/2}.
\end{align*}
Then by  \eqref{lem3-poly-bd},  we have   that for any $k\geq 1$:
 $\sqrt{ \mathbb{E}[\|x_{k }-x^*\|^2] }    \leq a^{k} \left(\sqrt{C}+C_{mc} \tfrac{e^v a^{-1}-1}{1-a }\right)+ \tfrac{ C_{mc} 2a^{-1}}{\ln(1/a)} k^{-v/2}.$
  Since $a\in (0,1)$ and $v>0$, by applying Lemma~\ref{bd-cqv}  with $u=1$,
   $a^k \leq e^{-v/2}\left(\tfrac{v}{\ln(1/a)}\right)^{v/2}  k^{-v/2}$ for $k\geq 1.$
   Then  $\sqrt{ \mathbb{E}[\|x_{k }-x^*\|^2] }  =\mathcal{O}\left( \ e^{v/2}v^{v/2} k^{-v/2} \right)$ and hence
$ \mathbb{E}[ \|x_k-x^*\|^2  ] =\mathcal{O}\left( {e^v v^v} k^{- v}\right)$. Similar to  Cor.~\ref{cor2},
we  may show that the iteration  and oracle complexity bounds  to obtain an $\epsilon-$NE are $\mathcal{O} ({v} (1/\epsilon)^{1/v})$ and    $\mathcal{O} \left( e^v v^v (1/\epsilon)^{1+1/v}  \right)$,  respectively. \hfill $\Box$

  \subsection{Proof of Lemma \ref{Lip-Tmap}}\label{ap:Lip-Tmap}

By the definition  $g_i(x_i,y)\triangleq \nabla_{x_i}f_i(x_i,x_i+ y )$ and appying the optimality condition on \eqref{def-map-T},  $T_i(y_i,z )   $ is a fixed point of the map ${\rm prox}_{\alpha r_i} \big[ x_i-\alpha\big( g_i(x_i,z -y_i )+ \mu(x_i-y_i)\big)\big]$ for any $\alpha>0$.  Then by  the triangle inequality  and the nonexpansivity property, we have the following for any $z_1,z_2 \in \mathbb{R}^d$:
\begin{align*}
 & \quad \|T_i(y_i,z_1 )   -T_i(y_i,z_2 )  \|
 \leq \Big\|    T_i(y_i,z_1 )  -\alpha\big( g_i\left(T_i(y_i,z_1 )  ,z_1    -y_i  \right)+ \mu(T_i(y_i,z_1 )  -y_i)\big)  \notag
  \\& \quad -\left( T_i(y_i,z_2 ) -\alpha\big( g_i\left(T_i(y_i,z_2 ) ,z_2   -y_i \right)+ \mu(T_i(y_i,z_2 ) -y_i)\big)\right)\Big \| \notag
\\&=\Big\| (1-\alpha \mu)\left(T_i(y_i,z_1 )   -T_i(y_i,z_2 )\right)
  -\alpha\big(g_i\left(T_i(y_i,z_1 )  ,z_1  -y_i\right) - g_i\left(T_i(y_i,z_2 ) ,z_1   -y_i \right)  \big)  \notag
\\&\quad  -\alpha\Big ( g_i\left(T_i(y_i,z_2 ) ,z_1   -y_i\right) - g_i\left(T_i(y_i,z_2 ) ,z_2    -y_i \right)  \Big) \Big \|  \notag
\\& \leq \Big\| (1-\alpha \mu)\left(T_i(y_i,z_1 )   -T_i(y_i,z_2 )\right)
  -\alpha\big(g_i\left(T_i(y_i,z_1 )  ,z_1  -y_i\right) - g_i\left(T_i(y_i,z_2 ) ,z_1   -y_i \right)  \big) \Big\|
\\&  \quad +\alpha \Big \|g_i\left(T_i(y_i,z_2 ) ,z_1   -y_i\right) -g_i\left(T_i(y_i,z_2 ) ,z_2    -y_i \right)    \Big \|.
 \end{align*}
By   Assumptions \ref{ass-agg-noise}(ii)   and   \ref{ass-agg-noise}(iii), we obtain that for any $\alpha <1/\mu: $
 \begin{equation*}
 \begin{split}
  & \Big\| (1-\alpha \mu)\left(T_i(y_i,z_1 )   -T_i(y_i,z_2 )\right)  -\alpha\big( g_i\left(T_i(y_i,z_1 )  , z_1   -y_i\right)
- g_i\left(T_i(y_i,z_2 ) , z_1     -y_i \right)  \big)  \Big\|^2
\\ &  \leq (1-\alpha \mu)^2\big\| T_i(y_i,z_1 )   -T_i(y_i,z_2 ) \big  \|^2
\\& \quad +\alpha^2 \big\|g_i\left(T_i(y_i,z_1 )  , z_1   -y_i\right)
- g_i\left(T_i(y_i,z_2 ) , z_1     -y_i \right)  \big\|^2
\\ &  \quad -2\alpha (1-\alpha \mu)\left(T_i(y_i,z_1 )   -T_i(y_i,z_2 )\right)^T   \big(g_i\left(T_i(y_i,z_1 )  , z_1   -y_i\right)
- g_i\left(T_i(y_i,z_2 ) , z_1     -y_i \right)  \big)
\\& \leq   \left(   1-2\alpha \mu + \alpha^2 (\mu ^2+L_a ^2)\right)\big\|  T_i(y_i,z_1 )   -T_i(y_i,z_2)  \big\|^2  .
\end{split}
\end{equation*}
 By combining the above two inequalities and using  Assumption \ref{ass-agg-noise}(iv), we obtain that
   \[  \|T_i(y_i,z_1 )   -T_i(y_i,z_2 )  \|  \leq
 \left(   1-2\alpha\mu  + \alpha^2 (\mu ^2+L_a ^2)\right)^{1\over 2} \|T_i(y_i,z_1 )   -T_i(y_i,z_2 )  \| + \alpha  L_{gi} \| z_1-z_2\|.\]
  Then by setting  $\alpha={\mu \over  \mu^2+L_a^2},$ we obtain  \eqref{T-Lipschitz}. \hfill $\Box$

  \subsection{Proof of Proposition \ref{prp3}}\label{app:prp3}
By noting that  $x_i^* = \wh{x}_i (x^*)$,  using    the triangle inequality   we obtain that
\begin{align}\label{prp3-bd}
\|x_{i,k+1}- x_i^*\| 	& \leq  \|x_{i,k+1} -T_i(x_{i,k}, n\hat{v}_{i,k})   \| \notag
  \\& + \| T_i(x_{i,k}, n\hat{v}_{i,k})  - \wh{x}_i(x_k)\|	+\| \wh{x}_i(x_k)- \wh{x}_i(x^*  )\| ,
\end{align}
where $\wh{x}_i(\bullet)$ is defined by \eqref{prox-BRi}. By  the definition of  $T_i(\cdot,\cdot)$ in \eqref{def-map-T} and   $y_k=  \sum_{i=1}^n x_{i,k}/n$ from  Lemma \ref{lem-pre}(ii), we have   {that for   any $i\in \mathcal{N}$, $\wh{x}_i(x_k)= T_i\left(x_{i,k}, \sum_{i=1}^n x_{i,k}\right)=T_i(x_{i,k},n y_k),$ and hence
$ \| T_i(x_{i,k}, n\hat{v}_{i,k})  - \wh{x}_i (x_k)\| \leq   n L_t\| \hat{v}_{i,k} -y_k\| $ by  \eqref{T-Lipschitz}.}
Then by using \eqref{prp3-bd},   {there holds $\|x_{i,k+1}- x_i^*\|   \leq  \|\varepsilon_{i,k+1} \| +  n L_t\| \hat{v}_{i,k} -y_k\| 	+\| \wh{x}_{i}(x_k) - \wh{x}_i (x^*)\| .$}   Then by   the triangle inequality,    using \eqref{cont-prox-best-resp}  and  $a  \triangleq \| \Gamma\|<1$,
\begin{align*}
&v_{k+1} \triangleq  \left\|   \pmat{ \|x_{1,k+1} -{x}_{1}^*\| \\
				\vdots \\
 \|x_{n,k+1}  - {x}_{n}^*\|} \right \|    \leq
av_k + \left\|  \pmat{
	\|\varepsilon_{1,k+1} \|\\
				\vdots \\
	\|\varepsilon_{n,k+1} \|} \right\|+ n L_t  \left\|  \pmat{
	\| \hat{v}_{1,k} -y_k\| \\
				\vdots \\
	\| \hat{v}_{n,k} -y_k\| } \right\|.
 \end{align*}
 Therefore, the following holds:
\begin{equation}\label{ineaxt-contractive}
\begin{split}
& v_{k+1}^2  \leq a^2 v_k^2+ \sum_{i=1}^n \|\varepsilon_{i,k+1} \|^2+ n^2L_t^2 \sum_{i=1}^n \| \hat{v}_{i,k} -y_k\|^2
+ 2av_k \left\|  \pmat{
	\|\varepsilon_{1,k+1} \|\\
				\vdots \\
	\|\varepsilon_{n,k+1} \|} \right\|\\&+ 2an L_tv_k  \left\|  \pmat{
	\| \hat{v}_{1,k} -y_k\| \\
				\vdots \\
	\| \hat{v}_{n,k} -y_k\| } \right\|+2n L_t  \left\|  \pmat{
	\|\varepsilon_{1,k+1} \|\\
				\vdots \\
	\|\varepsilon_{n,k+1} \|} \right\|  \left\|  \pmat{
	\| \hat{v}_{1,k} -y_k\| \\
				\vdots \\
	\| \hat{v}_{n,k} -y_k\| } \right\|
\end{split}
\end{equation}
From   \eqref{bd-consensus-err} and $\beta\in (0,1)$ it follows that  {for any $k\geq 0$,}
\begin{align*} \left\|  \pmat{
	\| \hat{v}_{1,k} -y_k\| \\
				\vdots \\
	\| \hat{v}_{n,k} -y_k\| } \right\| \leq  \sqrt{n} \left(  C_1\beta^{(k+1)(k+2)/2} + C_2  \beta^{ k+1 } \right)  \leq
\sqrt{n} \left(  C_1 + C_2 \right)    \beta^{ k+1 } .
\end{align*}
Note by $S_k=\left\lceil { C_r^2  \max_i\nu_i^2\over  \eta^{2(k+1)}}\right\rceil$ and Lemma \ref{agg-variance} that  $\mathbb{E}[\| \varepsilon_{i,k+1} \|^2] \leq \eta^{2(k+1)} $ for any $ i \in \mathcal{N}$.
 Then
by taking expectations   of   {the  inequality \eqref{ineaxt-contractive}}, using the H\"{o}lder's inequality
$\mathbb{E}[\|XY\|] \leq   ( \mathbb{E}[\|X \|^2] )^{1\over 2} ( \mathbb{E}[\|Y \|^2] )^{1\over 2}
$,  we obtain that
\begin{equation*}
\begin{split}
 & \mathbb{E}[\|x_{k+1}-x^*\|^2] \leq a^2  \mathbb{E}[\|x_{k }-x^*\|^2] +n\eta^{2(k+1)}+ n^3 L_t^2  \left(  C_1+ C_2 \right)^2\beta^{ 2(k+1)  } \\& + 2a  \sqrt{n}\eta^{k+1} \sqrt{ \mathbb{E}[\|x_{k }-x^*\|^2] }+ 2a n^{3\over 2} L_t    \left(  C_1  + C_2 \right)    \beta^{ k+1}\sqrt{ \mathbb{E}[\|x_{k }-x^*\|^2] }
 \\&+2n^2 L_t  \eta^{ k+1}\left(  C_1\beta^{1/2} + C_2 \right)    \beta^{ (k+1) /2}
\\& =\left(a \sqrt{ \mathbb{E}[\|x_{k }-x^*\|^2] }+\sqrt{n}\eta^{ k+1} + n^{3\over 2} L_t\left(  C_1  + C_2 \right)    \beta^{  k+1 }\right)^2 \\&
\Rightarrow \sqrt{ \mathbb{E}[\|x_{k+1 }-x^*\|^2] } \leq  a \sqrt{ \mathbb{E}[\|x_{k }-x^*\|^2] }+\sqrt{n}\eta^{ k+1} + n^{3\over 2} L_t \left(  C_1 + C_2 \right)    \beta^{  k+1 }.
\end{split}
\end{equation*}
Thus,  by the definitions  $\gamma=\max\{\eta,  \beta\}$ and $C_4 = \sqrt{n}  + n^{3\over 2} L_t \left(  C_1  + C_2 \right) $,
 {we obtain that for any $k\geq 0,$  $ \sqrt{ \mathbb{E}[\|x_{k+1 }-x^*\|^2] } \leq  a \sqrt{ \mathbb{E}[\|x_{k }-x^*\|^2] }+C_4\gamma^k
$.  Based on which, by  using Lemma \ref{lem-recur} we obtain the results.}
\hfill $\Box$
\def\cprime{$'$}

\end{document}